\documentclass[leqno]{amsart}

\title[Probabilistic Representation for Nonlinear Dirac Equations]{Probabilistic representation and Monte Carlo method for nonlinear Dirac equations}

\author{Jean-Pierre Fouque}
\address{Department of Statistics and Applied Probability, South Hall, University of California, Santa Barbara, CA 93106, USA}
\email{fouque@pstat.ucsb.edu}
\thanks{Email for Fouque: \href{mailto:fouque@pstat.ucsb.edu}{fouque@pstat.ucsb.edu}}

\author{Tomoyuki Ichiba}
\address{Department of Statistics and Applied Probability, South Hall, University of California, Santa Barbara, CA 93106, USA}
\email{ichiba@pstat.ucsb.edu}
\thanks{Email for Ichiba: \href{mailto:ichiba@pstat.ucsb.edu}{ichiba@pstat.ucsb.edu}}

\author{Ka Lok Lam}
\address{Department of Statistics and Applied Probability, South Hall, University of California, Santa Barbara, CA 93106, USA}
\email{kalok@ucsb.edu}
\thanks{Email for Lam: \href{mailto:kalok@ucsb.edu}{kalok@ucsb.edu}}

\usepackage{paper}
\RequirePackage[authoryear, square]{natbib}
\renewcommand{\cite}{\citep}

\usepackage{graphicx} 

\begin{document}
\begin{abstract}
 Dirac equations describe the relativistic evolution of quantum particles. Using the point-process stochastic calculus of piecewise-deterministic Markov processes (PDMPs) on general state spaces, which naturally allows additional Feynman-Kac potentials, we develop a probabilistic representation for the nonlinear Dirac equations in $1+1$ spacetime dimension with cubic nonlinearity and $L_\infty$-type problem data using a multi-type branching telegraph process together with both continuous and jump-type multiplicative functionals. Our main theorem leads to a classification result which, when applied to the specific Soler model, gives a sharp finite-time admissibility of the probabilistic representation. By extending the one-step representation to a multi-step, time-marching representation using an analytically-derived a priori $L_\infty$ estimate, we observe a ``local probabilistic representation, but global PDE well-posedness" phenomenon. We proceed to propose a time-marching Monte Carlo algorithm as a global numerical solver for the Soler system with a rigorous convergence proof and establish an additional stability result on the PDE side. Finally, while the underlying process parameters do not affect the validity of the probabilistic representation, we observe their potential for variance reduction in the Monte Carlo methods. 
\end{abstract}

\maketitle
\section{Introduction and Literature Review}\label{sec: introduction}
In 1928, almost a century ago, Paul Dirac proposed the equation that bears his name in an attempt to combine quantum mechanics with special relativity to accurately model the dynamics of electrons, ultimately predicting the existence of anti-particles (positrons in the case of electrons) and laying the foundation of quantum field theory to incorporate creation and annihilation of particles when describing dynamics in its natural state \cite{Dirac1928}. 

To model more accurately the effect of interaction between particles and the external field, researchers,  including \cite{Thirring1958, Soler1970, GrossNeveu1974} introduced self-interaction terms to the evolutionary dynamics and started to look at nonlinear Dirac equations. For an introduction to the Dirac equations from the perspective of mathematical Physics, we refer readers to the monograph \cite{Thaller1992}. Meanwhile, for the nonlinear case, we  refer readers to the survey \cite{CuevasMaraverEtAl2018} and the recent book \cite{Yang2023} with references therein for a mathematical introduction from the PDE perspective. For an account on its numerical analysis, we refer readers to the survey \cite{XuShaoTang2013} and the more recent paper \cite{CaiWang2022} and references therein. 

The association of paths with equations of motions from quantum mechanics originated from Richard Feynman's idea of path integrals \cite{Feynman1948}. In the famous Problem 2.6 of the classical textbook \cite{FeynmanHibbs1965}, it was already suggested that the Dirac equations are related to paths with ``zigzags", which are precisely how the trajectories of a telegraph process look like,  together with their number of corners. Several rigorous mathematical accounts then came along in the 1980s around the same time to connect the complex-valued solutions to linear Dirac equations with the telegraph process as a genuine stochastic process (as a probabilistic measure on paths instead of a complex measure as in the path integral formulation). These include several Markov process theory approaches \cite{Gaveau1984}, \cite{BCSS86}  and  a functional analytic approach \cite{Ichinose1984}, who also showed the existence of a path measure from the Dirac operator using a duality argument in addition to having the process constructed a priori as in the Markov theory approach. Meanwhile, at the PDE level, both approaches connected the probabilistic representation with  a priori $L_\infty$ estimates to the equations.  We refer readers to the paper \cite{BCSS1990} and references therein for further development at that time. We remark that, independent from the previous approaches, there is also a very recent work by the current authors to establish similar probabilistic representations using an asymptotic analysis of quantum walks, primarily motivated by an exploration to quantum computing. We refer readers to \cite{FIL2025} and references therein.

The current work extends the aforementioned legacy of  finding probabilistic representations for the Dirac equations to the nonlinear case with cubic nonlinearity, which will be introduced in \eqref{eqns: sect5_general_cubic_nonlinear_Dirac}. Our main representation theorem (Theorem \ref{thm: general_cubic_main_theorem}) leads to a classification result on the admissibility time of such representation based on explicit bounds on the PDE coefficients (Corollary \ref{corollary: classification_of_stochastic_representation}). In the case of the Soler model, the probabilistic representation has a sharp finite-time admissibility (Theorem \ref{thm: sec4_main_theorem}, Corollary \ref{corollary: sharpeness results_constant_data}). Our approach uses the point-process stochastic calculus of piecewise deterministic Markov processes (PDMPs) pioneered by Mark H.A. Davis \cite{davis1984}, which we review in Section \ref{sec. PDMP} for the case of general state spaces. More precisely,  we formulate a multi-type branching telegraph process, which we term the \textit{Dirac branching system} (introduced in Definition \ref{defn: Dirac branching system}) as a measure-valued PDMP, together with a lifted PDMP which contains additional information on its jump events.  We then study the stochastic differential equations (SDEs) satisfied by its associated multiplicative functionals and multiplicative test functions. 

The theory of probabilistic representation for evolutionary equations has a long history. We refer readers to \cite{LPS1998} (with the English translation \cite{LPS_Eng03}) and the short survey \cite{PierreNicolas10} for a general account. Specifically for transport equations, to which the Dirac equations belong, we refer readers to \cite{SW2005} for a treatment of the Boltzmann equation (which involves particle collisions), and the recent work of \cite{CoxHarrisKyprianouWang2022} together with the book \cite{HortonKyprianou2023} and references therein for neutron transport equations, whose simulation gave birth to the Monte Carlo algorithm \cite{MetropolisUlam1949}.  For diffusive equations, we refer readers to the very recent book \cite{Pages2026} and references therein from the point of view in financial mathematics, which is traditionally based on Brownian randomness \cite{Bachelier1900}, as opposed to jump-type randomness in PDMPs. 

Branching structures come into play because they can represent nonlinear terms in an evolutionary equation. Seminal examples include the short paper \cite{McKean1975} for the  Kolmogorov-Petrovskii-Piskunov (KPP) equation on its asymptotic speed, and the work \cite{LeJanSznitman1997} for the  Navier-Stokes equation on its well-posedness with respect to $L_\infty$ conditions in the Fourier variable. We also mention and refer  readers to references therein, \cite{Dunbar1988} and \cite{Ratanov2006} that extend \cite{McKean1975} to the case of nonlinear hyperbolic equations, making use of an underlying stochastic process in the spirit of our Dirac branching system. For a numerical account, we refer readers to the recent work of \cite{HenryLabordereOudjaneTanTouziWarin2019} and \cite{HT21} together with the recent combination with deep learning \cite{NYGN2024} and references therein. 

While the branching Monte Carlo method of \cite{HT21} includes semilinear hyperbolic PDEs, to which our nonlinear Dirac equation belongs, we remark the difference in approaches of our work with theirs, which makes use of the Duhamel principle and a renewal argument, producing a branching process without intermediate movement but additional non-local branching rule (that is, offspring  are born with an additional shift from their parents' position) determined by the Green function of the differential operator. Meanwhile, describing our Dirac branching system as a PDMP emphasizes its piecing-out mechanism to naturally accommodate a branching process with intermediate movements and local branching rules, which are explicitly described by its associated SDEs. This is particularly interesting in the presence of a potential term in the equation which is handled by a multiplicative exponential of an integral of the potential along the trajectory of the branching process in the spirit of the Feynman-Kac formula \cite{Kac1949}, which can be described without full knowledge of the jump events. That will be the case for the equations considered in Section \ref{sec: state_dependent_coefficients}.

With the probabilistic representation at hand, our second contribution is then the promotion of it as a tool to study the analytic property of the PDE solution with $L_\infty$-type data and as a numerical tool using Monte Carlo. At the PDE level, the finite-time admissibility arising in the Soler case (Theorem \ref{thm: sec4_main_theorem}) is consistent with the literature of using branching (which is often termed stochastic cascades in the PDE community) to describe nonlinearity, providing mostly a local-in-time probabilistic representation without additional assumption; see, for instance, the aforementioned Navier-Stokes work \cite[Remark 2.5]{LeJanSznitman1997} and the branching Monte Carlo paper \cite[Theorem 3.12]{HenryLabordereOudjaneTanTouziWarin2019}, which termed the phenomenon the ``small maturity" issue. We also mention \cite{BM2010} which works in a similar setting to ours and uses the renewal argument as in  \cite{HT21}, as we pointed out in the previous paragraph, to study the local well-posedness of nonlinear hyperbolic PDEs. Here again, our approach differs by using a point-process stochastic calculus, which is accompanied by a martingale argument.

Interestingly for the Soler model, we also observe the same phenomenon described in \cite{BlomkerRomitoTribe2007} and \cite{Morandin2005} that the integrability of the stochastic process, and hence the probabilistic representation is only valid locally, but the PDE is still globally well-posed. We address this phenomenon by extending the one-step probabilistic representation to a multi-step time marching representation in Theorem \ref{thm: global_well_posednesss_L_inf} via an $L_\infty$ a priori estimate modified from that of \cite{Huh2013GrossNeveu}, which can be proven in an analytic flavor. We then propose at the level of numerical method a time-marching Monte Carlo algorithm which is described in Algorithm \ref{alg:global-monte-carlo-time_marching} as a global numerical solver for the Soler model and  rigorously verify its convergence in Theorem \ref{thm: convergence_of_global_alg}. Lastly, we demonstrate preliminary numerical experiments and remark that while the probabilistic representation admits a degree of freedom of the Dirac branching parameters of the underlying Dirac branching system, such choice can have a variance reduction effect on our Monte Carlo implementation.

The remainder of the paper is organized as follows. In Section \ref{sec. PDMP}, we review the necessary preliminaries for PDMPs on general state spaces. From them, we introduce in Section \ref{sec: nonlinear_case} the Dirac branching system, discuss its stochastic calculus and semigroup properties, and establish our first main result, a probabilistic representation for the Soler model. In Section \ref{sec: state_dependent_coefficients}, we extend this probabilistic representation to  nonlinear Dirac equations which include more general Dirac operators and more general cubic nonlinearities as well as potential terms. With applications of the probabilistic representation in mind, in Section \ref{sec: global_section}, we first develop a one-step Monte Carlo scheme for general nonlinear Dirac equations, and then establish for the Soler model the corresponding global well-posedness result and global-in-time Monte Carlo scheme via a time-marching principle. We further derive additional stability results, establish rigorous convergence analysis for the global algorithm and accompany these results with numerical experiments. In Section \ref{sec. conclusion}, we summarize the paper and discuss future directions. Finally, Appendices \ref{appendix: main_proofs_inSec3} to \ref{appendix. Global arguments} collect proofs and technical details deferred from Sections \ref{sec: nonlinear_case} to \ref{sec: global_section}. 
 
\section{Preliminaries for PDMPs of general state spaces}\label{sec. PDMP}
Piecewise deterministic Markov processes (PDMPs) were introduced by Davis in \cite{davis1984} who also wrote the textbook \cite{davis93}, focusing on the case of finite-dimensional state spaces. In the following, we extract preliminary results from \cite{JS96_JMP} and \cite{Jacobsen2007} who dealt with the case of general state spaces so that the space of measures, which are in general infinite-dimensional, is included. We remark that PDMPs are used significantly in Mathematical biology and we point the reader to the survey \cite{Cloez2018} and references therein from that perspective.  Readers interested in the main results may proceed directly to Section \ref{sec: nonlinear_case} and refer back to this section as needed.

\subsection{Notations and definitions}
Let $E$ be a Blackwell measurable space (which includes all topological subspaces of a Polish space with its Borel sigma algebra). A PDMP over $E$ can be described  by a flow map $\phi: [0, \infty)\times E\to E$, a rate function $\ell:[0, \infty)\times E\to [0, \infty)$ and a probability kernel $\Gamma = \Gamma(x, A)\in [0, 1]$ for all $x\in E$ and measurable sets $A$.  The triple $(\phi, \ell, \Gamma)$ is an admissible local characteristic if:
\begin{enumerate}
    \item the flow function $\phi$ is a measurable function from $[0, \infty)\times E$ which satisfies:
\begin{align}
    \phi(0, x) &= x, \\
    \phi(t+s, x) &= \phi(s, \phi(t, x)),\quad s, t\geq 0, \quad x\in E.
\end{align}
    \item the rate function $\ell$ is measurable from $[0, \infty)\times E$ such that $\ell(0, x) = 0$ and  $t\mapsto \ell(t, x)$ is continuously increasing for every $x\in E$; in addition, we have the additivity property:
\begin{align}
    \ell(t+s, x) = \ell(t, x) + \ell(s, \phi(t, x)),\quad s,t\geq 0, \quad x\in E.
\end{align}
    \item the probability kernel $\Gamma$ has zero diagonal, that is, 
\begin{align}
    \Gamma(x, \{x\}) = 0,\quad x\in E.
\end{align}
\end{enumerate}
With an admissible local characteristic $(\phi, \ell, \Gamma)$, one can define the state-dependent survival functions $H(x,[t, \infty)) := e^{-\ell(t, x)}$ for all $x\in E$ and $t\geq 0$. If $\lim_{t\to \infty}\ell(t, x) =\infty$ for all $x\in E$, then the inter-arrival time associated with the survival function $H(x, \cdot)$ is almost surely finite, and we can define the probability kernel $G$ on the time-augmented state space $[0, \infty)\times E$ by:
\begin{align}\label{eqns: PDMP_embedded_Markov_chain_kernel}
    G(s, x; dt, dy) = H(x, dt) \Gamma(\phi (t, x), dy), \quad x\in E,\quad s,t\geq 0,
\end{align}
where we remark that the transition kernel $G$ does not depend on the current time $s$. If we fix a sufficiently rich probability space $(\Omega, \F, \P)$, then the transition kernel $G$ defines a Markov chain $(\sigma_n, Y_n)_{n=0}^\infty$ on it over $[0, \infty)\times E$ with:
\begin{align}
    \P(\sigma_{n+1}\in dt, Y_{n+1}\in dy \mid \sigma_n = s, Y_n = x) = G(s, x; dt, dy), \quad n\geq 0. 
\end{align}
Define $\tau_0 = 0, \tau_n:= \sigma_1+\cdots + \sigma_n$ for all $n\geq 1$, we can define the continuous-time process $(X_t)$ on $E$ by:
\begin{align}\label{eqns: PDMP_definition_process}
    X_t = \sum_{n=0}^\infty \phi (t - \tau_n,Y_n) \I_{\tau_n \leq t<\tau_{n+1}},
\end{align}
together with the marked point process $N$ over $[0, \infty)\times E$ defined by:
\begin{align}\label{eqns: MPP_definition_process}
    N(dt, dx) = \sum_{n=1}^\infty \delta _{(\tau_n, Y_{n})} (dt, dx),
\end{align}
and call
\begin{align}\label{eqns: embedded Markov chain}
    (\tau_n, Y_n)_{n=1}^\infty ,
\end{align}
the embedded Markov chain. 
We remark that for every $n\geq 1$, $Y_n$ is the \textit{post-jump} state at the time $\tau_n$. In particular every trajectory  of $(X_t)_{t\geq 0}$ defined in \eqref{eqns: PDMP_definition_process} starts at $X_0:= Y_0$. Then for every $k\geq 0$, it flows deterministically along $\phi$ over $t\in [\tau_{k}, \tau_{k+1})$ with the state $\phi(t - \tau_k, Y_{k})$ and then jumps to $Y_{k+1}$, continuing in the next time interval. 

The construction and stochastic calculus below apply on all of $[0, \infty)$ when the process is regular, that is, $\lim_n \tau_n = \infty$ almost surely, which is our case by Lemma \ref{lemma: finite_jump_of_Dbs}. 

With all these setups, here comes the foundational result that identifies the randomness behind PDMPs as marked point processes:
\begin{theorem}\label{thm: fundamental theorem of PDMP}
    The process $(X_t)$ is a time-homogeneous, strong Markov process on $E$ whose natural filtration coincides with the filtration generated by the marked point process $N$, whose compensator is given by:
\begin{align}\label{eqns: definition_of_nu}
    \nu(dt, dx) = \sum_{n=0}^\infty \frac{G(\tau_n, Y_n;  dt - \tau_n, dx)}{H(Y_n,  [t-\tau_n, \infty))} \I_{(\tau_n, \tau_{n+1}]}(t).
\end{align}
In particular, if $\ell $ is differentiable in $t$ then:
\begin{align}
    \nu(dt, dx) = \frac{d}{ds}\ell(s, X_{t-})\mid_{s= 0}\Gamma(X_{t-}, dx)dt.
\end{align}
\end{theorem}
\begin{proof}
    See Lemma 1, Theorem 8, the un-numbered Equation in Section 4.1 and Equation 44 of \cite{JS96_JMP}. 
\end{proof}

\begin{definition}[PDMP] \label{defn: PDMP}
Let $(\phi, \ell, \Gamma)$ be an admissible triple of local characteristics. We call $(X_t)$ defined in \eqref{eqns: PDMP_definition_process} the PDMP with local characteristics $(\phi, \ell, \Gamma)$, and $N$ defined in \eqref{eqns: MPP_definition_process} its associated marked point process. 
\end{definition}

\subsection{Stochastic calculus of PDMPs}
A central step in the proof of our probabilistic representation is to identify certain local martingales to be true martingales. The machinery is the stochastic calculus of PDMPs. In the following, we review semi-martingales associated to a PDMP together with the SDEs that they satisfy almost surely, that is, in the strong sense. In the following, we fix a PDMP $(X_t)$ with local characteristics $(\phi, \ell, \Gamma)$ together with its associated marked point process $N$. 
\subsubsection{The It\^o formula}
The first result is an It\^o formula, that is a chain rule. To this end for every measurable function $g:E\to \C$ we introduce a time-dimension to define the function $\hat g:[0, \infty)\times E\to \C$ by:
\begin{align}\label{eqns: lifting_via_flow}
    \hat g(t, x):= g(\phi(t, x)).
\end{align}
In addition, we let $\tilde N:= N - \nu$ be the compensated process of the associated marked point process $N$ defined in \eqref{eqns: MPP_definition_process}  with intensity $\nu$ \eqref{eqns: definition_of_nu}. 

What follows is the It\^o Lemma for PDMPs which first establishes when the scalar-valued process $(g(X_t))$ is a semi-martingale, and if that happens, the SDE satisfied by it. We also identify the semimartingale decomposition using the compensated process $\tilde N$ under additional  assumptions:
\begin{proposition}\label{prop: semi-martingale_of_PDMP}~
\begin{enumerate}
    \item The process $(g(X_t))$ is a semi-martingale if and only if  for every $x\in E$, the map $t\mapsto \hat g(t, x)$ is c\`{a}d\`{a}g and of finite variation locally over $[0, \infty)$. 
    \item If in addition, $t\mapsto \hat g(t, x)$ is Lipschitz over compact sets of $[0, \infty)$ and is differentiable at $0$, then $(g(X_t))$ satisfies the stochastic differential equation almost surely:
\begin{align}\label{eqns: pdmp_ito_formula}
    dg(X_t) = \frac{\p}{\p s}\hat g (s, X_t)\mid_{s = 0} dt + \int_{E} [g (y) - g(X_{t-})]N(dt, dy);
\end{align}
    \item if in addition 
\begin{align}\label{eqns: local_martingale_PDMP_compensator}
   \int_0^t\int_{E} \abs{g (y) - g(X_{s-})}\nu (ds, dy)<\infty\quad \mbox{almost surely},
 \end{align}
 for every $t>0$, then we can decompose:
\begin{align}
  \int_{E}[g(y) - g(X_{t-})]N(dt, dy) = \int_{E} [g (y) - g(X_{t-})]\nu (dt, dy) 
  +  \int_{E} [g (y) - g(X_{t-})]\tilde N(dt, dy),
\end{align} 
where:
\begin{align}\label{eqns: compenstated_local_martingale}
   t\mapsto  \int_0^t \int_E [g(y) - g(X_{s-})]\tilde N(ds, dy)
\end{align}
is a local martingale over $[0, \infty)$. 
    \item if in addition:
\begin{align}
    \E \int_0^T \int_E \abs{g(y) - g(X_{t-})} \nu(dt, dy) <\infty
\end{align}
for some $T>0$, then \eqref{eqns: compenstated_local_martingale} is a true $L_1$ martingale over $[0, T]$.
\end{enumerate}
\end{proposition}
\begin{proof}
    See Theorem 19 (a) of \cite{JS96_JMP} for (1) and (2). For (3) and (4), See Theorem 4.6.1 (iii) of \cite{Jacobsen2007}. A concise idea of proof follows from Theorem T3 and Corollary C4 in Chapter VIII.1 of \cite{bremaud1981}. 
\end{proof}

\subsubsection{Multiplicative functionals and their stochastic differential equations}\label{sec. Multiplicative_functional}
 In this section, we consider multiplicative functionals of the PDMP $(X_t)_{t\geq 0}$ and derive the SDEs they satisfy almost surely. To this end, we denote $(\theta_t)$ the usual shift operator for the Markov process $(X_t)$ and let   $(\tau_n, Y_{n})$ be the embedded Markov chain defined in \eqref{eqns: embedded Markov chain} associated to $(X_t)$. Note that  $X_{\tau_n} = Y_n$ for every $n\geq 1$. 
 
 We call a c\`adl\`ag $\C-$valued process $(A_t)$ adapted to $(\F_t^X)$, the natural filtration of $(X_t)$ a (strong) multiplicative functional if for every stopping time $T$ and $t>0$, we have:
\begin{align}
    A_{t+T} = A_T\cdot A_t \circ \theta _T,
\end{align}
almost surely on $\{T<\infty\}$. 
\begin{definition}\label{defn: definition_of_MF_PDMP}
Let $\tilde \fu:E\to \C$ be measurable and let  $\tilde \fm:E\times E\to \C$ be measurable.
\begin{enumerate}
    \item We call the process $(E_t)$ defined by: \begin{align}
    E_t := \exp(\int_0^t \tilde \fu(X_s)ds),
\end{align}
a continuous multiplicative functional and $\tilde \fu$ its potential, whenever the integral is well-defined for $t\geq 0$. 
    \item We call the process $(M_t)$ defined by:
\begin{align}
\begin{cases}
    M_t &= M_{\tau_n},\quad t\in [\tau_n, \tau_{n+1}),\quad \\
    M_{\tau_n} &= M_{\tau_n-} \, \tilde \fm(X_{\tau_n-}, X_{\tau_n})    
\end{cases},\quad n\geq 1; \quad M_t = 1, \quad 0\leq t<\tau_1,
\end{align}
a jump multiplicative function and $\tilde \fm $ its jump potential. 
\end{enumerate}
\end{definition}
It is not difficult to check that $(E_t)$ and $(M_t)$ are actually multiplicative functionals of $(X_t)$. 

By the standard chain rule, the process $(E_t)$  satisfies the stochastic differential equation almost surely:
\begin{align}
    dE_t = \tilde \fu (X_t) E_t dt, \quad E_0 = 1. 
\end{align}
When $\hat {\tilde \fu }$ (defined in \eqref{eqns: lifting_via_flow}) is of local finite variation, it is a continuous semi-martingale.

 Meanwhile, the jump multiplicative function $(M_t)$  can be written as the integral against the marked point process:
\begin{equation} \label{eq:Mt}
    M_t =1+ \int_0^t\int_E M_{s-} (\tilde \fm(X_{s-}, y)-1) N(ds, dy).
\end{equation}
To emphasize the multiplicative structure, we can also write:
\begin{equation}\label{eq:Mt2}
    M_t = \prod_{0<s\leq t} \tilde \fm(X_{s-}, y)^{N(ds, dy)}. 
\end{equation}
Similarly to \eqref{eqns: local_martingale_PDMP_compensator} in Proposition \ref{prop: semi-martingale_of_PDMP}, if the process $(M_t)$ satisfies:
\begin{align}
    \P \Big(\int_0^t \int_E \abs{M_{s-} (\tilde \fm(X_{s-}, y) -1)}\nu(ds, dy)<\infty \Big) = 1,
\end{align}
for every $t>0$, then $(M_t)$ satisfies the following  semi-martingale decomposition almost surely:
\begin{align}
    dM_t = \int_E M_{t-} (\tilde \fm(X_{t-}, y)-1) \nu(dt, dy)+\int_E M_{t-} (\tilde \fm(X_{t-}, y)-1) \tilde N(dt, dy), \quad M_0 = 1.
\end{align}

\newpage 
\section{Probabilistic representation for the Soler model}\label{sec: nonlinear_case}
We first consider the specific case of nonlinear Dirac equations called the massive Soler model with mass $m>0$, which is given by the evolutionary equation:
\begin{equation}\label{eqns: sec4_nonlinear_Dirac}
\begin{cases}
    \p_t \psi(t, x, z) = -z\p_x \psi(t, x, z)-im \psi(t, x, -z)+i \mathcal N(\psi, \overline{\psi})\psi(t, x, -z);\\
     \psi(0, x, z) = g(x, z), \quad (t, x, z)\in [0, \infty)\times \R \times \{\pm 1\};\\
    \mathcal N(\psi, \overline{\psi}) = \psi(t, x, z)\overline{\psi(t, x, -z)}+\psi(t, x, -z)\overline{\psi(t, x, z)}, 
\end{cases}
\end{equation}
where $g$ is the boundary condition and $\overline{\psi}$ represents the complex conjugate of $\psi$. In this work, the problem setup consists of $L_\infty$-type function spaces, which is consistent with the literature of probabilistic representation from, for instance, \cite{Ichinose1984} and \cite{BM2010}. More precisely, we denote the function space 
\begin{equation}\label{eqns: cb1_defn}
    C_b^1:= \{g: \R\times \{\pm 1\}\to \C \mid g, \p_x g\in L_\infty(\R\times \{\pm 1\})\},
\end{equation}
of bounded continuous functions $g(x, z)$ differentiable with respect to $x$ with a bounded derivative, which is not necessarily continuous. We always  consider the initial condition $g$ to satisfy $g\in C_b^1$ and the solution $\psi$ to be a \textit{classical} solution over $[0, T)$ where $T\in [0, \infty]$ so that  $\psi(t, x, z)$ with $t\in [0, T),  x\in \R, z= \pm 1$  is differentiable with respect to both $t$ and $x$ such that
\begin{align}\label{eqns: classical_solution_defn}
    \psi(t, \cdot)\in C_b^1,
\end{align}
for every  $t\in [0, T)$.  Moreover, we require that for every $t\in [0, T)$:
\begin{align}\label{eqns: classical_solution_additional_boundedness_condition}
    \sup_{s\in [0, t]}\bigg( \norm{\psi(s, \cdot)}_\infty + \norm{\p_x\psi(s, \cdot)}_\infty\bigg) <\infty. 
\end{align}
We call the case for $T = \infty$ a global classical solution. 
\begin{remark}
    We will see in Theorem \ref{thm: sec4_main_theorem} that our probabilistic representation will automatically satisfy the notion of a classical solution with the condition at \eqref{eqns: classical_solution_additional_boundedness_condition}.  
\end{remark}

 Typical problem setups for nonlinear Dirac equations consider well-posedness in the fractional Sobolev spaces $H^s(\R\times \{\pm 1\})$ with $s\geq 0$ (with $H^0(\R\times \{\pm 1\}) = L_2(\R\times \{\pm 1\})$; we refer the reader to \cite{HuhMoon2015} and the survey \cite{Pelinovsky2011} together with references therein. For treatments that consider non-decaying problem data, we refer readers to \cite{BournaveasZouraris2012} that considers periodic boundary conditions and the more recent \cite{Trapasso2020} that considers modulation spaces from time frequency analysis.

This section contains our first main result of establishing a  probabilistic representation for nonlinear Dirac equations in the case of  the massive Soler model (see Theorem \ref{thm: sec4_main_theorem}).  In the following, we first construct the underlying process, which we call the \textit{Dirac branching system}, develop its stochastic calculus and finally give our probabilistic representation theorem.  The treatment of the case for general cubic nonlinearity is postponed to Section \ref{sec: state_dependent_coefficients}, yet we emphasize that  it follows from a slight modification to machineries developed in this special case for the massive Soler model. 

\subsection{The Dirac branching system}\label{subsec: dirac branching system} 
\subsubsection{The measure-valued process}

We introduce the multi-type branching telegraph process, which we term the \textit{Dirac branching system}, underlying our probabilistic representation to \eqref{eqns: sec4_nonlinear_Dirac}. The process can be described as a PDMP over the space $E:= M_f(\X)$ of finitely supported measures of the enlarged state space $\X:= \R\times \{\pm 1\}\times \{\pm 1\}$. Heuristically, the extra  component $\{\pm 1\}$ is added to the state space of $\R\times \{\pm 1\}$ to handle complex conjugation in the nonlinear term in \eqref{eqns: sec4_nonlinear_Dirac}. To be precise, $\mu\in M_f(\X)$ if and only if:
\begin{align}
    \mu = \sum_{i=1}^N \delta_{(x_i, z_i, \theta_i)}, \quad N\in \N, \quad (x_i, z_i, \theta_i)\in \X,\quad i = 1,\cdots, N.
\end{align}
Alternatively, we identify each element $\bx \in \mathcal X$ as $\bx = (x_\bx, z_\bx, \theta_\bx)$ and we also write 
\begin{equation} \label{eq:mu-def}
    \mu = \sum_{\bx\in \mu}\delta_{\bx} = \sum_{\bx\in \mu}\delta_{(x_\bx, z_\bx, \theta_\bx)}
\end{equation}
for simplicity. For every measurable function $\fh :\X\to \C$, we define 
    $\inprod{\mu , \fh}:= \sum_{\bx\in \mu} \fh (\bx)$ for every $\mu\in E$ and define accordingly the total mass  $\abs{\mu}:= \inprod{\mu, 1} = \mu(\X)$  for $\mu\in E$. Note that $\X$ is a Polish space, so $E$ is naturally equipped with the topology of weak convergence. The induced Borel $\sigma-$algebra then makes $E$ a Blackwell measurable space, so the results on PDMPs in Section \ref{sec. PDMP} are applicable. 

Now, we follow Definition \ref{defn: PDMP} to describe the local characteristics of the PDMP on $E$. To describe the transition function of the PDMP, we will need to define the jump symbols:
\begin{align}
    J:= \{\mathsf{f}, \mathsf{b_1}, \mathsf{b_2}\},
\end{align}
as well as the jump transition function $\Phi:E\times \X\times J\to E$ defined by: 
\begin{align}\label{eqns: jump transition function}
\Phi(\mu; \bx, \mathsf{f}) &= \mu - \delta_{\bx}
 + \delta_{(x_{\bx}, -z_{\bx}, \theta_{\bx})},\\
 \Phi(\mu; \bx, \mathsf{b_1}) &= \mu- \delta_{\bx}
 + \delta_{(x_{\bx}, -z_{\bx}, \theta_{\bx})}
 + \delta_{(x_{\bx}, 1, 1)}+ \delta_{(x_{\bx}, -1, -1)},\\
 \Phi(\mu; \bx, \mathsf{b_2}) &= \mu - \delta_{\bx}
 + \delta_{(x_{\bx}, -z_{\bx}, \theta_{\bx})}
 + \delta_{(x_{\bx}, -1, 1)}+ \delta_{(x_{\bx}, 1, -1)},
\end{align}
for every $\bx\in \mu$. We are then ready to define: 
\begin{definition}[Dirac branching system]\label{defn: Dirac branching system} Let $\lambda>0$, $\{p(e)\}_{e\in J}$ be a probability mass function over $J$ (so $\sum_{e\in J}p(e) = 1$ and $p(e)> 0$ for all $e\in J$). Then the Dirac branching system $(\mu_t)$ with parameters $\lambda>0$ and $\{p(e)\}_{e\in J}$ is the PDMP in Definition \ref{defn: PDMP} on $E =  M_f (\mathcal X)$ associated to the local characteristics $(\phi, \ell, \Gamma)$ defined by:
\begin{align}
    \phi(t; \mu) &= \sum_{\bx \in \mu} \delta_{(x_\bx - tz_\bx, z_\bx, \theta_\bx)}, \quad \mu \in M_f(\mathcal X), \\ \label{eqns: rate_function_of_Dbs}
    \ell(t; \mu) &= t\lambda \abs{\mu}, \\
    \label{eqns: jump_transition_kernel}
    \Gamma(\mu; d\mu') &= \sum_{e\in J, \bx \in \mu} \frac{p(e)}{\abs{\mu}} \delta_{\Phi(\mu; \bx, e)}(d\mu'),\quad  \sum_{e\in J}p(e) = 1, \quad p(e)\geq 0.
\end{align}
\end{definition}
\begin{remark}\label{remark: concrete_description_of_dbs}
For a more concrete description, the Dirac branching system can be thought of as an evolution on a family of particles described by $\bx\in \X$ with position $x_\bx$, velocity $z_\bx$ and type $\theta_\bx$.   At each event time, which is driven by an independent exponential clock per particle, our process first chooses one particle with state $\bx\in \X$, among all particles and then chooses among the $\abs{J} = 3$ actions in $J$ which have the effect described in the next paragraph. Meanwhile, between events, each particle moves in a straight line in the $-z_\bx$ direction. 

When the action is $\mathsf{f}$, which corresponds to flipping, the particle flips the sign of its z-component. Otherwise, in addition to flipping its z-component, the particle gives birth to two new particles based on the $j$th branching rule, which corresponds to $\mathsf{b}_j$ with $j = 1, 2$. When $j = 1$, the states of the offspring are $(x_\bx, 1, 1), (x_\bx, -1, -1)$; when $j = 2$, the states of the offspring are $(x_\bx, -1, 1), (x_\bx, 1, -1)$.  This concrete description underpins the sampling procedure of the Dirac branching system and its associating functionals in Algorithm \ref{alg:sample_trajectories} of Section \ref{sec: global_section}. 
\end{remark}

\subsubsection{Lifting to a larger filtration}
In the particle-level description of the Dirac branching system in Remark \ref{remark: concrete_description_of_dbs}, each jump involves the selection of a particle $\bx\in \X$ and a jump type $e\in J$. This additional information is in fact not determined by the trajectory of the measure-valued process $(\mu_t)$: one can find examples of $\mu\in E $ such that different choices of $(\bx, e)$ give the same post-jump state $\Phi(\mu; \bx, e)$ under the transition function $\Phi$ in \eqref{eqns: jump_transition_kernel}. In particular, the jump multiplicative functionals defined later in \eqref{eqns: defn_of_measure_jump_MF} explicitly depend on such information. This motivates us to introduce a lifted PDMP on an enlarged state space to retain both the selected particle and jump type at jump events.

To this end we introduce the enlarged state space :
\begin{align}
    E^\lift:= E \times \X\times J,
\end{align}
which is equipped with the product topology (with $J$ equipped with the discrete topology) and remains to be a subspace of some Polish space (hence a Blackwell measurable space). On $E^\lift$, we define the lifted local characteristics using that of the Dirac branching system defined in Definition \ref{defn: Dirac branching system}, namely, 
\begin{align}\label{eqns: definition_lifted_local_characteristics}
    \phi^\lift (t; (\mu, \bx, e)) &= (\phi(t, \mu), \bx, e), \quad (\mu, \bx, e) \in E^{\text{lift}} , \\
    \ell^{\lift} (t; (\mu, \bx, e)) &= \lambda t \abs{\mu}, \\
    \Gamma^\lift(\mu, \bx, e; d\mu', d\bx', de') &= \frac{1}{\abs{\mu}} \sum_{\by\in \mu}\sum_{a\in J} p(a) \delta_{(\Phi(\mu; \by, a), \by, a)} (d\mu', d\bx', de').
\end{align}
We note that $\ell^\lift$ and $\Gamma^\lift$ in \eqref{eqns: definition_lifted_local_characteristics} do not depend on the current mark $(\bx, e)$. 

Let $(\tau_n, (\rho_n, \bx_n, e_n))$ be the embedded Markov chain (defined as in \eqref{eqns: embedded Markov chain}) given by the lifted local characteristics at \eqref{eqns: definition_lifted_local_characteristics} with some deterministic initial condition. If we  define 
\begin{align}\label{eqns: definition of lifted dbs}
    \xi_t := \sum_{n=0}^\infty (\bx_n, e_n) \I_{\tau_n \leq t<\tau_{n+1}},\quad t\geq 0
\end{align}
and 
\begin{align}\label{eqns: measure-valued_dbs_in_MarkovChain}
    \mu_t := \sum_{n = 0}^\infty \phi(t - \tau_n, \rho_n) \I_{\tau_n \leq t <\tau_{n+1}},\quad t\geq 0
\end{align}
then we have:
\begin{align}\label{eqns: lifted_PDMP_dbs}
 (\mu_t, \xi_t) = \sum_{n=0}^\infty \phi^\lift (t - \tau_n; (\rho_n, \bx_n, e_n)) \I_{\tau_n \leq t<\tau_{n+1}}, \quad t \ge 0 .    
\end{align}
Therefore, by definition \ref{defn: PDMP}, $(\mu_t, \xi_t)$ is a PDMP with local characteristics defined at \eqref{eqns: definition_lifted_local_characteristics}. Meanwhile, it is clear that the process $(\tau_n, \rho_n)$ is an embedded Markov chain of the local characteristic of the non-lifted process defined in Definition \ref{defn: Dirac branching system}, so $(\mu_t)_{t\geq 0}$ is a Dirac branching system.

 We denote the natural filtration of $(\mu_t)$ and its lifted process by:
\begin{align}
    \bF := (\F_t)_{t\geq 0}, \quad \bF^\lift:= (\F_t^\lift)_{t\geq 0}
\end{align}
respectively. As described in the introduction of this section, the inclusion:
\begin{align}\label{eqns: filtration inclusion}
    \bF \subset \bF^\lift 
\end{align}
is strict in general.  A useful observation is that the embedded Markov chain $(\tau_n, \rho_n, \bx_n, e_n)$ is determined completely by the marked jump sequence $(\tau_n, \bx_n, e_n)$ and an initial condition since we have the evolution: 
\begin{align}\label{eqns: PDMP_expansion_of_DBS}
    \rho_n = \Phi(\phi( \tau_{n}-\tau_{n-1}, \rho_{n-1}), \bx_n, e_n), \quad n\geq 1, \quad \rho_0 = \mu_0,
\end{align}
for any given initial condition $\mu_0\in E$. Hence, if we define the point process $\widehat N$ on $[0, \infty)\times \X\times J$ using only the marked jump sequence by:
\begin{align}\label{eqns: MPP_of_DBS}
    \widehat N(dt, d\bx, de) := \sum_{n=1}^\infty \delta _{(\tau_n, \bx_n, e_n)}(dt, d\bx, de),
\end{align}
and denote $\bF^{\widehat N} := (\F_t^{\widehat N})_{t\geq 0}$ to be its natural filtration, we have:
\begin{align}\label{eqns: equality of filtration}
    \bF^\lift = \bF^{\widehat N}, 
\end{align}
where the non-trivial inclusion $\bF ^\lift\subset \bF^{\widehat N}$ holds because of \eqref{eqns: PDMP_expansion_of_DBS}. We note that the point process $\widehat N$ has the compensator given by:
\begin{align}\label{eqns: MPP_DBS_Compensator}
    \widehat \nu(dt, d\bx, de) = \lambda p(de) \mu_{t-}(d\bx)dt.
\end{align}
The equality of filtration at \eqref{eqns: equality of filtration} allows us to simplify the subsequent point-process calculus by working directly with $\widehat N$ and its compensator $\widehat \nu$. 
\begin{convention}\label{convention: lifted_PDMP}
Whenever a Dirac branching system $(\mu_t)_{t\geq 0}$ is fixed, it should be understood to be the first component of some lifted process $(\mu_t, \xi_t)$ defined using \eqref{eqns: lifted_PDMP_dbs}, which in turn follows from a pre-defined embedded Markov chain $(\tau_n, (\rho_n, \bx_n, e_n))$ on the state space $[0, \infty)\times E^\lift$. The embedded Markov chain then defines the marked point process $\widehat N$ by \eqref{eqns: MPP_of_DBS} with $\widehat \nu$ by \eqref{eqns: MPP_DBS_Compensator}. 
\end{convention}
\subsubsection{Multiplicative functionals}\label{subsubsection: Mutliplicative_functionals}

We fix a Dirac branching system $(\mu_t)$ and introduce two multiplicative functionals associated to it, or more precisely, its lifted process. We first introduce the continuous multiplicative functionals $(E_t)$ defined by a measurable function $\fu: \X\to \C$ by 
\begin{align}\label{eqns: defn_of_measure_continuous_MF}
    E_t:= \exp\bigg(\int_0^t \inprod{\mu_s, \fu}ds\bigg), \quad t\geq 0.
\end{align}
We  call $\fu$ the  potential of $(E_t)$. Then, we define using a measurable function $\fm: J\times \X\to \C$  the jump multiplicative functionals $(M_t)$ by 
\begin{align}\label{eqns: defn_of_measure_jump_MF}
    M_t:= \prod_{e\in J, 0<s\leq t, \bx\in \X} \fm(e, \bx) ^{\widehat N(ds, d\bx, de)}
\end{align}
in the notation of \eqref{eq:Mt2}. We call such $\fm$ the associated jump potential.    It is not hard to see that $\fu, \fm$ can be naturally extended to the form required in Definition \ref{defn: definition_of_MF_PDMP}. Meanwhile, we make the important observation that while $(E_t)_{t\geq 0}$ is in fact adapted to the measure-valued filtration $\bF$ and is also a multiplicative functional of the measure-valued process $(\mu_t)$, its jump counterpart $(M_t)_{t\geq 0}$ is in general only adapted to the strictly larger filtration $\bF^\lift$ of the lifted process.

\subsection{Stochastic calculus of the Dirac branching system}
\label{subsec: pathwise calculus of DBS}
Let $(\mu_t)$ be the Dirac branching system with parameter $\lambda>0$ and $\{p(e)\}_{e\in J}$ defined in Definition \ref{defn: Dirac branching system}, $(E_t)$ be a continuous multiplicative functional with a potential $\fu:\X\to \C$ defined in \eqref{eqns: defn_of_measure_continuous_MF}, $(M_t)$ be a jump multiplicative function with potential $\fm:J\times \X\to \C$  defined in \eqref{eqns: defn_of_measure_jump_MF}. We shall call  
\begin{align}\label{eqns: standard paramter set}
    \lambda, \{p(e)\}_{e\in J}, \fu, \fm,
\end{align}
the \textit{standard parameter set}. Next, we develop in this section our stochastic calculus approach. 
\subsubsection{Finite jumps over finite time}
To start with, we want to show that $(\mu_t)$ is a regular process, that is, it never has infinite jump in finite time, so that the drifts of processes associated to $(\mu_t)$ can be of finite-variation. To this end, we define:
\begin{align}\label{eqns: pure birth process}
     N_t:= \widehat N([0, t]\times  \X\times J) ,\quad N_b(t) :=  \widehat  N([0, t]\times\X\times J\backslash\{\mathsf{f}\}), \quad t\geq 0,
\end{align}
which are the accumulated number of events and the accumulated number of branchings respectively.
\begin{lemma}\label{lemma: finite_jump_of_Dbs}
For every $t>0$, the process $(\mu_t)$ has at most finitely many jumps over $[0, t]$ almost surely. Equivalently, for every $t > 0 $. 
\begin{align}\label{eqns: pure_birth_process_of_dbs}
    N_t<\infty ,\quad \mbox{almost surely}. 
\end{align}
\end{lemma}
\begin{proof}
Observe that the number of flips is finite between consecutive branching times. Hence, it is equivalent to showing that the number $N_b(t)$ of branching is finite almost surely, which is a pure birth process with linear rate, similar to a Yule process. The result follows from standard results, for instance, Theorem 2.2 of \cite{BansayeMeleard2015}. 
\end{proof}

\subsubsection{It\^o formula and semi-martingale decompositions for the Dirac branching system}
To establish the It\^o formula, that is, the chain rule, we define the class of test functions to be used: for every  measurable function $\fh: \X\to \C$, we define the nonlinear functional  $G_\fh: E\to \C$ on the space of measures by:
\begin{align}\label{eqns: defn_of_nonlinear_functional}
    G_\fh(\mu):= \prod_{\bx\in \mu}\fh(\bx).
\end{align}
As we will use the assumption frequently, we shall write $\fh\in C_b^1$ to mean  
\begin{align}\label{eqns: lifted_regularity}
    \fh(\cdot, 1), \fh(\cdot, -1)\in C_b^1
\end{align}
in the sense of \eqref{eqns: cb1_defn} without any ambiguity. 

Thanks to Lemma \ref{lemma: finite_jump_of_Dbs}, we can now safely use Proposition \ref{prop: semi-martingale_of_PDMP} to establish the following It\^o formulae: 

\begin{proposition}\label{prop: semi-martingality_of_dbs}
Let $\fh:\X\to \C$ be such that $\fh\in C_b^1$ and suppose that the potentials $\fu, \fm$  introduced in \eqref{eqns: defn_of_measure_continuous_MF} and \eqref{eqns: defn_of_measure_jump_MF} are bounded. Then the processes $(E_t), (M_t), (G_\fh(\mu_t))$ are semi-martingales and satisfy the following stochastic differential equations almost surely for all $t\geq 0$:
\begin{align}
    dE_t &= \inprod{\mu_t, \fu}E_t dt, \quad E_0 = 1,\\ 
    dM_t &= \lambda M_{t-} \sum_{e\in J, \bx\in \mu_{t-}}p(e) [\fm(e, \bx) - 1]dt + \int_{\X\times J}M_{t-} (\fm(e, \bx) - 1) \tilde N (dt, d\bx, de),\\
    dG_\fh (\mu_t) &= \sum_{\bx\in \mu_t} \bigg\{\bigg[-z_\bx \p_x \fh (\bx) + \lambda (\sum_{e\in J} p(e)G_\fh (\Phi(\delta_\bx, \bx, e)) - \fh (\bx))\bigg]   \prod_{\by \in \mu_{t-} - \delta_\bx }\fh(\by)\bigg\}dt \\ \label{eqns: semi_martingale test function_dbs}
    &\quad + \int_{\X\times J} [G_\fh (\Phi(\mu_{t-}, \bx, e)) - G_\fh(\mu_{t-})]\tilde N(dt, d\bx, de), 
\end{align}
where 
\begin{align}
    \tilde N(dt, d\bx , de) := \widehat N(dt, d \bx, de) - \widehat \nu (dt, d\bx, de),
\end{align}
is the compensated marked point process (see \eqref{eqns: MPP_of_DBS} and \eqref{eqns: MPP_DBS_Compensator}) and stochastic integral terms against it are local martingales. In addition, the process $(E_t M_t G_\fh(\mu_t))$ given by their product satisfies almost surely:
\begin{align}\label{eqns: dbs_product_formula}
&d(E_t M_t G_\fh (\mu_t)) \\
= &\quad E_{t-}M_{t-}\sum_{\bx\in  \mu_t} G_\fh (\mu_{t} - \delta_\bx)\bigg[-z_\bx\p_x \fh(\bx) +(\fu(\bx)- \lambda)\fh(\bx)  + \lambda\sum_{e\in J}p(e) \fm(e, \bx) G_\fh (\Phi(\delta_\bx; \bx, e)) \bigg] dt \\
&\quad + E_{t-}M_{t-}\int_{\X\times J} G_\fh (\mu_{t-} - \delta_\bx)\bigg[\fm(e, \bx) G_{\fh}(\Phi(\delta_\bx; \bx, e)) - \fh(\bx)\bigg]\tilde N(dt, d\bx, de). 
\end{align}
\end{proposition}
\begin{proof}
Throughout the proof, we work with the natural filtration, identified as $\bF^{\widehat  N}$ following \eqref{eqns: equality of filtration}, of the lifted PDMP introduced at \eqref{eqns: definition of lifted dbs}. Proposition \ref{prop: semi-martingale_of_PDMP} is applied to the lifted process, with the processes $(E_t), (M_t), (G_\fh(\mu_t))$ viewed as functionals of its trajectory. 

From Lemma \ref{lemma: finite_jump_of_Dbs}, the process $(\mu_t)$ has only finitely many events over a finite time. It is straightforward to see that $(E_t), (M_t)$ has a drift term of finite variation and $(M_t)$ has a local martingale part as the stochastic integral with respect to the compensated measure $\tilde N$ from Proposition \ref{prop: semi-martingale_of_PDMP}. For the SDE  satisfied by $(M_t)$, it follows from:
\begin{align}
    dM_t =   \int_{\mathcal X \times J} M_{t-}( \fm(e, \bx)-1)\widehat N( dt, d\bx, de).
\end{align}
For $(G_\fh (\mu_t))$, we follow Proposition \ref{prop: semi-martingale_of_PDMP} with \eqref{eqns: lifting_via_flow} to consider the function 
\begin{align}
   \widehat{G_\fh}(t, \mu):=  G_\fh (\phi(t, \mu)) = \prod_{\bx \in \mu}\fh (x_\bx-z_\bx t, z_\bx, \theta_\bx),
\end{align}
and find that
\begin{align}\label{eqns: derivative_of_flow_lift}
    \frac{d}{dt}\widehat{G_\fh}(t, \mu) = \sum_{\bx\in \mu}(-z_\bx \p_x \fh(x_\bx-z_\bx t, z_\bx, \theta_\bx))\prod_{\by \in \mu-\delta_\bx} \fh (x_\by - z_\by t, z_\by, \theta_\by),
\end{align}
for every $t>0$ and $\mu\in E$. Hence $ \widehat{G_\fh}(t, \mu)$ is of finite variation due to the boundedness assumption of $\fh$ and  Lemma \ref{lemma: finite_jump_of_Dbs}, which justifies the process $(G_\fh(\mu_t))$ being a semimartingale.  Setting $t= 0$ in \eqref{eqns: derivative_of_flow_lift} gives the first part of the drift term, thanks to \eqref{eqns: pdmp_ito_formula}.  By Lemma \ref{lemma: finite_jump_of_Dbs} again, the integral
\begin{align}
    \int_0^t\int_{\X\times J} [G_\fh (\Phi(\mu_{s-}, \bx, e)) - G_\fh (\mu_{s-}) ]\widehat \nu (ds, d\bx, de)
\end{align}
is finite almost surely for all $t>0$. This gives both the expression of the local martingale term against $\tilde N$ and the remaining part of the drift in \eqref{eqns: semi_martingale test function_dbs}. 
 The final product formula is an application of It\^o's product rule. 
\end{proof}
\begin{remark}
From \eqref{eqns: dbs_product_formula}, we can immediately  formulate a restricted local martingale problem in the sense of P. 69 of \cite{Pinsky1991} for the process $(\mu_t)$ and the product $(E_t M_t)$ of multiplicative functionals using an operator defined from the drift of \eqref{eqns: dbs_product_formula} over the space of test functions of the form $G_\fh$ with $\fh\in C_b^1$. Nevertheless, we take a stochastic calculus approach in this paper, explicitly describing the underlying noise of the Dirac branching system $(\mu_t)$ and its multiplicative functionals using marked point processes. 
\end{remark}

We proceed with the standard next step in stochastic calculus to see when the local martingales in \eqref{eqns: dbs_product_formula} are true martingales. 

\begin{corollary}\label{corollary: martingale of product of dbs}
Let $T>0$ and $\fh(\cdot, 1), \fh(\cdot, -1)\in C_b^1$. The process  
\begin{equation*}
\begin{split}
   E_t M_t G_\fh{\mu_t}  - G_\fh (\mu_0)
    - \int_0^t E_s M_s \sum_{\bx\in \mu_s} G_{\fh}(\mu_s - \delta_\bx) & \Big[ -z_\bx \p_x \fh (\bx) + (u(\bx)-\lambda)\fh(\bx) \\
    & \qquad  {} + \lambda \sum_{e\in J} p(e)\fm (e, \bx) G_\fh (\Phi(\delta_\bx; \bx, e))\Big]ds
\end{split}
\end{equation*}
is a true $L_1$ martingale over $[0, T]$,  if
\begin{align}
    \E\Big[\int_0^T\sum _{\bx\in \mu_t} \sum_{e\in J} p(e)\abs{E_{t}M_t G_\fh(\mu_t - \delta_\bx)[ \fm(e, \bx) G_\fh (\Phi(\delta_\bx; \bx, e)) - \fh (\bx)]}dt\Big]<\infty. 
\end{align}
\end{corollary}
\begin{proof}
    It follows from standard sufficient condition for $L_1$ martingality of stochastic integrals against marked point processes. See for instance Section 4.3.2 of \cite{Applebaum09}. 
\end{proof}

\subsection{Semigroup theory of the Dirac branching system}\label{subsec: local_semigroup_quantitative_bound}
Using the stochastic calculus developed in Section \ref{subsec: pathwise calculus of DBS}, we look at the Markov process properties of the Dirac branching system including two fundamental properties to begin with, namely a branching property that extends to its lifted process and a translation homogeneity. 
\begin{proposition}\label{prop: weak_branching_of_dbs}
Let $(\mu_t^1, \widehat N_t^1)$ and $(\mu_t^2, \widehat N_t^2)$ be independent pairs of Dirac branching processes and its associated point process on a common probability space with the same Dirac branching parameters. Then $(\mu_t^1 + \mu_t^2, \widehat N_1 + \widehat N_2)$ is also a pair of the Dirac branching process and its associated point process with the same parameter. 

\end{proposition}
\begin{proof}
It suffices to work at the level of embedded Markov chains: let $(\tau_n^j, \rho_n^j, \bx_n^j, e_n^j)$ be the underlying embedded Markov chain of the pair indexed by $j = 1, 2$. Let $(\hat \tau_n, \hat \rho_n, \hat \bx_n, \hat e_n)$ be the chronological interweaving of the individual chain based on their first component, with the exception of defining:
\begin{align}
    \hat \rho_n := \mu_{\hat \tau_n}^1 +  \mu_{\hat \tau_n}^2.
\end{align}
Such chronological interweaving is well-defined because one can show that the jumping time of the two processes are disjoint almost surely. It is then straightforward to see that:
\begin{align}
    \mu_t^1 + \mu_t^2= \sum_{n = 0}^\infty \phi( t - \hat \tau_n, \hat \rho_n) \I_{\hat \tau_n \leq t <\hat \tau_{n+1}}, \quad 
    \widehat N^1 + \widehat N^2 =  \sum_{n = 1}^\infty \delta_{(\hat \tau_n, \hat \bx_n, \hat e_n)}, \quad t\geq 0.
\end{align}
 It remains to show that $(\hat \tau_n, \hat \rho_n, \hat \bx_n, \hat e_n)$ is a Markov chain with kernel induced by the lifted local characteristic \eqref{eqns: definition_lifted_local_characteristics}, which we defer to Appendix \ref{appendix: sec3_branching}.  
\end{proof}

Next for every $a\in \R$, we define the translation operation $t_a:\X\to \X$ by $t_a \bx := (x_\bx + a, z_\bx, \theta_\bx)$ and its lifting $T_a:E\to E$ by:
\begin{align}\label{eqns: translation_operator_over_measure}
    T_a(\sum_{\bx\in \mu} \delta_\bx) := \sum_{\bx\in \mu} \delta_{t_a\bx }, \quad \mbox{for all }\mu\in E.
\end{align}
The second fundamental distributional property of the Dirac branching system is its  homogeneity with respect to the translation operator:
\begin{proposition}\label{prop: translation_homogeneity}
Let $(\mu_t)$ be a Dirac branching system with its lifted point process:    $\hat N_t = \sum_{n=1}^\infty \delta_{(\tau_n, \bx_n, e_n)}$. Then $(T_a \mu_t)_{t\geq 0}$ is also a Dirac branching process of the same parameter with $\hat N_t^a:= \sum_{n=1}^\infty \delta_{(\tau_n, t_a \bx_n, e_n)}$ to be its lifted point process. 
\end{proposition}
\begin{proof}
Let $(\tau_n, \rho_n, \bx_n, e_n)$ be the associated embedded Markov chain. The proof follows from straight-forward verification that the translated chain $(\tau_n, T_a \rho_n, t_a\bx_n, e_n)$ has a Markov kernel given by \eqref{eqns: definition_lifted_local_characteristics}. 
\end{proof}

\subsubsection{The function $V$ and well-definedness of the expectation functional}
In the following we fix a Dirac branching process $(\mu_t)$ with its multiplicative functionals $(E_t), (M_t)$ parameterized by the standard parameter set \eqref{eqns: standard paramter set} such that the potentials $\fu, \fm$ are bounded measurable functions. For its semigroup theory, we shall only focus on functions of the form $G_\fh$ (see \eqref{eqns: defn_of_nonlinear_functional} for the definition) where $\fh:\X\to \C$ satisfies $\fh\in L_\infty$. To this end, we formally define for every $t>0$ and $\fh \in L_\infty $ the functions 
$\bP_t G_\fh : E\to \C$ and $\bP_t \fh: \X\to \C $ by
\begin{align}\label{eqns: formal_definition_semigroup_dbs}
    \bP_t G_{\fh} (\mu) &:= \E [E_t M_t G_\fh (\mu_t)\mid \mu_0 = \mu], \\
    \bP_t  \fh (\bx)&:= \bP_t G_\fh  (\delta_\bx).
\end{align}
By the branching property in Proposition \ref{prop: weak_branching_of_dbs}, whenever $\bP_t \fh$ is well-defined, we have:
\begin{align}\label{eqns: branching_property_nonlinear_function}
     \bP_t G_{\fh} (\mu) = \prod_{\bx \in \mu} \bP_t \fh(\bx) = G_{\bP_t \fh}(\mu),\quad \mu \in E.  
\end{align}
The study of $\bP_t G_\fh$ is then reduced to that of $\bP_t \fh$ which is well defined under the first moment condition:
\begin{align}\label{eqns: first_moment_semi-group_condition}
    \E[\abs{E_tM_tG_\fh(\mu_t)} \mid \mu_0 = \delta_\bx]<\infty.
\end{align}
for every $\bx\in \X$. To impose a sufficient condition for \eqref{eqns: first_moment_semi-group_condition}, we define:
\begin{align}\label{eqns: supremum_notation_Soler}
     m_e:= \sup_{\bx} \abs{\fm(e, \bx)}, \quad e\in J; \\
     \lambda_e := \lambda p(e), \quad e\in J; \quad \lambda_b:= \lambda(1-p(\mathsf{f}))
\end{align}
Furthermore, we let $\bbeta := (\beta_e)_{e\in J}$ where $\beta_e>0$ for every $e\in J$ and we define $(M_t^{\bbeta})$ to be the jump multiplicative functional with the jump potential:
\begin{align}\label{eqns: defn of majorant jump}
    \fm(e, \bx) := \beta_e, \quad &\mbox{for all }\quad \bx\in \X, \quad e\in J,
\end{align}
and define for every $\eta\in \R, t>0$ and $B>0$ that:
\begin{align}\label{eqns: definition_of_V}
    V(t, B; \eta; \bbeta ):= \E\bigg[\exp(\eta\int_0^t \abs{\mu_s}ds) M^{\bbeta}_t B^{\abs{\mu_t}}\mid \abs{\mu_0}=1\bigg]\in [0, \infty].
\end{align}
It is straight-forward to see that we have:
\begin{align}
    \E[\abs{E_t M_t G_\fh(\mu_t)}\mid \mu_0 = \delta_\bx] \leq V(t, \norm{\fh}_\infty, \sup_{\X}\Re \fu ; (m_e)_{e\in J}),
\end{align}
for every $t>0$. Hence a sufficient condition for \eqref{eqns: first_moment_semi-group_condition} would simply be:
\begin{align}\label{eqns: sufficient condition for first moment condition}
   V(t, \norm{\fh}_\infty, \sup_{\X}\Re \fu ; (m_e)_{e\in J})<\infty. 
\end{align}

The function $V$ defined in \eqref{eqns: definition_of_V} is key to the paper. Our next step is to characterize its finiteness, which, in turn, will give a sufficient condition for \eqref{eqns: first_moment_semi-group_condition} to be satisfied.  
\begin{proposition}\label{prop: property of V}
For every $t>0, B>0, \eta\in \R$ and $\bbeta:= (\beta_e)_{e\in J}$ where $\beta_e>0$ for all $e\in J$, we define:
\begin{align}\label{eqns: definition of C}
    C(\eta, \beta_f)&:= \eta + \beta_f \lambda_f - \lambda,\\ \label{eqns: definition of bbar}
    \bar\beta_b&:= \frac{1}{1-p(\mathsf{f})} \sum_{e\in J\backslash \{\mathsf{f}\}} p(e)\beta_e,\\
    \label{eqns: definition of r}
    r(t, B; \eta, \beta_f, \bar \beta_b)&:=
\begin{cases}
\displaystyle\frac{\bar \beta_b \lambda_b}{4C(\eta, \beta_f)} (e^{2C(\eta, \beta_f)t}-1)B^2, \quad & \text{ \rm if } C(\eta, \beta_f) \neq 0, \\
    \displaystyle\frac{1}{2}\bar \beta_b \lambda_b B^2 t, \quad & \text{ \rm if } C(\eta, \beta_f) = 0 \, . 
\end{cases}
\end{align}
Then, the function $V$ defined in \eqref{eqns: definition_of_V} is well-defined and given explicitly by:
\begin{align}\label{eqns: explicit_formula_of_V}
    V(t, B; \eta; \bbeta ) &= \frac{e^{tC(\eta, \beta_f)}}{\sqrt{1-4r(t, B; \eta, \beta_f, \bar \beta_b)}}B,
\end{align}
if and only if the function $r$ defined at \eqref{eqns: definition of r} satisfies:
\begin{align}\label{eqns: condition for r}
     r(t, B; \eta, \beta_f, \bar \beta_b)<1/4.
\end{align}

\end{proposition}
\begin{proof} 
Note that the estimator in the definition of $V$ in \eqref{eqns: definition_of_V} comes from the pure birth process $(N_b(t))_{t\geq 0}$ of jumps (defined in \eqref{eqns: pure birth process})  of linear rate which has two birth types together with its multiplicative functionals. Standard techniques for birth and death process apply, for instance, one can construct a PDE satisfied by $V$ with respect to the variable of $t$ and $B$. Here we use a probabilistic approach, making use of the distribution of the number of branchings and that of the branching time conditional on this number. We refer the reader to Appendix \ref{appendix: Sec3_PropertyV} for the detailed arguments. 
\end{proof}
It follows immediately that we have:
\begin{corollary}\label{corollary: t_star_condition}
Using all the quantities defined in Proposition \ref{prop: property of V}, we further introduce:
\begin{align}\label{eqns: defn_of_general_t_star}
   t^*(B ;\eta ; \beta_f, \bar \beta_b):=
 \begin{cases}
    \displaystyle  \frac{1}{2C(\eta, \beta_f)} \log \bigg(\frac{C(\eta, \beta_f)}{\bar \beta_b \lambda_b B^2} +1\bigg), \quad & \text{\rm if } C(\eta, \beta_f) >0, \\
    \displaystyle  \frac{1}{2 \bar \beta_b \lambda_b B^2}, \quad & \text{\rm  if } C(\eta, \beta_f) = 0, \\
   \displaystyle    \frac{1}{2C(\eta, \beta_f)} \log \bigg(\frac{C(\eta, \beta_f)}{\bar \beta_b \lambda_b B^2} +1\bigg), \quad & \text{\rm if } -B^2 \lambda_b \bar \beta_b<C(\eta, \beta_f)<0, \\
  \displaystyle     \infty , \quad & \text{\rm if } C(\eta, \beta_f) \leq -B^2 \lambda_b \bar \beta_b \, 
 \end{cases},
\end{align}
which takes values in $[0, \infty]$. Then the condition \eqref{eqns: condition for r} holds if and only if $t< t^*(B ;\eta ; \beta_f, \bar \beta_b)$. 
\end{corollary}
\begin{remark}
     Our proof of probabilistic representation shares with \cite[Theorem 3.2]{HT21} the strategy of  using a majorizing quantity, which is the function $V$ in our case.  Meanwhile, in our stochastic calculus approach, the additional exponential process 
     $(\exp(\eta \int_0^t \lvert \mu_s \rvert ds))_{t\geq 0}$, 
     which is adapted to the strictly smaller filtration $\bF$, arises naturally. The result is a delicate finiteness condition as reflected in the case-by-case formula of $t^*$ in Corollary \ref{corollary: t_star_condition}. 
\end{remark}

Since our main interest is whether \eqref{eqns: sufficient condition for first moment condition} is satisfied, whenever a standard parameter set  of the triple $(\mu_t, E_t, M_t)$ is fixed, we simplify notations to introduce for every $t, B\geq 0$:
\begin{align}\label{eqns: definition_of_V_sim}
    V(t, B)&:=  V(t, B, \sup_\X \Re \fu,(m_e)_{e\in J}),\quad t, B\geq 0;\\ \label{eqns: definition_of_t_star_sim}
      t^*(B)&:= t^*(B,\sup_\X \Re \fu, m_f, \bar m_b),\quad B>0,\quad \bar m_b:= \sum_{e\in J\backslash\{\mathsf{f}\}}\frac{p(e) m_e}{1-p(\mathsf{f})}.
\end{align}
 We can now address when the formally defined operator $\bP_t$  would be well-defined:
\begin{proposition}\label{prop: property of semigroup of dbs}
Let $\fh\in L_\infty$ and $s, t\geq 0$. Suppose $\fu, \fm \in L_\infty$. 
\begin{enumerate}
    \item Suppose $t<t^*(\norm{\fh}_\infty)$. Then $\bP_t \fh\in L_\infty$ is well-defined with:
\begin{align}\label{eqns: norm_control_semi_group}
    \norm{\bP_t \fh}_\infty \leq V(t, \norm{\fh}_\infty).
\end{align}
If in addition we have $\fh, \fu, \fm\in C_b^1$, then we have $\bP_t h\in C_b^1$. 
    \item Suppose $s+t<t^*(\norm{\fh}_\infty)$. Then we have the following equality of functions over $\X$:
\begin{align}\label{eqns: semi_group_over_X}
    \bP_{t+s} \fh = \bP_t (\bP_s \fh) = \bP_s (\bP_t \fh),
\end{align}
in which all terms are well-defined. In addition, we also have the following equality of functions over $E$:
\begin{align}\label{eqns: semi_group_over_E}
    \bP_{t+s} G_\fh = \bP_t (\bP_s G_\fh ) =  \bP_s (\bP_t G_\fh ). 
\end{align}
\end{enumerate}
    
\end{proposition}
\begin{proof}
    For the first statement, the well-definedness and the norm bound for $\bP_t \fh$ in \eqref{eqns: norm_control_semi_group} is an immediate result of Proposition \ref{prop: property of V}. For the differentiability, the homogeneity with respect to translations established in Proposition \ref{prop: translation_homogeneity} allows us to write:
\begin{align}
    \bP_t \fh(\bx) = \E[E_t M_t G_\fh(\mu_t)\mid \mu_0 = \delta_{\bx}] = \E[E_t^{x_\bx} M_t^{x_\bx} G_\fh (T_{x_\bx} \mu_t)\mid \mu_0 = \delta_{0, z_\bx, \theta_\bx}], 
\end{align}
where $T_{x_\bx}$ is the translation operator defined in \eqref{eqns: translation_operator_over_measure}  and $E_t^{x_\bx} M_t^{x_\bx}$ are the modified multiplicative functionals after translating the underlying process. The proof then follows once we can verify the interchange of the expectation and spatial differentiation with respect to $x_\bx$.  We leave the details to Appendix \ref{appendix: Sec3_Differentiability} which includes the precise definition of $E_t^{x_\bx}$ and $M_t^{x_\bx}$.
    The equalities \eqref{eqns: semi_group_over_X} of the second result follow from the first result and the equality:
\begin{align}
    V(t+s, \norm{\fh}_\infty) = V(t, V(s, \norm{\fh}_\infty)) = V(s, V(t, \norm{\fh}_\infty)),
\end{align}
which can be checked for instance by the observation that $V$ is the expectation over multiplicative functionals of a pure birth process in its definition \eqref{eqns: definition_of_V}, or simply by the explicit formula \eqref{eqns: explicit_formula_of_V} established in Proposition \ref{prop: property of V}. The equalities \eqref{eqns: semi_group_over_E} then follows from applying the weak branching property at \eqref{eqns: branching_property_nonlinear_function} on \eqref{eqns: semi_group_over_X}. 
\end{proof}
\subsubsection{The generator}
In the following we strengthen the assumption of $\fh$ and let $\fh\in C_b^1$.  We formally define $\bL G_{\fh}: E\to \C$ by:
\begin{align}\label{eqns: point-wise generator - nonlinear}
    \bL G_\fh(\mu):=\lim_{t\to 0^+} \frac{\bP_{t} G_\fh(\mu) - G_\fh(\mu)}{t},\quad \mu\in E,
\end{align}
We shall see \eqref{eqns: point-wise generator - nonlinear} is always well-defined. 
\begin{proposition}\label{prop: measure_generator}
  Suppose $\fu, \fm\in C_b^1$. Then  for all $\fh \in C_b^1$, the function $\bL G_\fh:E\to \C$  is well-defined. In addition for every $\bx\in \X$, we have the explicit formula:
    \begin{align}
    \bL G_\fh (\delta_\bx)=-z_\bx \p_x \fh (x_\bx, z_\bx, \theta_\bx) +(\fu(\bx)-\lambda)\fh(\bx)+\lambda \sum_{e\in J}p(e) m(e, \bx) G_\fh (\Phi(\delta_\bx; \bx, e)), \quad \bx\in \X. 
\end{align}
Therefore the function $\bL \fh: \X\to \C$ defined by $\bL \fh (\bx) := \bL G_\fh (\delta_\bx)$ for every $\bx\in \X$ is bounded, that is, $\norm{\bL \fh}_\infty<\infty $. 
\end{proposition}
\begin{proof}
 By Proposition \ref{prop: property of semigroup of dbs}, $\bP_t G_\fh$ is well-defined for small $t>0$. Hence, we can write the numerator in \eqref{eqns: point-wise generator - nonlinear} as a difference of expectation. We can then apply \eqref{eqns: dbs_product_formula} from Proposition \ref{prop: semi-martingality_of_dbs} to rewrite the numerator as an expectation over sums of integrals. The result follows by applying the martingale condition in Corollary \ref{corollary: martingale of product of dbs} over a sufficiently small time interval and justifying the interchange of limit and expectation. Since the domain of the operators $(\bP_t)$ are changing, in particular shrinking as $t\to \infty$, we do not have the usual strong continuity to exchange the limit and expectation. However, one can still justify this last step via a dominated convergence argument, and we present the full proof of the result in Appendix \ref{appendix: Sec3_Generator}. 
\end{proof}
\subsection{Probabilistic representation for the Soler model}\label{subsec: proof_of_main_theorem}
Our strategy to write the probabilistic representation for the Soler model \eqref{eqns: sec4_nonlinear_Dirac} is to first rewrite the Soler model so that we can \textit{directly see}  how we specify  parameters for a Dirac branching system together with its multiplicative functionals. The probabilistic representation then follows from the Markov process theory established in Section \ref{subsec: local_semigroup_quantitative_bound}. 

\subsubsection{Rewriting the Soler system}
We first establish the notation and define for every $f:\R\times \{\pm 1\}\to \C$, the lifted function $\bf :\X\to \C$ by:
\begin{align}\label{defn: boldface}
    {\pmb f}(\bx):= f(x_\bx, z_\bx)\I_{\theta_\bx = 1} + \overline{f(x_\bx, z_\bx)}\I_{\theta_\bx = -1},\quad \mbox{for every}\quad \bx\in \X,
\end{align}
which we shall call its boldface. Now, let $g\in C_b^1$.  Note that $g\in C_b^1$ if and only if $\bg \in C_b^1$. Taking the complex conjugation on \eqref{eqns: sec4_nonlinear_Dirac}, we can obtain the following lifted evolutionary equation over the enlarged state space $\X$:
\begin{equation}\label{eqns: lifted_non_linear_Dirac}
\begin{cases}
    \p_t \tilde{\psi}(t, \bx) = -z_\bx\p_x \tilde{\psi}(t,\bx)-i\theta_\bx m \tilde \psi(t, x_\bx, -z_\bx, \theta_\bx)+i\theta_\bx \mathcal N(\tilde{\psi})\tilde{\psi}(t, x_\bx, -z_\bx,\theta_\bx);\\
     \tilde{\psi}(0, \bx) = \bg(\bx), \quad (t, \bx)\in [0, \infty)\times \X;\\
    \mathcal N(\tilde{\psi}) = \tilde{\psi}(t, x_\bx, z_\bx, 1)\tilde{\psi}(t, x_\bx, -z_\bx, -1)+\tilde{\psi}(t, x_\bx, -z_\bx, 1)\tilde{\psi}(t, x_\bx, z_\bx, -1),
\end{cases}
\end{equation}
with $\tilde \psi(t,\cdot):\X\to \C$ being the solution.  It is straight-forward to see that if $\tilde \psi$ is a classical solution over $[0, T)$ with $T\in [0, \infty]$ to the lifted equation \eqref{eqns: lifted_non_linear_Dirac} in the sense that $\tilde \psi$ is differentiable with respect to $t, x_\bx$ such that $\tilde \psi(t, \cdot)\in C_b^1$ and
\begin{align}\label{eqns: lifted_classical_solution_defn}
\sup_{s\in [0, t]} \bigg(\norm{\tilde \psi(s, \cdot)}_\infty + \norm{\p_x\tilde \psi(s, \cdot)}_\infty \bigg)<\infty 
\end{align}
for all $t\in [0, T)$ where we denote $\p_{x_\bx}$ as $\p_x$, then $\psi(t, x, z):= \tilde \psi(t, x, z, 1)$ is a classical solution to \eqref{eqns: sec4_nonlinear_Dirac} in the sense defined in \eqref{eqns: classical_solution_defn} and \eqref{eqns: classical_solution_additional_boundedness_condition}. Conversely, if $\psi$ is a classical solution to \eqref{eqns: sec4_nonlinear_Dirac} then the boldface $\bpsi$ defined in \eqref{defn: boldface} is a classical solution to \eqref{eqns: lifted_non_linear_Dirac} over $[0, T)$ for $T\in [0, \infty]$. 

Furthermore, we observe that for any measurable function $\fh:\X\to \C$, using the jump transition defined in \eqref{eqns: jump transition function} would give: 
\begin{align}\label{eqns: non_linear_observation}
    G_\fh(\Phi(\delta_\bx, \mathsf{f})) &=\fh(x_\bx, -z_\bx, \theta_\bx),\\
    G_\fh(\Phi(\delta_\bx, \mathsf{b_1})) &=\fh(x_\bx, -z_\bx, \theta_\bx)\fh(x_\bx, 1, 1)\fh(x_\bx, -1, -1),\\
    G_\fh(\Phi(\delta_\bx, \mathsf{b_2})) &=\fh(x_\bx, -z_\bx, \theta_\bx)\fh(x_\bx, -1, 1)\fh(x_\bx, 1, -1).
\end{align}
Therefore, we can rewrite the lifted Soler model \eqref{eqns: lifted_non_linear_Dirac} using the identities \eqref{eqns: non_linear_observation} as:
\begin{equation}\label{eqns: transformed_lifted_non_linear_Dirac}
\begin{cases}
   \hfill  \p_t \tilde{\psi}(t, \bx) &= -z_\bx\p_x \tilde{\psi}(t, \bx)\\
   &\quad -im\theta_\bx  G_{\tilde \psi(t,\cdot)}(\Phi(\delta_\bx, \mathsf{f}))+i\theta_\bx G_{\tilde \psi(t,\cdot)}(\Phi(\delta_\bx, \mathsf{b_1}))+i\theta_\bx G_{\tilde \psi(t,\cdot)}(\Phi(\delta_\bx, \mathsf{b_2}));\\
    \hfill  \tilde{\psi}(0, \bx) &= \bg(\bx), \quad (t, \bx)\in [0, \infty)\times \X. 
\end{cases}
\end{equation}
At this point, it is clear that the original problem of solving for a classical solution of the Soler model with initial condition $g\in C_b^1$ is equivalent to solving for a classical solution in the sense of \eqref{eqns: lifted_classical_solution_defn} for \eqref{eqns: transformed_lifted_non_linear_Dirac} with initial condition $\bg$ induced by $g$. The way we see from \eqref{eqns: transformed_lifted_non_linear_Dirac} the assumptions we want to impose on the standard parameter set \eqref{eqns: standard paramter set} is by looking at its coefficients.
\begin{assumption}\label{assumption: non-linear_coefficients}
The parameters $\lambda, \{p(e)\}_{e\in J}, \fu, \fm$ defining the processes $(\mu_t), (E_t), (M_t)$ shall satisfy:
    \begin{align}
    \fu(\bx) - \lambda = 0, &\quad \lambda \fm(\mathsf{f}, \bx)p(\mathsf{f}) = -im\theta_\bx,\\
    \lambda \fm(\mathsf{b_1}, \bx)p(\mathsf{b_1}) = i\theta_\bx ,&\quad  \lambda \fm(\mathsf{b_2}, \bx)p(\mathsf{b_2}) = i\theta_\bx,
\end{align}
for every $\bx\in \X$. 
\end{assumption}
\begin{remark}\label{remark: degree of freedom}
    Note that under Assumption \ref{assumption: non-linear_coefficients}, as soon as the Dirac branching parameters $\lambda>0$ and $\{p(e)\}_{e\in J}$ are chosen, the potentials $\fu, \fm$ are determined by the coefficients of the Soler system (assuming the mass parameter $m>0$ of the Soler system is given). In particular, the flipping multiplier $\fm(\mathsf{f}, \cdot)$ is determined by the coefficient of the linear interaction term, the branching multipliers $\fm(\mathsf{b_j}, \cdot)$ with $j = 1, 2$ are determined by the coefficients of the cubic terms. In fact, the continuous potential $\fu$ is determined by the linear potential term, which vanishes in the case of the Soler system. The same correspondence will appear in Assumption \ref{assumption: general_cubic_assumption} of the general case in Section \ref{subsec: bounded cubic coefficients}. 
\end{remark}
An important observation is that under Assumption \ref{assumption: non-linear_coefficients}, we have:
\begin{align}
   C(\sup_\X \Re \fu, m_f) :=\sup_\X \Re \fu +\lambda_fm_f - \lambda  = m>0,\quad \lambda_b\bar \beta_b = 1+1=2,
\end{align}
where $C$ and $\bar \beta_b$ are defined at \eqref{eqns: definition of C} and \eqref{eqns: definition of bbar}. The expression depends only on the mass parameter of the Soler system. Hence, we are always in the first condition in the definition of $t^*$ defined in Corollary \ref{corollary: t_star_condition}, and the simplified notation in \eqref{eqns: definition_of_t_star_sim} gives:
\begin{align}\label{eqns: definition_of_t_star_Soler}
    t^*(B) = \frac{1}{2m}\log \bigg(\frac{m}{2B^2}+1\bigg), \quad B>0,
\end{align}
which is a bijective function of $B>0$ and so $(t^*)^{-1}:(0, \infty)\to (0, \infty)$ is well-defined. 
\subsubsection{Probabilistic representation for the Soler model}
In this section, we consider the Dirac branching system $(\mu_t)$ together with two of its multiplicative functionals $(E_t), (M_t)$ which we recall from the beginning of Section \ref{subsec: dirac branching system} such that the standard parameter set \eqref{eqns: standard paramter set} defining them satisfies  Assumption \ref{assumption: non-linear_coefficients}. The main goal of the section is to establish the probabilistic representation for the Soler system \eqref{eqns: sec4_nonlinear_Dirac}. More concretely, we establish the relation between the transformed  system \eqref{eqns: transformed_lifted_non_linear_Dirac} with initial condition $\bg\in C_b^1$ (with definition in \eqref{eqns: lifted_regularity}) and the expectation functional $\bP_t \bg $ defined in \eqref{eqns: formal_definition_semigroup_dbs}.  

To this end, we first establish the weaker result of probabilistic representation under an additional $L_\infty$ a priori estimate:
\begin{lemma}\label{lemma: non_linear_necessary_condition}
Let $g\in C_b^1(\R\times \{\pm 1\})$ and $\bg$ be its boldface. Let $(\mu_t, E_t, M_t)_{t\geq 0}$ be a triple such that their standard parameter set satisfies Assumption \ref{assumption: non-linear_coefficients}.  Let $\tilde \psi(t, \bx)$ be a classical solution to \eqref{eqns: transformed_lifted_non_linear_Dirac} with respect to $\bg$ over $[0, T)$ for some $T>0$ with the additional a priori  estimate:
\begin{align}\label{eqns: necessary_boundedness_assumption_nonlinear}
   \sup_{t\in [0, T)} \norm{\tilde \psi(t, \cdot)}_\infty<(t^*)^{-1}(T),
\end{align}
where $t^*$ is defined at \eqref{eqns: definition_of_t_star_Soler}. 
    Then we have 
\begin{align}
    \tilde \psi(t, \bx) = \bP_t \bg (\bx) =\E[E_t M_t G_\bg(\mu_t)\mid \mu_0 = \delta_\bx],\quad \bx \in \mathcal X , t\in [0, T) \,  . 
\end{align}
\end{lemma}
\begin{proof}
We use a standard It\^o formula proof and provide the details in Appendix \ref{appendix: Sec3_Necessary_Condition}. In short, we first fix $t<T$ and establish the reverse-time PDE satisfied by $\tilde \psi(t-s, \bx)$ over $s\in [0, t]$. Then we apply a time-augmented It\^o formula to those in Proposition \ref{prop: semi-martingality_of_dbs} with a test function given by $\tilde \psi(t-\cdot,\cdot )$. The drift of the resulting semi-martingale vanishes because of the PDE $\tilde \psi(t-\cdot, \cdot)$ satisfies and the local martingale term is a true martingale as we can  now apply Proposition \ref{prop: property of V} due to the a priori estimate \eqref{eqns: necessary_boundedness_assumption_nonlinear}. The result basically follows by taking expectations.
\end{proof}

Here comes the main theorem of the paper:
\begin{theorem}\label{thm: sec4_main_theorem}
Let $g\in C_b^1(\R\times \{\pm 1\})$ and $\bg$ be its boldface. Let $(\mu_t, E_t, M_t)_{t\geq 0}$ be a triple such that their standard parameter set satisfies Assumption \ref{assumption: non-linear_coefficients}. Then for every $T>0$ with 
\begin{align}\label{eqns: upper_bound_of_T_local_wellposedness}
    T< t^*(\norm{g}_\infty) = \frac{1}{2m}\log \bigg(\frac{m}{2\norm{g}_\infty^2}+1\bigg),
\end{align}
 the system \eqref{eqns: transformed_lifted_non_linear_Dirac} of equations  has a unique classical solution $\tilde \psi(t,\bx)$ over $[0, T)$  given by 
\begin{align}\label{eqns: the probabilistic representation_main_theorem}
    \tilde \psi(t, \bx) = \bP_t \bg (\bx) =\E[E_t M_t G_\bg(\mu_t)\mid \mu_0 = \delta_\bx].
\end{align} 
More concretely with the notation  \eqref{eq:Mt2}, under the Soler model \eqref{eqns: sec4_nonlinear_Dirac} with the initial condition $g\in C_b^1$ in \eqref{eqns: cb1_defn},  the same well-posedness result holds with a unique classical solution  
\begin{align}
    &\psi(t, x, z) = \bP_t \bg(x, z, 1)\\
    &=\E \Big[e^{\lambda\int_0^t\abs{\mu_s}ds}\prod_{\substack{0<s\leq t, \\e\in J, \bx\in \X}} \!\!\! \fm(e, \bx)^{\widehat N(ds, d\bx, de)} \prod_{\bx\in \mu_t}g(x_\bx, z_\bx)^{\I_{\{\theta_\bx = 1\}}}\, \overline{g(x_\bx, z_\bx)}^{\I_{\{\theta_\bx = -1\}}}\mid \mu_0 = \delta_{x, z, 1}\Big] \, . 
\end{align}
for $x \in \mathbb R$, $z \in \{\pm 1\}$, $t \in [0, T)$. 
\end{theorem}

\begin{proof}
First, we show that $\bP_t \bg(\bx)$ is a classical solution over $[0, T)$.  Let $\epsilon>0$ be such that $T+\epsilon <t^*(\norm{g}_\infty)$. Since the potentials $\fu, \fm$ of the multiplicative functionals $(E_t), (M_t)$ are independent of the position variable $x_\bx$ under Assumption \ref{assumption: non-linear_coefficients}, together with the choice of $T$ and $\epsilon>0$, we can apply Proposition \ref{prop: property of semigroup of dbs} to claim that $\bP_t \bg  \in C_b^1$ for every $t\in [0, T+\epsilon)$. Because $\bP_t \bg  \in C_b^1$, we can further apply Proposition \ref{prop: measure_generator} to obtain $\bL G_{\bP_t\bg }$ which satisfies 
\begin{align}\label{eqns: concrete_formula_in_proof_main_theorem}
    \bL G_{\bP_t\bg } (\delta_\bx) = -z_\bx \p_x \bP_t \bg(\bx) + \lambda \sum_{e\in J}p(e) \fm(e, \bx) G_{\bP_t \bg}(\Phi(\delta_\bx; \bx, e)), \quad t\in [0, T+\epsilon). 
\end{align}
By unraveling the definition of $\bL$ in $\bL G_{\bP_t g}$, we have 
\begin{align}
    \bL G_{\bP_t g} (\mu) = \lim_{h\to 0}\frac{\bP_h \bP_t G_\bg(\mu) - \bP_t G_\bg(\mu)}{h},
\end{align}
in which we used the fact that $G_{\bP_t \bg} = \bP_t G_\bg$ in \eqref{eqns: branching_property_nonlinear_function}. 
Considering small enough $h>0$ such that $t+h<T+\epsilon $, we can make use of the semi-group property at \eqref{eqns: semi_group_over_E} in Proposition \ref{prop: property of semigroup of dbs} and obtain that $\bP_t g(\bx)$ has a right-derivative equal to $\bL G_{\bP_t g}(\delta_\bx)$ for all $t\in [0, T)$. Left differentiability follows from a similar argument of Proposition \ref{prop: measure_generator} in Appendix \ref{appendix: Sec3_Generator}, applied over $[t-h, t]$ with sufficiently small $h$, with the only difference that the averaged drift at \eqref{eqns: difference quotient and drift} converges to its pre-jump value at time $t-$. This difference can be addressed by noting that every fixed $t$ is almost surely not a jump time as in the derivation of \eqref{eqns: left_limit_sam_as_original_time}. Combining differentiability from both sides, we have
\begin{align}
    \frac{d}{dt}\bP_t G_\bg (\delta_\bx) = \frac{d}{dt} \bP_t \bg(\bx)  =\bL G_{\bP_t g}(\delta_x), \quad t\in [0, T),
\end{align}
which is bounded due to the explicit formula in \eqref{eqns: concrete_formula_in_proof_main_theorem}, which shows that $\bP_t \bg$ satisfies \eqref{eqns: transformed_lifted_non_linear_Dirac}. The uniform boundedness on compact time intervals \eqref{eqns: classical_solution_additional_boundedness_condition} follows from \eqref{eqns: norm_control_semi_group}, \eqref{eqns: final_upper_bound} and the continuity of $V$ and its parameter derivatives.

For the uniqueness of the probabilistic representation, we  apply a continuation argument. Suppose we have a solution $\tilde \psi(t, \cdot)$ over $[0, T)$ that is a classical solution such that derivatives are bounded. 
   Using \eqref{eqns: lifted_non_linear_Dirac} to transfer the boundedness condition \eqref{eqns: classical_solution_additional_boundedness_condition} from the spatial derivative of $\tilde \psi$ to its time derivative, the map $t\mapsto \tilde\psi(t, \cdot)$ is continuous in $L_\infty$ and so the running supremum  $t\mapsto \sup_{t'\in [0, t]} \lVert{\tilde\psi(t',\cdot)}\rVert_\infty$ is continuously increasing from $\norm{\bg}_\infty $. Because $\lim_{t\to 0^+}(t^*)^{-1}(t) = \infty$, it follows that there must exist $\tau\in(0, T)$ such that $\sup_{t\in [0, \tau]} \norm{\psi(t, \cdot)}_\infty <(t^*)^{-1}(\tau)$. It follows from Lemma \ref{lemma: non_linear_necessary_condition} that $\tilde\psi(t, \cdot)$ must be given by the probabilistic representation \eqref{eqns: the probabilistic representation_main_theorem} over $[0, \tau]$. Hence we can define:
\begin{align}
    t^*:= \sup\{s\in [0, T):\psi(t,\cdot)\mbox{ has probabilistic representation over $[0, s]$} \}.
\end{align}
Suppose we had $t^* <T$. Note that $\norm{\psi(t^*,\cdot)}$ is bounded and so by the existence argument we have a probabilistic representation using $\psi(t^*, \cdot)$ as the initial condition over $[0, \epsilon]$ for some $\epsilon>0$ such that $t^*+\epsilon<T<t^*(\norm{g}_\infty)$. By the semigroup property in Proposition \ref{prop: property of semigroup of dbs}, we can extend the probabilistic representation of $\psi(t, \cdot)$ from $t^*$ to $t^* + \epsilon$ which is a contradiction.  
\end{proof}
\begin{remark}
Although the choice of $(\mu_t, E_t, M_t)_{t\geq 0}$ in the probabilistic representation \eqref{eqns: the probabilistic representation_main_theorem} is not unique and has a degree of freedom given by the Dirac branching parameter $\lambda>0$ and $p(e)>0$ for every $e\in J$,  the expectation expression there determines a unique value. In addition, the upper bound for the admissible time interval $[0, T)$ at \eqref{eqns: upper_bound_of_T_local_wellposedness} is also independent from such choice. We will see in Section \ref{subsec: one_step_MC} that such choice is related to the variance reduction in the Monte Carlo implementation of the representation. 
\end{remark}

The probabilistic representation for the Soler model in Theorem \ref{thm: sec4_main_theorem} in fact unveils a deeper question with the following simple observation:
\begin{corollary}\label{corollary: sharpeness results_constant_data}
Let $g\equiv B>0$ be a constant function and $T>0$. Then \eqref{eqns: transformed_lifted_non_linear_Dirac} has a unique classical solution  given by the probabilistic representation \eqref{eqns: the probabilistic representation_main_theorem} over $t\in [0, T)$ if \textbf{and only if} $T$ satisfies \eqref{eqns: upper_bound_of_T_local_wellposedness} with $\leq$ replacing $<$. 
\end{corollary}
\begin{proof}
    The if part follows from Theorem \ref{thm: sec4_main_theorem}, while the only if part follows from constancy of $g$ and that of the potentials which gives the \textit{equality}:
\begin{align}
    \E[\abs{E_t M_t G_\bg(\mu_t)}] = V\Big(t, B; \lambda ; \Big(\frac{m}{\lambda_f}, \frac{1}{\lambda_{b_1}}, \frac{1}{\lambda_{b_2}}\Big) \Big),
\end{align}
which is finite over $[0, T)$ if and only if \eqref{eqns: upper_bound_of_T_local_wellposedness} holds. 
\end{proof}
Thus, it is evident that the local-in-time probabilistic representation for the Soler model is an intrinsic limitation of our one-step representation with the $L_\infty$-type initial data. Nevertheless, global well-posedness as well as large-time numerical algorithms of the Soler model is well-known in the literature for standard function spaces including $L_2$ integrable functions (see work reviewed in Section \ref{sec: introduction}). It raises the question of whether the Soler system genuinely has a finite time span over $L_\infty$-type initial data, and if not, how  we should extend our probabilistic representation, together with the associated Monte Carlo methods, for the global solution. We will address these questions in Section \ref{sec: global_section}.

\section{Extension to more general nonlinear Dirac equations}\label{sec: state_dependent_coefficients}
In this section, we illustrate the robustness of our stochastic calculus method of establishing a probabilistic representation by extending it from the case of the Soler model (which is the content of Section \ref{sec: nonlinear_case}) to the case of the Soler system with more general Dirac operators and the case of more general cubic nonlinearity.

\subsection{Soler system with more general Dirac operators}
Already in the foundational paper \cite{Dirac1928},  Dirac equations were first formulated as a class of differential equations parameterized by a matrix pair 
\begin{align}\label{eqns: Clifford pair}
    (\gamma_0, \gamma_1),
\end{align}
 where $\gamma_0, \gamma_1$ are $2\times 2$ complex-valued matrices that satisfy the anti-commutativity relation: 
\begin{align}\label{eqns: anti_commute_relation}
    \gamma_i \gamma_j + \gamma_j \gamma_i = 2\eta_{ij} I_2,\quad i, j = 0, 1,
\end{align}
where $I_2$ is the two dimensional identity matrix and $\eta := (\eta_{ij})_{0\le i,j \le 1}$ forms another $2\times 2$ matrix, called the Lorentz signature defined by 
\begin{align}
    \eta = 
\begin{bmatrix}
    \eta_{00} & \eta_{01}\\
    \eta_{10} & \eta_{11}
\end{bmatrix}=
\begin{bmatrix}
    1 & 0 \\
    0 & -1
\end{bmatrix}.
\end{align}
Under the Soler model, the parameterized family is given by:
\begin{align}\label{eqns: matrix_pair_Soler}
    i\gamma_0 \p_t \psi(t, x) + i\gamma_1 \p_x \psi(t, x) = m\psi(t, x) - \psi(t, x)^\dagger\gamma_0 \psi(t, x) \psi(t, x),\quad \psi(0, x) = g(x),
\end{align}
where $m>0$, $g(x) = 
\begin{bmatrix}
    g(x, 1)\\
    g(x, -1)
\end{bmatrix}:\R\to \C^2
, \psi(t, x) = 
\begin{bmatrix}
    \psi(t, x, 1) \\
    \psi(t, x, -1)
\end{bmatrix}
:[0,  \infty)\times \R\to \C^2$, and $\dagger$ denotes the conjugate transpose. As usual, we look at the case where the initial condition $g$ satisfies $g\in C_b^1$ and the solution $\psi$ is a classical solution in the sense given by \eqref{eqns: classical_solution_defn} and \eqref{eqns: classical_solution_additional_boundedness_condition}. In Section \ref{sec: nonlinear_case}, we were looking at the case:
\begin{align}\label{eqns: our_choice_matrix_pair}
    \gamma_0 := \sigma_1, \quad \gamma_1:= -i\sigma_2,
\end{align}
where $(\sigma_1, \sigma_2, \sigma_3)$ are the Pauli matrices defined by:
\begin{align}
    \sigma_1:=
\begin{bmatrix}
    0 & 1 \\
    1 & 0
\end{bmatrix}, \quad 
\sigma_2:=
\begin{bmatrix}
    0 & -i\\
    i & 0
\end{bmatrix},\quad \sigma_3:=
\begin{bmatrix}
    1 & 0 \\
    0 & -1
\end{bmatrix}.
\end{align}
We made such choice because the particular matrix pair at \eqref{eqns: our_choice_matrix_pair} is more natural for a  probabilistic representation theory, and in fact the choice of such matrix pair is simply off by a change of coordinates:

\begin{lemma}\label{lemma: CliffordPair_Change_of_coordinates}
Suppose $(\gamma_0, \gamma_1)$ and $(\gamma_0', \gamma_1')$ are matrix pairs that satisfy \eqref{eqns: anti_commute_relation}. Then there exists an invertible matrix $S$ such that:
\begin{align}
    S^{-1} \gamma_j S = \gamma_j',\quad j = 0, 1.
\end{align}
In the case where both matrices in the matrix pairs are either Hermitian or anti-Hermitian\footnote{A $2\times 2$ complex-valued matrix $A$ is Hermitian if the conjugate transpose $A^\dagger$ satisfies: $A^\dagger = A$; it is anti-Hermitian if $A^\dagger = -A$.}, $S$ can be chosen to be unitary. 
\end{lemma}
\begin{proof}
While most references claim this two-dimensional result as a folklore or cite the fundamental theorem by Pauli \cite{Pauli1936}, which indeed proved only the four-dimensional case (see for instance Lemma 2.1 of \cite{ComechGuanGustafson2014} and Remark 2 of \cite{OzawaYamauchi2004}), we are providing a simple self-contained proof in Appendix \ref{appendix: Sec4_CliffordPair}. 
\end{proof}

\begin{theorem}\label{thm: well-posedness_different_matrix_pair}
For  any matrix pair $(\gamma_0, \gamma_1)$ satisfying \eqref{eqns: anti_commute_relation} that are either Hermitian or anti-Hermitian, we let $U$ be the unitary matrix (due to Lemma \ref{lemma: CliffordPair_Change_of_coordinates}) such that:
\begin{align}
    \gamma_0 = U^\dagger\sigma_1 U, \quad \gamma_1 = U^\dagger (-i\sigma_2) U.
\end{align}
Let $g\in C_b^1$ and define $h\in C_b^1$ by: $h(x, \cdot):= Ug(x,\cdot)$ for every $x\in \R$. Then for any $T>0$ with:
\begin{align}
    T<t^*(\norm{h}_\infty), 
\end{align}
where $t^*$ is defined at \eqref{eqns: definition_of_t_star_Soler}, the general Soler model \eqref{eqns: matrix_pair_Soler} has a unique classical solution $\psi(t, x, z)$ over $[0, T)$ given by:
\begin{align}\label{eqns: matrix_pair_probabilistic_representation}
    \psi(t, x, \cdot) := U^\dagger \bP_t \bh (x, \cdot, 1)=U^\dagger\begin{bmatrix}
        \E[E_t M_t G_{\pmb {h}}(\mu_t)\mid \mu_0 = \delta_{x, 1, 1}]\\
         \E[E_t M_t G_{\pmb {h}}(\mu_t)\mid \mu_0 = \delta_{x, -1, 1}]
    \end{bmatrix},
\end{align}
for every  $t\in [0, T)$, where $(\mu_t), (E_t), (M_t)$ is a Dirac branching system with its multiplicative functionals whose standard parameter set (defined in \eqref{eqns: standard paramter set}) satisfies Assumption \ref{assumption: non-linear_coefficients} and $\bP_t$ is defined in \eqref{eqns: formal_definition_semigroup_dbs} with the bold face $\bh$ defined in \eqref{defn: boldface}.  
 \end{theorem}
\begin{proof}
We observe if $(\gamma_0, \gamma_1)$ and $(\gamma_0', \gamma_1')$ are matrix pairs such that:
\begin{align}
    \gamma_i= U^\dagger \gamma_i'U, \quad i = 0, 1,
\end{align}
then a function $\psi$ is a classical solution to the Dirac equation with initial condition $g$ and matrix pair $(\gamma_0, \gamma_1)$ if and only if $U\psi$ is a classical solution to the Dirac equation with initial condition $Ug$ and matrix pair $(\gamma_0', \gamma_1')$. Since under such condition of $T$, the function $\bP_t \bh(x, z, 1)$ over $t\in [0, T)$ is the unique solution to the matrix pair $(\sigma_1, -i\sigma_2)$ by Theorem \ref{thm: sec4_main_theorem}, it follows that $U^\dagger \bP_th(x)$, which is the right-hand-side of \eqref{eqns: matrix_pair_probabilistic_representation}, is the required unique solution to the matrix pair $(\gamma_0, \gamma_1)$. 
\end{proof}

\subsection{Cubic nonlinearity with bounded coefficients}\label{subsec: bounded cubic coefficients}
The goal of the section is to demonstrate the robustness of the techniques developed in Section \ref{sec: nonlinear_case} by establishing the probabilistic representation of the general nonlinear cubic Dirac equations, which take the form:
\begin{equation}\label{eqns: sect5_general_cubic_nonlinear_Dirac}
\begin{split}
   \p_t \psi(t, x, z)&= -z\p_x \psi(t, x, z)+ \alpha_{-1}(x, z) \psi(t, x, -z) +\alpha_1(x, z) \psi(t, x, z),\\
   &\quad +\sum_{\abs{\bk} = 3} \alpha_{\bk}(x, z)\psi(t, x, z)^{k_{1,1}} \psi(t, x, -z)^{k_{1,-1}}\overline{\psi}(t, x, z)^{k_{-1,1}}\overline{\psi}(t, x, -z)^{k_{-1,-1}}\\
   \psi(0, x, z) &= g(x, z)
\end{split},
\end{equation}
where $\bk:= (k_{1,1}, k_{1,-1}, k_{-1,1}, k_{-1,-1})\in \mathbb{Z}_{\geq 0}^4$ are 4-tuples of non-negative integers and we define $\abs{\bk}$ to be its $\ell_1$ sum. We assume that the coefficients:
\begin{align}\label{eqns: sec5_intro_bounded_coefficients}
\balpha:= \{\alpha_1, \alpha_{-1}, \alpha_{\bk}:\abs{\bk} = 3\}\subset \{\alpha:\R\times \{\pm 1\}\to \C\}
\end{align}
are complex-valued  and are in $C_b^1$ so their derivatives are also bounded. We also require $\alpha_{\bk}$ not to be all vanishing for every $\abs{\bk} = 3$.  We assume the initial condition $g$ to satisfy $g\in C_b^1$ and consider the solution $\psi$ to be a classical solution in the sense of \eqref{eqns: classical_solution_defn} and \eqref{eqns: classical_solution_additional_boundedness_condition}.

The underlying stochastic process of \eqref{eqns: sect5_general_cubic_nonlinear_Dirac} is basically the same as the Dirac branching system defined in Definition \ref{defn: Dirac branching system} for the Soler model, except that we use the enlarged jump symbol:
\begin{align}
\widehat J:=\{\mathsf{f}, \mathsf{b}_{\bk}:\abs{\bk}=3\},
\end{align}
which has $1+20$ elements, and the jump transition function $\widehat \Phi:E\times \X \times \widehat J\to E$ defined by:
\begin{align}\label{eqns: extended jump transition function}
\Phi(\mu; \bx, \mathsf{f})&:= \mu - \delta_{\bx} + \delta_{(x_\bx, -z_\bx, \theta_\bx)},\\
\Phi(\mu; \bx, \mathsf{b}_{\bk})&:= \mu - \delta_{\bx} + k_{1,1}\cdot \delta_{x_\bx, z_\bx, \theta_\bx}+k_{-1,1}\cdot\delta_{x_\bx, z_\bx, -\theta_\bx}\\
&\quad +k_{1,-1}\cdot\delta_{x_\bx, -z_\bx, \theta_\bx}+k_{-1,-1}\cdot \delta_{x_\bx, -z_\bx, -\theta_\bx},
\end{align}
for every $\bx\in \mu$. This underlying process, which we call the generalized Dirac branching system, is thus parameterized by the event rate $\lambda>0$ and the event probabilities $\{p(e)\}_{e\in \widehat J}$. We remark that an analog of Convention \ref{convention: lifted_PDMP} regarding an associated lifted PDMP remains in effect and we will use the same notation for the associated objects. 

In the remainder of this section, $(\mu_t)$ is always a generalized Dirac branching system and we write  $(E_t)$ and $(M_t)$ for its multiplicative functionals which are defined analogously to \eqref{eqns: defn_of_measure_continuous_MF} and \eqref{eqns: defn_of_measure_jump_MF} respectively using bounded potentials $\fu$ and $\fm$ with the modification that the domain of $\fm$ is $\widehat J\times \X$ instead of $J\times \X$. We call the parameters:
\begin{align}\label{eqns: generalized_parameter_set}
    \lambda, \{p(e)\}_{e\in \widehat J}, \fu,\fm,
\end{align}
the generalized parameter set. 

Similar to how we obtained Assumption \ref{assumption: non-linear_coefficients}, we note that \eqref{eqns: sect5_general_cubic_nonlinear_Dirac} can be re-written using the jump transition function \eqref{eqns: extended jump transition function} as:
\begin{equation}\label{eqns: sect5_rewritten_general_cubic_nonlinear_Dirac}
\begin{split}
   \p_t \psi(t, x, z)&= -z\p_x \psi(t, x, z)+ \alpha_{-1}(x, z) \psi(t, x, -z) +\alpha_1(x, z) \psi(t, x, z) \\
   &\quad +\sum_{\abs{\bk} = 3} \alpha_{\bk}(x, z)G_{\bpsi(t,\cdot)}(\Phi(\delta_{x,z,1};b_{\bk})),\\
   \psi(0, x, z) &= g(x, z)
\end{split}.
\end{equation}
This motivates us to specify $(\mu_t), (E_t), (M_t)$ by imposing the generalized parameter set \eqref{eqns: generalized_parameter_set} to satisfy:
\begin{assumption}\label{assumption: general_cubic_assumption}
    \begin{align}\label{eqns: general_relation_between_process_and_PDE}
    \fu(\bx) - \lambda = \alpha_{1}(\bx),\quad \lambda p(\mathsf{f})\fm(\mathsf{f}, \bx) = \alpha_{-1}(\bx), \quad \lambda p(\mathsf{b_\bk})\fm(\mathsf{b}_{\bk}, \bx) =\alpha_{\bk} (\bx),\quad \bx\in \X,
\end{align}
where the coefficients $\alpha_1, \alpha_{-1}, \alpha_{\bk}$, when evaluated on $\X$, use the boldface convention introduced at \eqref{defn: boldface}. If $\alpha_{\bk}\equiv 0$ then we use the convention that $m(b_{\bk}, \cdot) \equiv 1$ with $p(b_{\bk}) = 0$.  
\end{assumption}

Similar to the discussion in Section \ref{subsec: local_semigroup_quantitative_bound}, we define $(m_e)_{e\in \widehat J}, (\lambda_e)_{e\in \widehat J}$ and $\lambda_b$ analogous to \eqref{eqns: supremum_notation_Soler}. In addition, for $\widehat\bbeta:= (\beta_e)_{e\in \widehat J}$ where $\beta_e>0$ for every $e\in J$, we define $M_t^{\widehat \bbeta}$ analogous to that in \eqref{eqns: defn of majorant jump} together with:
\begin{align}\label{eqns: definition_of_wideV}
    \widehat V(t, B; \eta; \widehat \bbeta ):= \E\bigg[\exp(\eta\int_0^t \abs{\mu_s}ds) M^{\widehat \bbeta}_t B^{\abs{\mu_t}}\mid \abs{\mu_0}=1\bigg]\in [0, \infty].
\end{align}
similar to  that of $V$ defined at \eqref{eqns: definition_of_V}. It is straight-forward to see that we have for every $t\geq 0$ and $\bx\in \X$ that:
\begin{align}\label{eqns: general first moment upper bound}
    \E[\abs{E_t M_t G_\fh(\mu_t)}\mid \mu_0 = \delta_\bx] \leq \widehat V(t, \norm{\fh}_\infty, \sup_{\X}\Re \fu ; (m_e)_{e\in \widehat J}). 
\end{align}

With obvious modifications to the definitions introduced in Proposition \ref{prop: property of V} and Corollary \ref{corollary: t_star_condition}, both results continue to hold if every instance of $V, J, \bbeta$ in the statements is replaced by $\widehat V, \widehat J$ and $\widehat \bbeta$, with basically the same proof in Appendix \ref{appendix: Sec3_PropertyV}. Now the crucial observation is that under Assumption \ref{assumption: general_cubic_assumption}, the quantities defining the upper bound in \eqref{eqns: general first moment upper bound} satisfy:
\begin{align}\label{eqns: general_case_quantity_simplification}
   C(\sup_\X \Re \fu, m_f)  &= \sup_\X \Re \alpha_1+\norm{\alpha_{-1}}_\infty,\\
   \lambda_b \bar \beta_b &= \sum_{\abs{\bk}=3} \norm{\alpha_\bk}_\infty, 
\end{align}
which depends only on the coefficients of the PDE.  Hence the quantity  $t^*(B)$ defined at \eqref{eqns: definition_of_t_star_sim} also depends only on the coefficients of the PDE. It follows that we obtain basically the same result as Theorem \ref{thm: sec4_main_theorem}:

\begin{theorem}\label{thm: general_cubic_main_theorem}
Let $g\in C_b^1(\R\times \{\pm 1\})$ and $\bg$ be its boldface. Let $(\mu_t, E_t, M_t)_{t\geq 0}$ be a triple such that their generalized parameter set \eqref{eqns: generalized_parameter_set} satisfies Assumption \ref{assumption: general_cubic_assumption}. Then for every $T>0$ with:
\begin{align}\label{eqns: upper_bound_of_T_local_wellposedness_general_cubic}
 T&< t^*(\norm{g}_\infty),
\end{align}
 the generalized equation  \eqref{eqns: sect5_general_cubic_nonlinear_Dirac} has a unique classical solution $\psi(t, x, z)$ over the finite interval $[0, T)$  given by:
\begin{align}\label{eqns: main_thm_general_expectation_formula}
    \psi(t, x, z) = \E\bigg[E_t M_t G_{\bg}(\mu_t)\mid \mu_0 = \delta_{x, z, 1}\bigg], \quad t\in [0, T).
\end{align}    
\end{theorem}
\begin{proof}
Basically the same proof in Theorem \ref{thm: sec4_main_theorem} runs through, except that in the  proof of existence, the potentials $\fu, \fm$ are no longer constants. Since we impose a differentiability assumption on the coefficients, we can still apply results analogous to Proposition \ref{prop: property of semigroup of dbs} so that the corresponding result analogous to \ref{prop: measure_generator} applies. 

Meanwhile, in the proof of uniqueness, the inverse $(t^*)^{-1}$ appearing  in the proof in Theorem \ref{thm: sec4_main_theorem} remains valid if $ C(\sup_\X \Re \fu, m_f)  \geq 0$. In the case of $ C(\sup_\X \Re \fu, m_f)<0$, we shall instead consider the inverse of the restriction of $t^*(B)$ in the open interval $B\in (\sqrt{\frac{-C(\sup_\X \Re \fu, m_f)}{\lambda_b \bar \beta_b}}, \infty )$ because $t^*(B) = \infty$ when $B$ is out of the interval as in the explicit formula in \eqref{eqns: definition_of_t_star_sim}. With such a  modification, the proof still goes through. 
\end{proof}

Unlike the Soler model where the quantity $ C(\sup_\X \Re \fu, m_f)$ defined at \eqref{eqns: general_case_quantity_simplification} is the always positive mass parameter, since $\sup_\X \Re \alpha_1$ can take any values, all cases of $t^*$ established in Corollary \ref{corollary: t_star_condition} can hold as already reflected in the proof of Theorem \ref{thm: general_cubic_main_theorem}. We shall make clear the phase transition phenomenon arising due to such difference:

\begin{corollary}\label{corollary: classification_of_stochastic_representation}
Let $g\in C_b^1(\R\times \{\pm 1\})$ and $\bg$ be its boldface. Let $(\mu_t, E_t, M_t)_{t\geq 0}$ be a triple such that their generalized parameter set \eqref{eqns: generalized_parameter_set} satisfies Assumption \ref{assumption: general_cubic_assumption}.
\begin{enumerate}
    \item Suppose $\sup_\X \Re \alpha_1+\norm{\alpha_{-1}}_\infty >  0$. If
\begin{align}\label{eqns: upper_bound_specific_general_cubic}
    T< \frac{1}{2(\sup_\X \Re \alpha_1 + \norm{\alpha_{-1}}_\infty )}\log \bigg( \frac{\sup_X \Re \alpha_1 + \norm{\alpha_{-1}}_\infty }{\norm{g}_\infty^2 \sum_{\abs{\bk} = 3}\norm{\alpha_{\bk}}_\infty}+1\bigg).
\end{align}
    then  the generalized equation  \eqref{eqns: sect5_general_cubic_nonlinear_Dirac} has a unique classical solution $\psi(t, x, z)$ over the finite interval $[0, T)$  given by \eqref
{eqns: main_thm_general_expectation_formula}
    \item Suppose $\sup_\X \Re \alpha_1+\norm{\alpha_{-1}}_\infty=  0$. If 
\begin{align}
    T< \frac{1}{2\sum_{\abs{\bk} = 3}  \norm{\alpha_\bk}_\infty \norm{g}_\infty ^2},
\end{align}
 then  the generalized equation  \eqref{eqns: sect5_general_cubic_nonlinear_Dirac} has a unique classical solution $\psi(t, x, z)$ over the finite interval $[0, T)$  given by \eqref
{eqns: main_thm_general_expectation_formula}.
\item Suppose $\sup_\X \Re \alpha_1+\norm{\alpha_{-1}}_\infty <0$. We have two cases: under the small norm condition:
\begin{align}\label{eqns: small_norm upper bound}
    \norm{g}_\infty \leq \sqrt{\frac{\sup_\X \Re \alpha_1 + \norm{\alpha_{-1}}_\infty }{-\sum_{\abs{\bk} = 3} \norm{\alpha_\bk}_\infty }}
\end{align}
the generalized equation  \eqref{eqns: sect5_general_cubic_nonlinear_Dirac} has a unique classical \textbf{global} solution $\psi(t, x, z)$ given by \eqref{eqns: main_thm_general_expectation_formula} for all $t\geq 0$. Otherwise, it has a unique classical solution $\psi(t, x, z)$ given by \eqref
{eqns: main_thm_general_expectation_formula} over the finite interval $[0, T)$  for all $T$ satisfying \eqref{eqns: upper_bound_specific_general_cubic}
\end{enumerate}
\end{corollary}
\begin{proof}
    It follows from writing the explicit formula of $t^*(B)$ at \eqref{eqns: definition_of_t_star_sim} using the relation \eqref{eqns: general_case_quantity_simplification}.
\end{proof}

\section{Monte Carlo method and global continuation for the Soler model}\label{sec: global_section}
The one-step probabilistic representation theorems (Theorem \ref{thm: sec4_main_theorem} and Theorem \ref{thm: general_cubic_main_theorem}) reveal a distinction that becomes central to the remainder of the paper. While our one-step representation theorems have, in most cases, led to local well-posedness results of the non-linear Dirac equations, existing literature, especially in the case of the Soler model, has established global well-posedness results \cite{Huh2013GrossNeveu, HuhMoon2015,ZhangZhao2015}, as well as numerics with a relatively long time, which we refer readers to the early paper \cite{AC81} and the survey \cite{XuShaoTang2013}. These results, however, do not directly provide the global $C_b^1$ theory required by the present probabilistic framework, which allows bounded and non-decaying initial data. As we mentioned at the end of Section \ref{sec: nonlinear_case}, it is evident that the local-in-time well-posedness is at best what we can derive from the one-step representation, this leads us to ask whether a global-in-time well-posedness is possible  in the $C_b^1$ framework, and in such case, establish the relation between  the two different notions of lifespan, one for the deterministic Soler model  and one for the probabilistic Dirac branching system.

In this section, we address the ``local probabilistic representation, but global PDE well-posedness" phenomenon for the Soler model. The basic idea is that in the case of the Soler model we can extend the one-step probabilistic representation to a multi-step probabilistic representation thanks to an $L_\infty$ a priori estimate. At the PDE level, this is demonstrated in the time-marching representation we term in Theorem \ref{thm: global_well_posednesss_L_inf} while at the numerical level, this is demonstrated in the time-marching Monte Carlo algorithm established at Algorithm \ref{alg:global-monte-carlo-time_marching} with a rigorous convergence analysis at Theorem \ref{thm: convergence_of_global_alg}. 
\subsection{One-step Monte Carlo solver}\label{subsec: one_step_MC}
We begin by studying the Monte Carlo implementation of the one-step probabilistic representation at \eqref{eqns: main_thm_general_expectation_formula} for the general nonlinear Dirac equation \eqref{eqns: sect5_general_cubic_nonlinear_Dirac}; the first step is to be able to sample from the joint distribution of the generalized Dirac branching system $(\mu_t)$ together with its multiplicative functionals $(E_t), (M_t)$.
As all processes are adapted to the lifted filtration $\bF^\lift$ generated by the point process associated with the embedded Markov chain $(\tau_n, \rho_n, \bx_n, e_n)$, the straight-forward way is to sample from the point process and build $(E_t), (M_t), (\mu_t)$ from the samples. We note that:

\begin{align}
    E_t &= \exp\bigg(\sum_{\tau_n<t} \int_{\tau_n}^{\min\{t, \tau_{n+1}\}} \!\!\!\!\!\!\!\!\!\inprod{\mu_s, \fu}ds\bigg) = \exp\bigg(\sum_{\tau_n<t}\sum_{\by\in \rho_n} \int_{\tau_n}^{\min\{t, \tau_{n+1}\}} \!\!\!\!\!\! \fu(x_\by-(s-\tau_n)z_\by, z_\by, \theta_\by)ds\bigg)\\
    &=\exp\bigg(\sum_{\tau_n<t}\sum_{\by\in \rho_n} \int_{0}^{\min\{t, \tau_{n+1}\}-\tau_n} \!\!\!\!\! \fu(x_\by-sz_\by, z_\by, \theta_\by)ds\bigg) \\
    &= \prod_{\tau_n<t}\prod_{\by\in \rho_n} \exp\bigg(\int_{0}^{\min\{t, \tau_{n+1}\}-\tau_n} \fu(x_\by-sz_\by, z_\by, \theta_\by)ds\bigg),
\end{align}
and 
\begin{align}
    M_t = \prod_{\tau_n\leq t,n\neq 0}\fm(e_n, \bx_n).
\end{align}
Meanwhile, the expression of $(\mu_t)$ written in terms of the embedded Markov chain is at \eqref{eqns: measure-valued_dbs_in_MarkovChain}.  

With these expressions, we can now formulate our sampler for the joint process $(\mu_t, E_t, M_t)$ with a given generalized parameter set (defined at \eqref{eqns: generalized_parameter_set}) and an initial condition $\mu_0\in E$, which we assume to be a singleton for simplicity. 
\begin{breakablealgorithm}
\caption{Sampler for $(\mu_t, E_t, M_t)$ starting at $\mu_0 = \delta_{x, z, 1}$}
\label{alg:sample_trajectories}
\begin{algorithmic}[1]
\Require Generalized parameter set: $\lambda, \{p(e)\}_{e\in \widehat J}, \fm, \fu$; an initial state $x, z$ and time $t$.
    \State $\tau\gets 0$, $\mu\gets\delta_{(x,z,1)}$,
    $E\gets 1$, and $M\gets 1$

    \While{$\tau<t$}
        \State Set $r\gets\lambda|\mu|$, sample
        $\sigma\sim\operatorname{Exp}(r)$, and set
        $\Delta\gets\min\{\sigma,t-\tau\}$

        \State Update $\widehat{\mu}\gets\phi(\Delta,\mu)$,
        and $\tau\gets \tau+\Delta$
        \State Update
    \begin{align}\label{eqns: algorithm_integrator_step}
        E\gets E \prod_{\by\in \mu} \exp\bigg(\int_{0}^{\Delta } \fu(x_\by-sz_\by, z_\by, \theta_\by)ds\bigg)
    \end{align}
        \If{$\tau =t$}
            \State $\mu\gets\widehat{\mu}$
            \State \textbf{break}
        \EndIf

        \State Sample $e\in \widehat J$ according to
        $\{p(e)\}_{e\in \widehat J}$ and choose $\bx$ uniformly
        from $\widehat{\mu}$

        \State Update $M\gets M\,\fm(e,\bx )$ and
        $\mu\gets\Phi(\widehat{\mu};\bx,e)$
    \EndWhile

    \State  $(\mu_t, E_{t}, M_{t})\gets (\mu, E, M)$

\State \Return $(\mu_t, E_{t}, M_{t})$
\end{algorithmic}
\end{breakablealgorithm}
For simplicity, we assume the integral step \eqref{eqns: algorithm_integrator_step} can be performed efficiently, for instance when the potential $\fu$ admits an  anti-derivative representable using elementary functions along the flow map.
\begin{remark}
    The update rules for $(E_t)$ and $(M_t)$ in algorithm \ref{alg:sample_trajectories} are consistent with the filtration that they are adapted to (see Section \ref{subsubsection: Mutliplicative_functionals}). In particular $(E_t)$ is updated before sampling both the jump types and particles in each loop, using strictly less information. 
\end{remark}

To build a numerical solver to the nonlinear Dirac equation \eqref{eqns: sect5_general_cubic_nonlinear_Dirac}, in the following we fix a set of coefficients $\balpha$ (defined in \eqref{eqns: sec5_intro_bounded_coefficients}) and an initial condition $g\in C_b^1$. Now let  $N\in \mathbb{N}$ and let $(\mu_{t,x,z}^{(i)},E_{t, x, z}^{(i)}, M_{t, x, z}^{(i)}, )_{i=1}^N$ for $N\in \mathbb{N}$ be independent copies of the joint process $( \mu_{t, x, z}, E_{t, x, z}, M_{t, x, z})$ which are $(\mu_t, E_t, M_t)$ starting from $\mu_0:= \delta_{x, z, 1}$ such that the generalized parameter set satisfy Assumption \ref{assumption: general_cubic_assumption}. 

With Algorithm \ref{alg:sample_trajectories}, we can now build the estimator:
\begin{align}\label{eqns: local_MonteCarlo_estimator}
    \hat P_{t, N} g (x, z):= \frac{1}{N}\sum_{i=1}^N E_{t, x, z}^{(i)}M_{t, x, z}^{(i)} G_\bg(\mu_{t, x, z}^{(i)}),
\end{align}
for every $t>0, x\in \R, z= \pm 1, g\in C_b^1$. Theorem \ref{thm: general_cubic_main_theorem} implies that whenever $t$ is small enough to satisfy \eqref{eqns: upper_bound_of_T_local_wellposedness_general_cubic}, the estimator \eqref{eqns: local_MonteCarlo_estimator} is an unbiased and consistent estimator for the solution $\psi(t, x, z)$ to \eqref{eqns: sect5_general_cubic_nonlinear_Dirac}. To ensure $\hat P_{t, N}g(x, z)$ is within an error of $\epsilon>0$ from $\psi(t, x, z)$ with constant probability, we consider:
\begin{align}
    \P(\abs{\hat P_{t, N}g(x, z) - \psi(t, x, z)}\geq \epsilon ) \leq \frac{1}{N\epsilon^2}\var(E_{t, x, z} M_{t,x, z} G_{\bg}(\mu_{t, x, z})).
\end{align}
Using the quantity $\widehat V$ defined in \eqref{eqns: definition_of_wideV}, we have the upper bound:
\begin{align}\label{eqns: variance_MC}
    \var(E_{t, x, z} M_{t,x, z} G_{\bg}(\mu_{t, x, z})) &\leq \E[\abs{E_{t, x, z} M_{t,x, z} G_{\bg}(\mu_{t, x, z})}^2] \\
    &\leq \widehat V(t, \norm{g}_\infty^2, 2\sup_\X \Re \fu ;  (m_e^2)_{e\in \widehat J}),
\end{align}
for every $x\in \R$ and $z = \pm 1$. 
 The quantities defining the upper bound of the second moment at \eqref{eqns: variance_MC} are given by 
\begin{align} \label{eqns: definition of second moment C}
    C(2\sup_\X \Re \fu, m_f^2) &= 2\sup_\X \Re \alpha_1 +\lambda + \frac{\norm{\alpha_{-1}}_\infty^2}{\lambda_f}, \\
    \lambda_b \bar \beta_b &= \sum_{\abs{\bk} = 3} \frac{\norm{\alpha_\bk}_\infty^2}{\lambda_{b_{\bk }}}.
\end{align}
Consistent with Assumption \ref{assumption: general_cubic_assumption}, we adopt the convention $0/0:=0$ for terms with zero event rates in \eqref{eqns: definition of second moment C}. Unlike the case for the first moment upper bound at \eqref{eqns: general first moment upper bound}, \eqref{eqns: definition of second moment C} also depends on the Dirac branching parameters $\lambda>0$ and $\{p(e)\}_{e\in \widehat J}$. This suggests that the freedom in choosing the parameters could be exploited for variance reduction in the Monte Carlo estimation. One way to exploit this freedom is to minimize the upper bound at \eqref{eqns: variance_MC}, for which we have an explicit formula analogous to \eqref{eqns: explicit_formula_of_V}. Although minimizing the upper bound does not necessarily minimize the variance itself,  we will follow this observation in our numerical experiments in Section \ref{subsec: numerical experiments}. To emphasize the dependence on the Dirac branching parameters, we use the equivalent event-rate parameterization:
\begin{align}
    \Lambda:= (\lambda_e)_{e\in \widehat J}\in [0, \infty)^{\abs{\widehat J}},\quad \lambda_e:=\lambda p(e),\quad e\in \widehat J,
\end{align}
and write the upper bound at \eqref{eqns: variance_MC} as a function of $\Lambda$:
\begin{align}\label{eqns: second_moment_upper_bound}
    V_t(\Lambda):=\widehat V(t, \norm{g}_\infty^2, 2\sup_\X \Re \fu ;  (m_e^2)_{e\in \widehat J}),
\end{align}
assuming a fixed initial condition $g\in C_b^1$. 

We summarize the section with the following pseudocode of the Monte Carlo algorithm:
\newpage 
\begin{breakablealgorithm}
\caption{One-step Monte Carlo solver for $\psi(t,x,z)$}
\label{alg:local_MC_solver}
\begin{algorithmic}[1]
\Require The initial condition $g\in C_b^1$; the coefficients of the PDE: $\balpha$ (defined at \eqref{eqns: sec5_intro_bounded_coefficients}); target state $t, x, z$; number of samples $N$; event rates $\Lambda$
\Ensure $t$ satisfies \eqref{eqns: upper_bound_of_T_local_wellposedness_general_cubic} and $\inf_\Lambda V_t(\Lambda)<\infty $
\State Choose the event rates $\Lambda$ to minimize the variance $\var(E_t M_t G_\bg(\mu_t))$ and the corresponding Dirac branching parameter $\lambda>0, p(e)$ with $e\in \widehat J$. 
\State Choose the generalized parameter set: $\lambda, \{p(e)\}_{e\in \widehat J}, \fm, \fu$ to satisfy Assumption \ref{assumption: general_cubic_assumption}.  

\State Apply \textbf{Algorithm \ref{alg:sample_trajectories}} with $N$ samples
starting from $\mu_0=\delta_{x,z,1}$ to obtain
\[
\left(
E_{t,x,z}^{(i)},
M_{t,x,z}^{(i)},
\mu_{t,x,z}^{(i)}
\right)_{i=1}^N
\]

\State Compute
\[
\widehat P_{t,N}g(x,z)
\gets
\frac{1}{N}
\sum_{i=1}^N
E_{t,x,z}^{(i)}
M_{t,x,z}^{(i)}
G_{\bg}\!\left(\mu_{t,x,z}^{(i)}\right)
\]

\State \Return $\widehat P_{t,N}g(x,z)$
\end{algorithmic}
\end{breakablealgorithm}

\subsection{Global continuation for the Soler model and the time-marching Monte Carlo} \label{subsec: global_analysis}
At this point, it should be clear that there are two distinct obstructions to numerically implementing the probabilistic representation of Theorem \ref{thm: sec4_main_theorem} and Theorem \ref{thm: general_cubic_main_theorem} over a long time horizon. First, the representation itself is subject to a first-moment integrability constraint as seen from Corollary \ref{corollary: sharpeness results_constant_data}. Meanwhile, to have a Monte Carlo error of constant probability, a second-moment constraint is seen at \eqref{eqns: variance_MC}. Neither obstruction, however, necessarily reflects a finite-time blowup of the underlying deterministic PDE. For the Soler model, an $L^\infty$ a priori estimate allows us to overcome both difficulties of long-time application of the local probabilistic representation by repeatedly restarting it over sufficiently short time intervals. 

In the following, we focus ourselves to the Soler model \eqref{eqns: sec4_nonlinear_Dirac} and begin by establishing the aforementioned $L^\infty$ a priori estimate, which is a modification from that of \cite{Huh2013GrossNeveu} who also uses a similar estimate to establish the global well-posedness in the Sobolev initial data regime $H^s$ with $s\geq 0$. After that we use this observation to obtain a global time-marching representation, and then turning the same construction into the global Monte Carlo method of Algorithm \ref{alg:global-monte-carlo-time_marching}.

\subsubsection{An $L_\infty$ a priori estimate}
 First we define the quantity that we will use frequently below: let $t>0$ and $B>0$, we define:
\begin{align}
   L(t, B):=\bigg[(1+2mt) B^2\exp(t(m+8B^2))\bigg]^{1/2}.
\end{align}

\begin{lemma}\label{lemma: apriori_L_inf_estimate}
   Suppose $\psi(t, x, z)$ is a classical solution  to \eqref{eqns: sec4_nonlinear_Dirac} with initial condition $g\in C_b^1$ over $[0, T)$ for some $T>0$. Then for every $t\in [0, T)$: 
\begin{align}\label{eqns: apriori_L_inf_bound}
    \norm{\psi(t,\cdot)}_\infty \leq L(t, \norm{g}_\infty),
\end{align}
\end{lemma}
\begin{proof}
We are following the calculations in Section 2 of \cite{Huh2013GrossNeveu} with slight modifications. For the convenience of the readers, we state and prove the following modified result from Section 2 of \cite{Huh2013GrossNeveu} in Appendix \ref{sec: proof of Lemma_L_inf_estimate}. 
\end{proof}

\subsubsection{An $L_\infty$-type global well-posedness of the Soler model}
While it is standard to apply a continuation argument to establish a global well-posedness result with a local well-posedness result at hand together with an a priori estimate, we establish instead a constructive proof, showing the availability of a time-marching probabilistic representation that we define in the next result. 

In all remaining subsections of Section \ref{sec: global_section}, we fix a Dirac branching system and its multiplicative functionals satisfying Assumption \ref{assumption: non-linear_coefficients} and let $\bP_t$ denote its expectation functional defined in \eqref{eqns: formal_definition_semigroup_dbs}. To distinguish the probabilistic representation from the solution semi-group, we shall define for every $g\in C_b^1$ and its bold face $\bg\in C_b^1$  the function:
\begin{align}
    P_t g(x, z):= \bP_t \bg (x, z, 1),\quad x\in \R, z= \pm 1,
\end{align}
for every $t\geq 0$ whenever the right-hand side is well-defined.  

\begin{theorem}\label{thm: global_well_posednesss_L_inf}
 For the initial data $g\in C_b^1$, 
 there exists a unique global classical solution $\psi(t,\cdot)$ of the Soler model \eqref{eqns: sec4_nonlinear_Dirac} that is defined for every $t>0$. In particular, for every $T>0$ and for every $\tau>0$ with:
\begin{align}\label{eqns: theoretical_upper_bound_for_tau}
    \tau<t^*(L(T, \norm{g}_\infty)),
\end{align}
the function:
\begin{align}\label{eqns: power_time_marching_representation}
    \psi(t, \cdot) = P_{t-[t/\tau]\tau}P_\tau^{[t/\tau]}g,\quad t\in [0, T)
\end{align}
defines a classical solution over $[0, T)$. 
We call \eqref{eqns: power_time_marching_representation} the time marching representation for $\psi(t, \cdot)$ with time marching parameter $\tau$ over $[0, T)$.  
\end{theorem}
\begin{proof}
First we verify the candidate function $\psi(t, \cdot)$ defined by \eqref{eqns: power_time_marching_representation} is well-defined and is a classical solution to the Soler model with an inductive argument. Let $\tau>\epsilon>0$ be such that:
\begin{align}
    \tau + \epsilon <t^*(L(T, \norm{g}_\infty)).
\end{align}
Since $t^*$ is decreasing, we have
\begin{align}
    \tau+\epsilon <t^*(L(T, \norm{g}_\infty)) \leq t^*(\norm{g}_\infty).
\end{align}
Hence by Theorem \ref{thm: sec4_main_theorem}, the function $P_s g$ is the unique classical solution of the Soler system over $[0, \tau+\epsilon]$. In particular $\psi(\tau, \cdot) = P_\tau g$ satisfies:
\begin{align}
    \tau+\epsilon <t^*(L(T, \norm{g}_\infty)) \leq t^*(\norm{\psi(\tau, \cdot) }_\infty).
\end{align}
We can apply the uniqueness in Theorem \ref{thm: sec4_main_theorem} to claim that $\psi(s,\cdot)$ is well-defined over $[0, \tau+\epsilon)$ and is equal to $P_s g$ for $s\in [\tau, \tau+\epsilon]$. Hence $\psi$ is differentiable in time at the glueing point $\tau$, which is an interior point. 
The boundedness condition \eqref{eqns: classical_solution_additional_boundedness_condition} follows by combining the bounds for the two representations on $[0, \tau]$ and $[\tau, \tau+\epsilon)$ respectively. Now suppose $\psi(t, \cdot)$ is well-defined up to $[0, k\tau+\epsilon ]$ for some $1\leq k\leq [T/\tau]$ and is a classical solution of the Soler system on the interval. We repeat the argument in the first step using $\psi(k\tau, \cdot)$ as the initial condition so that the formula defining $\psi$ gives a classical solution over $[0, (k+1)\tau+\epsilon)$. Hence finally the representation \eqref{eqns: power_time_marching_representation} gives a classical solution to the Soler model over $[0, ([T/\tau]+1)\tau +\epsilon)$ which contains $[0, T]$. 

Next, for uniqueness, we suppose $\tilde \psi$ is another classical solution over $[0, T]$. Then by a very similar inductive argument and repeatedly applying the uniqueness result in Theorem \ref{thm: sec4_main_theorem}, $\tilde \psi$ must be given by the formula at \eqref{eqns: power_time_marching_representation}.  
\end{proof}

\subsubsection{The time-marching Monte Carlo algorithm}\label{subsubsec: time-marching MC algorithm}
To numerically approximate the global solution $\psi(t, \cdot)$ using the time marching representation \eqref{eqns: power_time_marching_representation}, there are two technical issues. First, except for  the initial condition $g$, we do not have the precise functions $P_{\tau}^k g$ for $k\geq 1$, but only Monte Carlo approximations from the local solver. Second, such Monte Carlo approximations are just defined over a specific spatial grid of $\R$ instead of a function from $\R$. The first problem is addressed via a stability result demonstrated in Lemma \ref{lemma: a local derivative bound and stability} for $L_\infty$ data which are not necessarily differentiable, while the second can be solved by including an interpolation mechanism, making use of the homogeneity of the Dirac branching system with respect to translations, which was established in Proposition \ref{prop: translation_homogeneity}.  We demonstrate the mechanism here and summarize it in Algorithm \ref{alg:global-monte-carlo-time_marching}. 

Fix a spatial boundary $[-L, L]$ defined by some interpolation radius $L>0$. Let $M$ be the number of grid points. Then  we can define the spatial increment  $\Delta x:= 2L/(M-1)$ and the set of grid points $(x_j)_{j=1}^M$ defined by $x_j:= -L + (j-1)\Delta x$ for all $j = 1, 2, \cdots, M$.

Let $m>0$ be the mass and $g\in C_b^1$ be the initial condition of the Soler model. Let $(E_t, \mu_t, M_t)$ be a Dirac branching system and multiplicative functionals such that their parameter set satisfies Assumption \ref{assumption: non-linear_coefficients}.  The first goal is to obtain approximations for $\psi(t, x_i, z)$ for all $i= 1, \cdots, M$ and $z = \pm 1$. 

Fix $N\in \mathbb{N}$. First,
we can sample $N$ independent copies of $(E_t, M_t, \mu_t)\mid \mu_0:= \delta_{0, z, 1}$ for every $z = \pm 1$, which we denote by:
\begin{align}
    E_t^{(i)}(z), \quad M_t^{(i)}(z), \quad \mu_{t}^{(i)}(z),\quad z= \pm 1, \quad i = 1, \cdots, N.
\end{align}
Motivated by the homogeneity to translation established in Proposition \ref{prop: translation_homogeneity}, we define the random operator $\hat P_{t, N, M}$ such that for every $g\in L_\infty $ (not necessarily differentiable) with $P_t g$ to be well-defined, we output the random vector $\hat P_{t, N, M}g \in \C^{[M]\times \{\pm 1\}}$ where $[M]:= \{1, \cdots, M\}$ defined by:
\begin{align}\label{eqns: definition_of_the_random_operator}
   \hat P_{t, N,M} g(j, z):=\frac{1}{N}\sum_{i=1}^N E_t^{(i)}( z) M_t^{(i)}( z) G_{\bg}(T_{x_j}(\mu_t^{(i)}(z)), \quad j = 1, \cdots, M,
\end{align}
where $T_x$ is the translation operator over $E$ for all $x\in \R$ defined back in \eqref{eqns: translation_operator_over_measure}. 

Next to extend the random vector $\hat P_{t, N, M}$ to a function on all of $(x, z)\in \R\times \{\pm 1\}$, we first define the notion of interpolation we use: for a vector $\vec{g}\in \C^{[M]\times \{\pm 1\}}$, we define $I_{L,M}[\vec{g}]:\R\times \{\pm 1\}\to \C$ by:
\begin{equation*}
I_{L, M}[\vec{g}](x, z) = 
\begin{cases}
    \vec{g}(j,z), \quad \text{ if } x= x_j, \\
    \vec{g}(1,z), \quad \text{ if }  x\leq -L, \\
    \vec{g}(M,z), \quad \text{ if }  x\geq L,\\
    \frac{\vec{g}(j+1,z)-\vec{g}(j,z)}{x_{j+1}-x_j} (x-x_j) + \vec{g}(j,z),\quad \text{ if }  x\in [x_j, x_{j+1}], 
\end{cases} \!\!\!\! z =\pm 1, \quad j = 1,\cdots, M.
\end{equation*}
With the notion of the interpolation, we then define the random function $\hat \psi_{N, M}(t, \cdot):\R\times \{\pm 1\}\to \C$ by:
\begin{align}
    \hat \psi_{N, M}(t,\cdot) = I_{L,M}[\hat P_{t, N, M}g](\cdot).
\end{align}

Now we are able to describe the time marching Monte Carlo algorithm we propose. 
To this end, fix $t>0$. Then for every $0<\tau <t$ (which we will specify how small it is) and $K:=[t/\tau]$, we let $\hat P_{\tau, N, M}^{(k)}$ with $k\in [K]$ be independent copies of the random operator $\hat P_{\tau, N, M}$, which in turn are defined by the independent families:
\begin{align}
    \{E_\tau^{(i,k)}(z), M_\tau^{(i,k)}(z), \mu_\tau^{(i,k)}(z)\}_{i=1, k = 1}^{N, K}.
\end{align}
If $t/\tau\notin \N$, then additionally,  we let $\hat P^{(K+1)}_{\tau, N, M}$ to be a copy of the independent operator $\hat P_{t - K\tau, N, M}$, which are defined by the random independent collection which is independent to everything defined before:
\begin{align}
    \{E_{t-K\tau}^{(i,K+1)}(z), M_{t-K\tau}^{(i,K+1)}(z), \mu_{t-K\tau}^{(i,K+1)}(z)\}_{i=1, }^{N}
\end{align}
Our global algorithm is then such that for every $g\in C_b^1$, we choose a small enough $\tau<t$ and define inductively the random sequence by $(\hat \psi_k)_{k=1}^K$ with $\hat \psi_0 := g$, and inductive relation
\begin{align}
    \hat \psi_k:= I_{L,M}[\hat P^{(k)}_{\tau, N, M} \hat\psi_{k-1}],\quad k\in [K].
\end{align}
We then output $\hat \psi_{K+1}$ (or $\hat \psi_K$ when $t/\tau\in \N$) to approximate $\psi(t, \cdot)$. For simplicity we consider only the case $t/\tau \in \N$ in the following pseudocode.
\begin{breakablealgorithm}
\caption{Global Monte Carlo solver for $\psi(t, \cdot)$ by time-marching}
\label{alg:global-monte-carlo-time_marching}
\begin{algorithmic}[1]
\Require Spatial interval $[-L, L]$, number of spatial point $M$, terminal time $t$,
number of samples $N$, initial condition $g$, the time marching parameter $\tau$
\Require $\tau<t$ to be small enough with $t/\tau \in \N$. 
\State Define $K:= [t/\tau]$. Set $\hat g_0\gets g$

\For{$k = 1, \cdots, K$}
\For{$z\in \{\pm 1\}$}
\State Use \textbf{Algorithm \ref{alg:sample_trajectories}} to sample $N$ trajectories starting from $\mu_0 = \delta_{0,z,1}$:
\begin{align}
    \bigg\{E_\tau^{(i,k)}(z), M_\tau^{(i, k)}(z), \mu^{(i, k)}(z)\bigg\}
\end{align}
\For{$j = 1, \cdots, M$}
\For{$i = 1, \cdots, N$}
\State Compute the product:
\begin{align}
    D_{\tau, j}^{(i, k)}(z) \gets E_\tau^{(i,k)}(z) M_\tau^{(i,k)}(z) G_{\pmb {\hat g_{k-1}}}\bigg(T_{x_j}\bigg(\mu_\tau^{(i,k)}(z)\bigg)\bigg)
\end{align}
\EndFor
\State Set the random vector:
\begin{align}
    \hat P_{\tau, N, M}^{(k)}(j, z) \gets \frac{1}{N}\sum_{i=1}^N D_{\tau,j}^{(i,k)}(z)
\end{align}
\EndFor
\State Set the function:
\begin{align}
    \hat \psi_{k,N, M} (\tau, x, z)\gets I_{L,M}[\hat P^{(k)}_{\tau, N, M}](x, z)
\end{align}
\EndFor
 \State Set:
 \begin{align}
     \hat g_k \gets \hat \psi_{k, N, M}
 \end{align}
\EndFor
\State \Return $\hat g_K$
\end{algorithmic}
\end{breakablealgorithm}

\subsection{Convergence analysis of Algorithm \ref{alg:global-monte-carlo-time_marching}}\label{subsec: convergence_analysis}

We start by setting up notations: for any $D\subset \R$, we write 
\begin{align}
    \norm{f}_{\infty, D}
\end{align}
for the supremum norm of some complex-valued function $f$ with the supremum taken over spatial coordinates $x\in D$ together with all other variables. In addition we also define for every $D\subset \R$ and $t>0$ that:
\begin{align}
    D^{(t)}:= \{x\in \R: d_\R(x, D)\leq t\},
\end{align}
the closed collection of points $t$ distance away from $D$ with respect to the standard absolute value distance $d_\R$.

The proof of the convergence analysis of Algorithm \ref{alg:global-monte-carlo-time_marching} will use the following Lemma, which is of independent interest and will be expanded to more results in Section \ref{subsec: global derivative bound and stability}.  
\begin{lemma}\label{lemma: a local derivative bound and stability}
Let $g, h\in L_\infty$ (not necessarily differentiable) and $D\subset \R$. For every $t>0$ with:
\begin{align}
   V(t,\max\{\norm{g}_{\infty, D^{(t)}}, \norm{h}_{\infty, D^{(t)}}\}) <\infty, 
\end{align}
the functions $P_t g, P_t h$ are defined on $D$ and we have:
\begin{align}\label{eqns: sec6. stability_local}
    \norm{P_tg - P_t h}_{\infty, D} \leq \p_B V(t,\max\{\norm{g}_{\infty, D^{(t)}}, \norm{h}_{\infty, D^{(t)}}\}) \norm{g-h}_{\infty, D^{(t)}}. 
\end{align}
\end{lemma}
\begin{proof}
The result follows from the fact that for every pair of finite lists $a:=(a_1, \cdots, a_N)\in \C^N$ and $b:=(b_1, \cdots, b_N)\in \C^N$ where $N\in \mathbb{N}$, we have:
\begin{align}
    \prod_{i=1}^N a_i - \prod_{i=1}^N b_i = \sum_{i=1}^N (a_i - b_i) \prod_{j<i}a_j \prod_{k>i}b_k,
\end{align}
by adding telescoping terms with one different term at a time and so:
\begin{align}
    \abs{  \prod_{i=1}^N a_i - \prod_{i=1}^N b_i} \leq  N \norm{a - b}_\infty(\max\{\norm{a}_\infty, \norm{b}_\infty\})^{N-1}.
\end{align}
 We consider the enlarged domain $D^{(t)}$ in both upper bounds because spatial parameters of particles in $(\mu_t)$ is at most $t$ distance away from the spatial coordinate of the initial particle.
\end{proof}

\begin{theorem}\label{thm: convergence_of_global_alg}
For every $t\in (0, L]$ and $g\in C_b^1$, there exists $t>\tau>0$ such that for every $\epsilon>0$,
\begin{align}
   \lim_{M\to \infty}\limsup_{N\to \infty}\P(\norm{\hat g_K - \psi(t, \cdot)}_{\infty, [-L+t, L-t]} \geq \epsilon) = 0
\end{align}    
\end{theorem}
\begin{proof}
The proof proceeds inductively along the time-marching steps and is deferred to Appendix \ref{appendix: sec5_convergence_Proof}. For its heuristic, at each step, the probability of the approximation error bounded away from zero can be decomposed into three parts: the spatial interpolation error, controlled by the uniform derivative bound \eqref{eqns: classical_solution_additional_boundedness_condition} in our definition  of a classical solution; the stability error, controlled by the stability estimate in \eqref{eqns: sec6. stability_local} of Lemma \ref{lemma: a local derivative bound and stability}; and the Monte Carlo error, controlled by a union bound together with the choice of the time step $\tau$, which is chosen sufficiently small so that the corresponding Monte Carlo estimator gives a finite variance. Iterating these estimates over the finitely many time steps up to time $t$ and sending $N\to \infty$ then $M\to \infty$ yields the convergence result. 

The restriction of the error estimate to $[-L+t, L-t]$ is due to the flow of the Dirac branching system underlying the probabilistic representation for $t$ units of time.
\end{proof}

\subsection{Global stability of the Soler model}\label{subsec: global derivative bound and stability}
The result in Lemma \ref{lemma: a local derivative bound and stability} was local in time and can be extended to a global stability result using the a priori estimate \ref{lemma: apriori_L_inf_estimate} and the time-marching representation in Theorem \ref{thm: global_well_posednesss_L_inf}. 
\begin{proposition}\label{prop: global_independent_results}
For any initial data $g\in C_b^1$, we let $\psi_g(t, x, z)$ be the corresponding unique global solution of the Soler model. Then for every $t>0$ the solution map:
\begin{align}
    g\mapsto \psi_g(t,\cdot)
\end{align}
is locally Lipschitz in the uniform norm, that is, for every $g\in C_b^1$, for every $R>0$, there exists a constant $C = C(t,g, R)$ such that for every $h_1, h_2\in C_b^1$ with $\norm{h_i-g}_\infty <R$ we have:
\begin{align}\label{sec6_global_stability}
    \norm{\psi_{h_1}(t,\cdot) - \psi_{h_2}(t, \cdot)}_\infty \leq C(t,g, R)\norm{h_1 - h_2}_\infty.
\end{align}
\end{proposition}
\begin{proof}
The proof follows from iterating Lemma \ref{lemma: a local derivative bound and stability} for $D = \R$ along the time-marching representation of Theorem \ref{thm: global_well_posednesss_L_inf}. The details are deferred to Appendix \ref{appendix: sec5_global_stability_proof}. 
\end{proof}

\subsection{Numerical experiments}\label{subsec: numerical experiments}
In this section, we illustrate the numerical performance of the time-marching Monte Carlo method introduced in Algorithm \ref{alg:global-monte-carlo-time_marching} on the Soler model \eqref{eqns: sec4_nonlinear_Dirac}, which we recall to have a mass parameter $m>0$. We consider three benchmark solutions of the Soler model for which explicit global solutions are available: a standing solitary wave, a traveling solitary wave and a solution with constant initial condition. These examples allow us to compare the Monte Carlo approximation directly with the exact solution, both through the shape of their spatial profiles and through the corresponding $L_\infty$ error over time. We first recall the three benchmark solutions and then present numerical results. 

\subsubsection{Benchmark solutions}\label{subsec: benchmark solutions}
Standing and traveling solitary waves of the one-dimensional Soler model are classical and have been extensively studied numerically; see, for instance \cite{AC81, XuShaoTang2013} and references therein. We recall below the explicit forms used in our numerical experiments. 
\paragraph{\bf Standing solitary wave solution}
 A standing solitary wave solution with frequency parameter $\omega\in (0, m)$ is a solution of the form:
\begin{align}
    \psi_{S, \omega}(t, x, z) = e^{-i\omega t} g_{S,\omega}(x, z),\quad t\geq 0, x\in \R, z= \pm 1,
\end{align}
where $g_{S, \omega}$ is the initial condition defined by:
\begin{align}\label{eqns: standing_initial_condition}
    g_{S,\omega}(x)  =
\begin{bmatrix}
     g_{S,\omega}(x, 1) \\
      g_{S,\omega}(x, -1)
\end{bmatrix} = H\begin{bmatrix}
    A(x)\\
    iB(x)
\end{bmatrix},
\end{align}
where
\begin{align}\label{eqns: definition of A, B}
    A(x) = \sqrt{2}\frac{\kappa (m+\omega)^{1/2} \cosh(\kappa x)}{m+\omega \cosh(2\kappa x)}, \quad B(x) = \sqrt{2}\frac{\kappa (m-\omega)^{1/2} \sinh(\kappa x)}{m+\omega \cosh(2\kappa x)},
\end{align}
with 
\begin{align}
    \kappa:= \sqrt{m^2 - \omega^2}, \quad H:=
\frac{1}{\sqrt{2}}
\begin{bmatrix}
    1 & 1 \\
    1 & -1
\end{bmatrix}.
\end{align}
The matrix $H$ is commonly called the Hadamard matrix. The supremum norm of the initial condition \eqref{eqns: standing_initial_condition} is given by 
\begin{align}
    \norm{\psi_{S, \omega}(0, \cdot)}_\infty &=
\begin{cases}
    \frac{m}{2\sqrt{\omega}},&\quad \text{ if } 0<\omega\leq m/2, \\
    \sqrt{m-\omega}, & \quad \text{ if } m/2\leq \omega <m\, . 
\end{cases}
\end{align}

\begin{remark}\label{remark: relation with AC81}
In the convention of \cite{AC81},  the initial standing-wave profile is written as 
$\begin{bmatrix}
    A(x)\\
    iB(x)
\end{bmatrix}$, defined by \eqref{eqns: definition of A, B}. It is because they considered the Soler system under the matrix pair (defined in \eqref{eqns: Clifford pair}) $(\sigma_3, i\sigma_2)$, which is different from our choice of $(\sigma_1, -i\sigma_2)$ (pointed out at \eqref{eqns: our_choice_matrix_pair}) by a unitary transformation via the Hadamard matrix $H$. Hence, the initial condition at \eqref{eqns: standing_initial_condition} is a result of  Theorem \ref{thm: well-posedness_different_matrix_pair} on the initial condition in \cite{AC81}. 
    
\end{remark}
\paragraph{\bf Traveling solitary wave solution} 
A traveling wave solution is parameterized by a velocity $v\in(-1, 1)$ and location $x_0\in \R$, in addition to a frequency parameter $\omega\in (0, m)$. It is given by 
\begin{align}\label{eqns: traveling_wave_profile}
    \psi_{T,\omega, v, x_0}(t, x,\cdot) = e^{-i\omega t'}H
\begin{bmatrix}
\sqrt{\frac{\gamma+1}{2}} & \sgn(v) \sqrt{\frac{\gamma-1}{2}}\\
\sgn(v) \sqrt{\frac{\gamma-1}{2}} & \sqrt{\frac{\gamma+1}{2}}
\end{bmatrix}
\begin{bmatrix}
    A(x') \\
    iB(x')
\end{bmatrix},
\end{align}
where 
\begin{align}
    \gamma :=\frac{1}{\sqrt{1 - v^2}},\quad x'= \gamma((x-x_0) - vt),\quad t' = \gamma(t - v(x-x_0)),
\end{align}
with $A, B$ defined in \eqref{eqns: definition of A, B} and $H$ being the Hadamard matrix. 
In particular, its initial condition is given by:
\begin{align}\label{eqns: traveling_initial_condition}
    \psi_{T,\omega, v, x_0}(0, x, \cdot) =e^{i\omega \gamma v(x-x_0)}
    H\begin{bmatrix}
\sqrt{\frac{\gamma+1}{2}} & \sgn(v) \sqrt{\frac{\gamma-1}{2}}\\
\sgn(v) \sqrt{\frac{\gamma-1}{2}} & \sqrt{\frac{\gamma+1}{2}}
\end{bmatrix}
\begin{bmatrix}
    A(\gamma(x-x_0))\\
    iB(\gamma(x-x_0))
\end{bmatrix}.
\end{align}
The supremum norm of the initial conditions of the solutions \eqref{eqns: traveling_initial_condition} satisfies:
\begin{align}
\norm{\psi_{T, \omega, v, x_0}(0,\cdot)}_\infty &= \sqrt{\gamma(1+\abs{v})}   \norm{\psi_{S, \omega}(0, \cdot)}_\infty\\
&=\sqrt{\gamma(1+\abs{v})}\begin{cases}
    \frac{m}{2\sqrt{\omega}},&\quad \text{ if } 0<\omega\leq m/2\\
    \sqrt{m-\omega}, & \quad \text{ if } m/2\leq \omega <m
\end{cases}.
\end{align}
We remark that the relation between the traveling wave profile at \eqref{eqns: traveling_wave_profile} and that in \cite{AC81} is off by a transformation via $H$ for the same reason pointed out in Remark \ref{remark: relation with AC81}.

\paragraph{\bf Constant initial condition}
In contrast with the solitary-wave solutions above, constant solutions do not appear to be commonly used as benchmark solutions in the numerical literature on the Soler model. We include this example because it provides an explicit non-decaying solution and, in particular, an initial condition in $C_b^1$ which is not in $L_2$. Let $c\in \C$. By explicit calculation, one can verify that if we start from the constant initial condition:
\begin{align}
    g_{c}(x, z):= c, \quad (x, z)\in \R\times \{\pm 1\},
\end{align}
then the unique solution to the Soler system at \eqref{eqns: sec4_nonlinear_Dirac} is given by:
\begin{align} \label{eqns: constant_global_solution_soler}
    \psi_c(t, x, z):= e^{i (2\abs{c}^2 - m)t}c.
\end{align}
\begin{remark}
    To see how one can obtain the formula \eqref{eqns: constant_global_solution_soler} without knowing it a priori, indeed, by the homogeneity of the Dirac branching system (Proposition \ref{prop: translation_homogeneity}) and the time-marching representation at Theorem \ref{thm: global_well_posednesss_L_inf}, one can infer that $\psi_c(t, x, z) = q_c(t)$ for some functions $q_c$ depending only on time. Substituting it into \eqref{eqns: sec4_nonlinear_Dirac} gives:
    \begin{align}
        q_c'(t) = i(2\abs{q_c(t)}^2 - m)q_c(t), \quad q_c(0) = c,
    \end{align}
    which gives the desired formula.

\end{remark}

\subsubsection{Numerical solutions}
Now, we apply Algorithm \ref{alg:global-monte-carlo-time_marching} to the three benchmark solutions introduced in  Section \ref{subsec: benchmark solutions} (which are standing solitary waves, traveling solitary waves and constant profiles). Throughout the experiments, we consider Soler systems with mass parameter $m = 1$ and compute its solution up to $T = 30$. For each example, we will calculate the initial supremum norm $\norm{\psi(0, \cdot)}_\infty $ and the two upper bounds $t^*_{1st}(B)$ and $t_{var}^*(B)$ as a function of $B>0$, required in the one-step Monte Carlo solver at Algorithm \ref{alg:local_MC_solver}. While $t^*_{1st}$ is given by \eqref{eqns: upper_bound_of_T_local_wellposedness} with:
\begin{align}
    t^*_{1st}(B) := \frac{1}{2} \log (1+\frac{1}{2B^2}),
\end{align}
the second upper bound given by:
\begin{align}
    t_{var}^*(B):= \sup\{s>0: \inf_\Lambda V_s(\Lambda)<\infty\},
\end{align}
which we solve numerically using the explicit formula of $V_s(\Lambda)$ at \eqref{eqns: second_moment_upper_bound}  as well as the easily deduced fact that the unique optimizer of $V_s(\Lambda)$ has the form 
\begin{align}
    (\lambda_f, \lambda_{b_1}, \lambda_{b_2})= (m =1, r, r),
\end{align}
 for some $r>0$. We specify the time-marching parameter $\tau>0$ to be such that $\tau< t^*_{1st}(B), t_{var}^*(B)$ with $B = \norm{\psi(0, \cdot)}_\infty$ and the event-rate parameters $\Lambda$ to be the unique optimizer of $\inf_\Lambda V_\tau(\Lambda)$ obtained numerically. These parameters determine the underlying Dirac branching system and multiplicative functionals. We note that such parameter choice has the best effect to minimize the variance of the Monte Carlo estimator in the case of the constant profile.

Unless otherwise specified, we use an interpolation radius of $L:= 40$ over $M:= 401$ grid points and $N:= 10^5$ Monte Carlo samples at each time step. In addition, different seeds are used in different experiments.  As a remark, the upper bound $t^*(L(T, \norm{\psi(0,\cdot)}_\infty))$ given in \eqref{eqns: theoretical_upper_bound_for_tau} is very small and is impractical for numerical experiments, so we will not use it. After that we will report the error:
\begin{align}\label{eqns: definition_of_uniform_error}
    \mathcal{E}_\infty (t) := \norm{\widehat \psi(t, \cdot) - \psi(t, \cdot)}_{\infty, (-L+t, L-t)},\quad t\in [0, T],
\end{align}
where $\widehat \psi$ is the output of our algorithm and $\psi$ is the benchmark solution, on the interior region appearing in the convergence analysis in Theorem \ref{thm: convergence_of_global_alg}. In addition, we present spatial profile of the output of our algorithm and the benchmark solution alongside each other at several times.

\begin{figure}[p]
    \centering
    \begin{minipage}{0.6\textwidth}
        \centering
        \includegraphics[width=\textwidth]{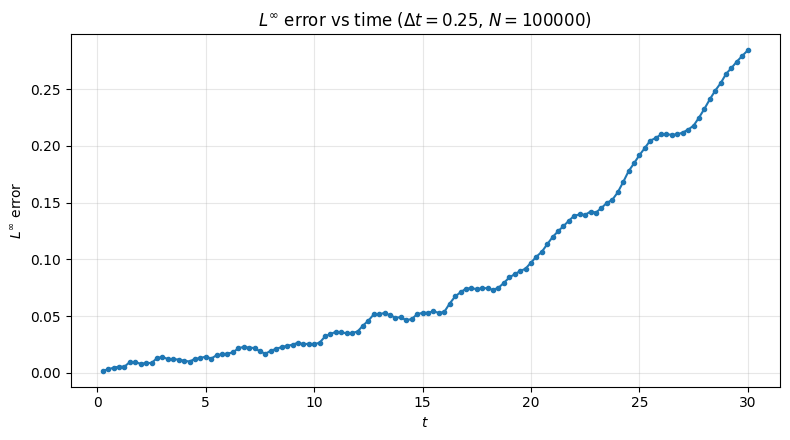}

        \smallskip
        (a). $L^\infty$ error over time.
    \end{minipage}

    \vspace{0.4cm}

    \begin{minipage}{\textwidth}
        \centering
        \includegraphics[width=0.9\textwidth]{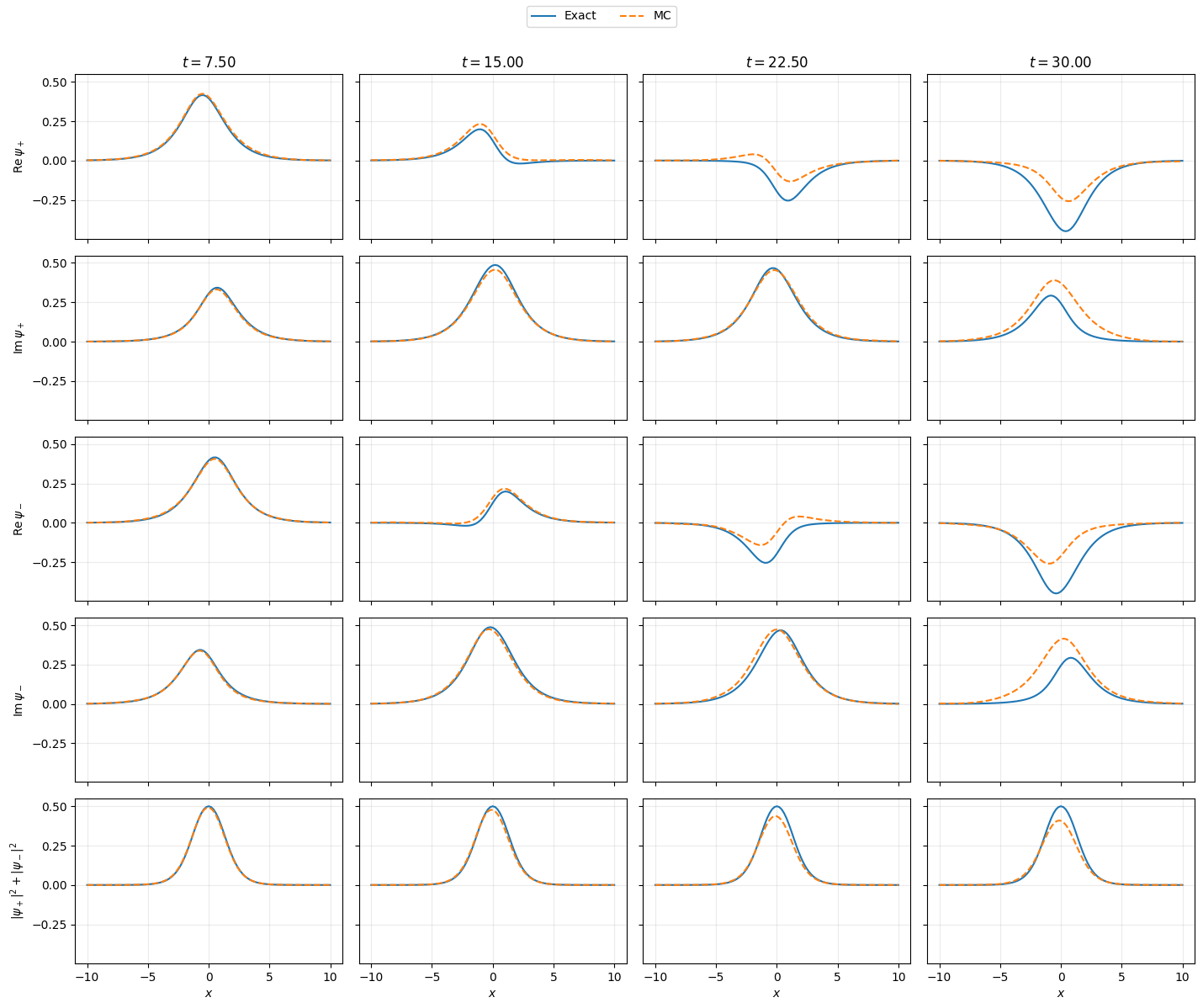}

        \smallskip
        (b). Spatial profiles at $t=7.5,15,22.5,30$.
    \end{minipage}

\caption{
Numerical results for the standing solitary wave with $\omega=0.75$.  A time-marching parameter of $\tau = 0.25$ is used. 
The localized charge profile (the last row of $(b)$) remains centered at the origin throughout the
evolution.
}
    \label{fig:standing_wave_experiment}
\end{figure}

\paragraph{\bf Standing solitary wave solution}
First we consider standing solitary wave solution with frequency parameter $\omega = 0.75$, which yields an initial supremum norm of:
\begin{align}
 \norm{g_{S, \omega = 0.75}}_\infty  = 0.5,
\end{align}
which gives:
\begin{align}
t^*_{1st} \approx 0.5493, \quad t^*_{var} \approx 0.4417.
\end{align}
We choose a time-marching parameter of $\tau = 0.25$, which gives $120$ timesteps with a Dirac branching parameter of:
\begin{align}
    \lambda \approx 1.7923, \quad p(\mathsf{f}) \approx 0.5579, \quad p(\mathsf{b_1}) =p( \mathsf{b_2}) \approx 0.2210
\end{align}
The resulting second moment upper bound at \eqref{eqns: second_moment_upper_bound} is given by:
\begin{align}
 V_\tau(\Lambda) \approx 0.62014
\end{align}
We plot our numerical results in Figure \ref{fig:standing_wave_experiment}.

\paragraph{\bf Traveling solitary wave solution}

\begin{figure}[p]
    \centering
    \begin{minipage}{0.6\textwidth}
        \centering
        \includegraphics[width=\textwidth]{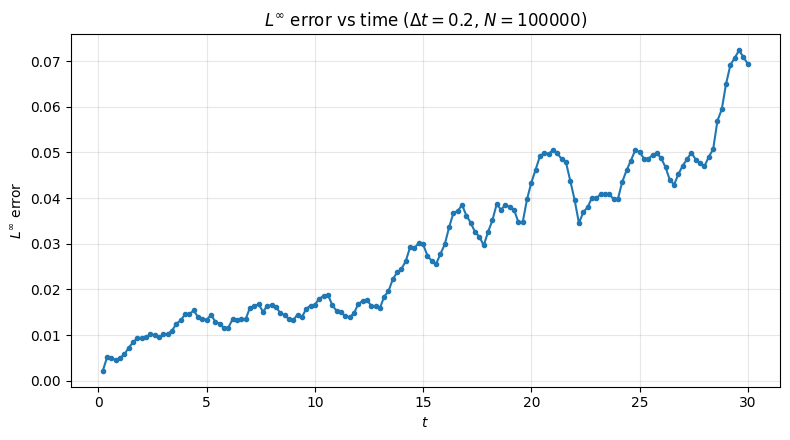}
        \smallskip
        (a). $L^\infty$ error over time.
    \end{minipage}

    \vspace{0.4cm}

    \begin{minipage}{\textwidth}
        \centering
        \includegraphics[width=0.9\textwidth]{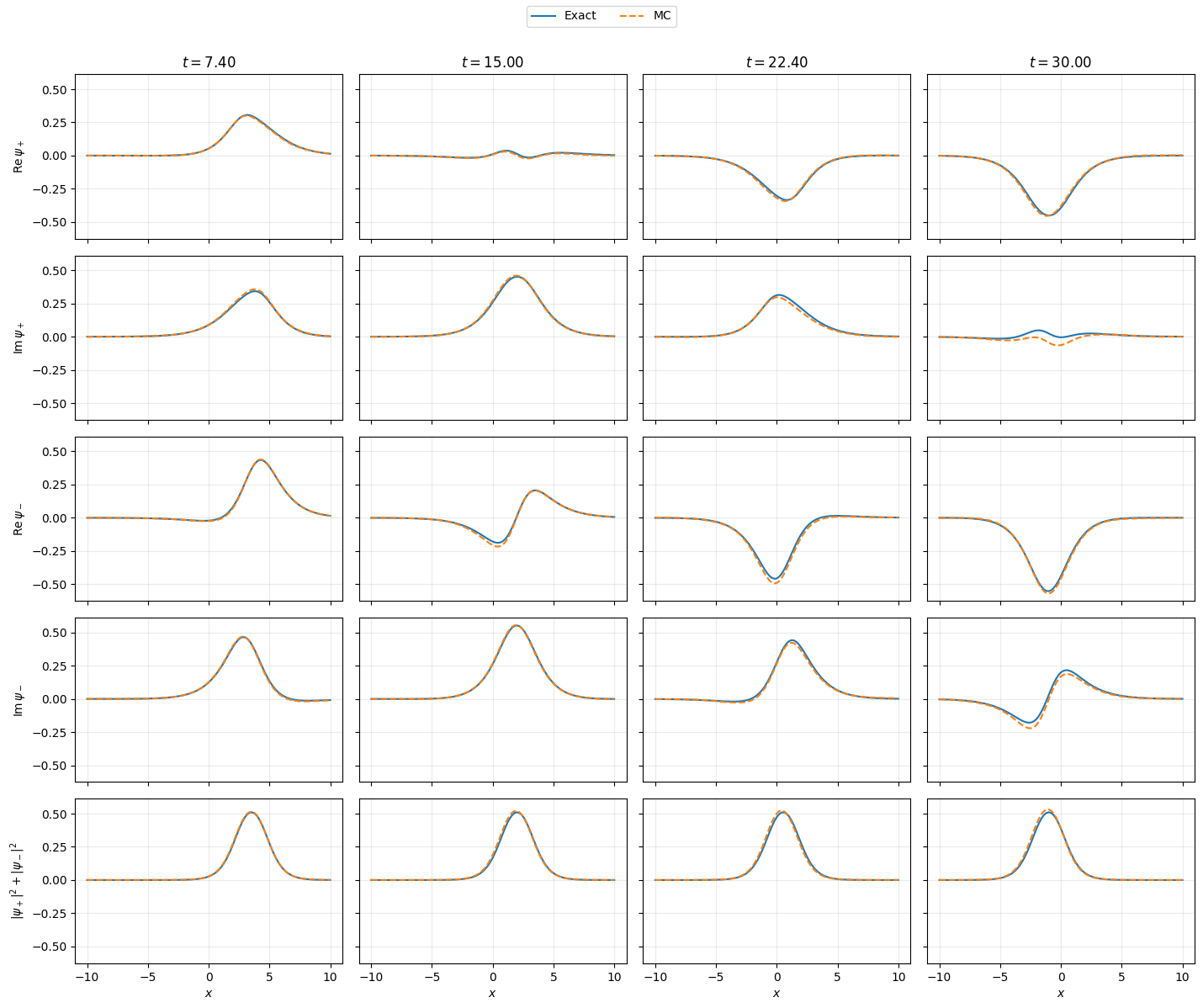}
        
        \smallskip
        (b). Spatial profiles at $t=7.4,15,22.4,30$.
    \end{minipage}

\caption{
Numerical results for the traveling solitary wave with
$\omega=0.75$, $x_0=5$, and $v=-0.2$. A time-marching parameter of $\tau = 0.2$ is used, which explains the reported times of $t=7.4, 22.4$ instead of $t = 7.5, 22.5$. The localized charge profile (the last row of (b)) propagates to the left while retaining its
shape with the change of center from $x = 5$ to $x = -1$. 
}
    \label{fig:traveling_wave_experiment}
\end{figure}

Next we consider traveling solitary wave solution with initial position $x_0 = 5$,  speed $v =-0.2$ and frequency parameter $\omega = 0.75$ which yields an initial supremum norm of:
\begin{align}
 \norm{g_{T, \omega = 0.75, v = -0.2, x_0 = 5}}_\infty  \approx 0.55334,
\end{align}
which gives:
\begin{align}
t^*_{1st} \approx 0.4841, \quad t^*_{var} \approx 0.3891.
\end{align}
We choose a time-marching parameter of $\tau = 0.2$, which gives $150$ timesteps with a Dirac branching parameter of:
\begin{align}
    \lambda \approx 1.9009, \quad p(\mathsf{f}) \approx 0.5261, \quad p(\mathsf{b_1}) =p( \mathsf{b_2}) \approx 0.2370
\end{align}
The resulting second moment upper bound at \eqref{eqns: second_moment_upper_bound} is given by:
\begin{align}
    \inf_\Lambda V_\tau(\Lambda) \approx 0.66057
\end{align}
We plot our numerical results in Figure \ref{fig:traveling_wave_experiment}.

\begin{figure}[p]
    \centering
    \begin{minipage}{0.6\textwidth}
        \centering
        \includegraphics[width=\textwidth]{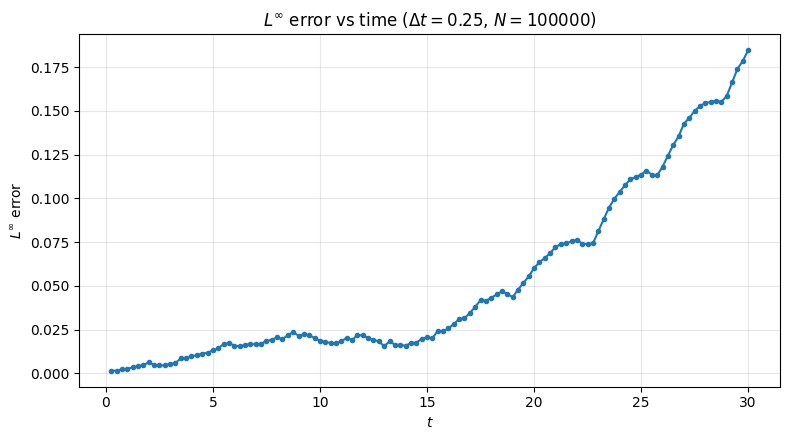}
        \smallskip
        (a). $L^\infty$ error over time.
    \end{minipage}

    \vspace{0.4cm}

    \begin{minipage}{\textwidth}
        \centering
        \includegraphics[width=0.9\textwidth]{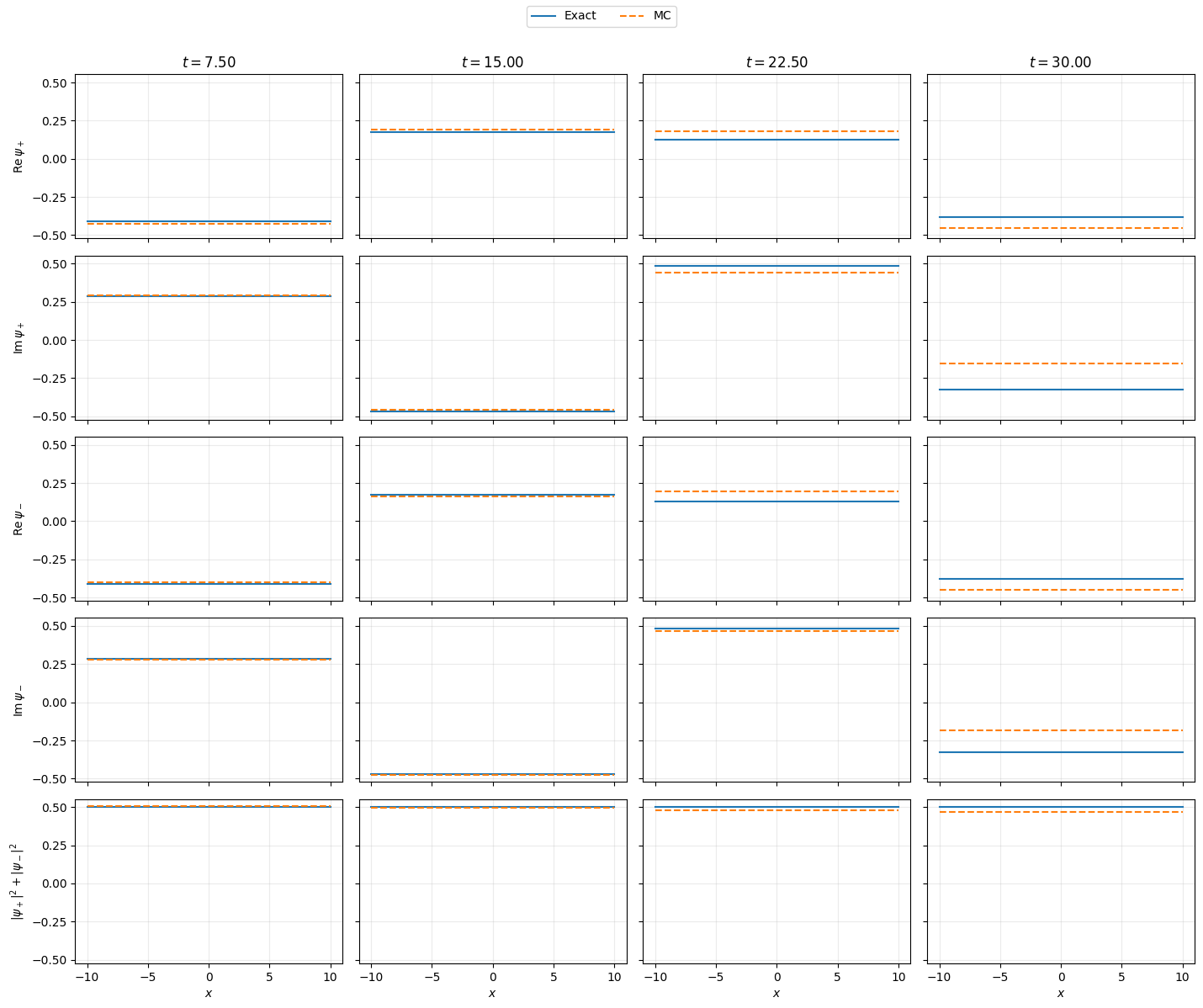}
        
        \smallskip
        (b). Spatial profiles at $t=7.5,15,22.5,30$.
    \end{minipage}

\caption{
Numerical results for the constant benchmark with $c=0.5$ with time marching-parameter $\tau = 0.25$. Constant profiles are natural in our $C_b^1$ problem data framework. 
}
    \label{fig:constant_experiment}
\end{figure}

\paragraph{\bf Constant initial condition}
Lastly, we consider an initial constant profile of $c = 0.5$, which shares the same $t^*_{1st}$ and $t^*_{var}$ as the standing solitary wave case. We choose the same time-marching parameter as well as Dirac branching parameter as in such case.
The plots are shown in Figure \ref{fig:constant_experiment}. 

In each figure, panel (a) shows the uniform error $\mathcal{E}_\infty(t)$ defined in \eqref{eqns: definition_of_uniform_error} over the interior interval $[-L+t, L-t]$ for all $t\in [0, 30]$. Panel (b) compares the Monte Carlo approximation with the corresponding exact solution at $t = 7.5, 15, 22.5, 30$ on the common spatial interval $[-10, 10]$. The solid and dashed curves represent the exact and the Monte Carlo solutions respectively. 
From top to bottom, the rows in panel (b) display $\Re \psi_+, \Im \psi_+, \Re \psi_-, \Im \psi_- $ and $\abs{\psi_+}^2+\abs{\psi_-}^2 $ with the notation:
\begin{align}
    \psi_\pm (t, x):= \psi(t, x, \pm 1).
\end{align}
The quantity 
\begin{align}
    Q(t):=\int_{-\infty}^\infty  (\abs{\psi_+(t, x)}^2+\abs{\psi_-(t, x)}^2)dx, \quad t\geq 0
\end{align}
is  called the total charge of the solution. For initial data with $Q(0)<\infty$, conservation of the total charge, that is, $Q(t) = Q(0)$ for all $t\geq 0$, is a standard result (see for instance \cite{AC81}). In fact, by a standard integration argument using \eqref{eqns: local_conservation_identity}, this conservation law extends to the global classical solution established in Theorem \ref{thm: global_well_posednesss_L_inf} for any initial data  $g\in C_b^1$, including the case when $Q(0) = \infty$. 

\clearpage


\section{Conclusion}\label{sec. conclusion}
\subsection{Summary} Using the stochastic calculus of PDMPs on general state spaces, we establish a probabilistic representation for nonlinear Dirac equations with $C_b^1$ coefficients and $C_b^1$ initial condition (Theorem \ref{thm: general_cubic_main_theorem}) where the driving noise of the representing process, which we term the Dirac branching system, is a marked point process. Since there are additional deterministic movements between events, we can make use of exponential functionals of the Dirac branching system; the result is a classification of admissible time of the stochastic representation illustrated in Corollary \ref{corollary: classification_of_stochastic_representation}. In the special case of the Soler model, which falls into the finite-time representation regime, we indeed have a sharp upper bound on the admissible time by considering constant initial profiles which are available in our $C_b^1$ regime (Corollary \ref{corollary: sharpeness results_constant_data}). In Section \ref{sec: global_section}, we then document the interesting "local probabilistic representation, but global well-posedness" phenomena of the Soler model by extending the one-step probabilistic representation to a time-marching representation in Theorem \ref{thm: global_well_posednesss_L_inf}. The key ingredient is a modification of an $L_\infty$ a priori estimate of \cite{Huh2013GrossNeveu} derived analytically. Alongside PDE applications, we develop the Monte Carlo method to implement both the one-step probabilistic representation (Algorithm \ref{alg:local_MC_solver}) and the time-marching representation (Algorithm \ref{alg:global-monte-carlo-time_marching}) with rigorous convergence analysis which makes use of stability results derived from the probabilistic representations. Lastly, numerical experiments are demonstrated in Section \ref{subsec: numerical experiments}. We remark the preliminary observations that the choice of the event rates of the Dirac branching system can have a variance reduction effect on the Monte Carlo implementation while not affecting the validity of the corresponding stochastic representation. 

\subsection{Outlook}
Of immediate interest is the extension of the current probabilistic representation of the nonlinear Dirac equations from $1+1$ to $3+1$ space-time dimension, the physically standard setting \cite{Dirac1928, Thaller1992}. From a PDE perspective, it would be interesting to relax the current  $C_b^1$ assumptions  to bounded Lipschitz initial conditions and coefficients, or even bounded measurable coefficients. The latter appears particularly promising, as most of our arguments carry over except for the application of Proposition \ref{prop: measure_generator}. From a numerical perspective, further variance reduction techniques and complexity analysis for both the one-step and multi-step algorithm would be of interest.

A particularly intriguing problem arising in the current paper is the lack of a probabilistic interpretation of the $L_\infty$ a priori estimate modified from \cite{Huh2013GrossNeveu}. While it is the key to open our global analysis of the Soler system from Section \ref{subsec: global_analysis} onward, it is derived in a geometric argument together with an algebraic manipulation to the specific cubic nonlinearity, independent of the probabilistic flavor of the rest of the paper. It is very interesting to see whether cancelations in complex-weighted expectations can establish $L_\infty$ estimates with a milder time dependence than the corresponding absolute-moment bounds. Furthermore, our framework follows a Kolmogorov backward-equation approach on the infinite-dimensional space of measures. The corresponding forward-equation viewpoint could potentially provide new connections to stochastic mechanics for nonlinear Dirac equations \cite{deAngelis1986}. In particular, it would be interesting to explore a forward-equation-type formulation at the level of the complex-valued solution itself, rather than only its modulus-squared density.

Lastly, motivated by the successful applications of probabilistic representations to the asymptotic analysis of nonlinear PDEs (\cite{McKean1975, Ratanov2006}), it would be interesting to investigate whether our probabilistic framework can provide further insights into limiting regimes of nonlinear Dirac equations such as  the non-relativistic regime \cite{Najman1992}.

\newpage 
\appendix
\section{Proof details of intermediate results in Section \ref{sec: nonlinear_case}}\label{appendix: main_proofs_inSec3}
This appendix collects the proof details deferred from Section \ref{sec: nonlinear_case}
\subsection{Proof of Markov property of the intervening chain in Proposition \ref{prop: weak_branching_of_dbs}} \label{appendix: sec3_branching}
 Indeed, given a current state of $(\hat \tau_n, \hat \rho_n, \hat \bx_n, \hat e_n)$ the interarrival time is sampled from the minimum of independent exponential clocks with rate $\lambda \abs{\mu_t^1}$ and $\lambda \abs{\mu_t^2}$, which depends only on $\abs{\hat \rho_n}$. After that there is a probability of $\abs{\mu_{\hat \tau_n}^j}/ \abs{\widehat\rho_n}$ for a particle to be chosen from  $\mu_{\hat \tau_n}^j$ for $j = 1, 2$, depending on the clock that goes off, and uniformly within that. An event is then sampled according to the same Dirac branching parameter. Hence, the transition depends only on the current state and it is straightforward to show that such transition is given by the lifted local characteristic at \eqref{eqns: definition_lifted_local_characteristics}.
 
\subsection{Proof of Proposition \ref{prop: property of V}}\label{appendix: Sec3_PropertyV}
We recall from \eqref{eqns: definition_of_V} that we have:
\begin{align}
    V(t, B; \eta ; \bbeta) = \E\bigg[\exp\bigg(\eta \int_0^t \abs{\mu_s}ds\bigg) M_t ^{\bbeta} B^{\abs{\mu_t}} \mid \abs{\mu_0} = 1\bigg)\bigg].
\end{align}
Consider the accumulated number of branching $(N_b(t))$ defined at \eqref{eqns: pure birth process}. Note that this is a pure birth process with linear self-exciting rate $\lambda_b(2N_b(t) +1)$. By standard calculations, see for instance P. 235 of \cite{ross1996stochastic} for the case of the Yule process, $N_b(t)$ is distributed as a negative binomial distribution with shape parameter $1/2$ and success probability $e^{-2\lambda_b t}$ with the following probability mass function:
\begin{align}\label{eqns: negative_binomial}
    \P(N_b(t) = k) = \binom{2k}{k} e^{-\lambda_b t} (\frac{1 - e^{-2\lambda_b t}}{4})^k,\quad k = 0, 1, 2, \cdots.
\end{align}
In addition, if we let $\bT:= (T_n)_{n=1}^\infty $ be the sequence of branching times and $\bT_k:= (T_n)_{n=1}^k $ denote the sequence of the first $k$ branching times where $k = 1, 2, \cdots$. The conditional distribution $\bT_k\mid N_b(t) = k$ is given by the order statistics of $k$ i.i.d. random variables supported over $[0, t]$ with density given by $f(x)= 2\lambda_b e^{2\lambda_b x}(e^{2\lambda_b t}-1)^{-1}$ for $x\in [0, t]$. More concretely, $\bT_k\mid N_b(t) = k$ has the joint density given by:
\begin{align}\label{eqns: appendix_order_statistics}
    f(t_1, \cdots t_k\mid N_b(t) = k) = k!(2\lambda_b)^ke^{2\lambda_b (t_1+\cdots + t_k)}(e^{2\lambda_b t}-1)^{-k}, \quad 0\leq t_1\leq \cdots \leq t_k\leq t.
\end{align}
Observe that, we have:
\begin{align}\label{eqns: V_proof_initial_decomposition}
    &V(t, B; \eta; \bbeta) = \E\bigg\{ \E\bigg[\exp\bigg(\eta \int_0^t \abs{\mu_s}ds\bigg) M_t ^{\bbeta} B^{\abs{\mu_t}} \mid \abs{\mu_0} = 1\bigg)\bigg]\mid \bT, N_b(t)\bigg\}\\
    &= \sum_{k=0}^\infty \E\bigg\{ \E\bigg[\exp\bigg(\eta \int_0^t \abs{\mu_s}ds\bigg) M_t ^{\bbeta} B^{\abs{\mu_t}} \mid \abs{\mu_0} = 1\bigg)\bigg]\mid \bT_k, N_b(t) = k\bigg\} \P(N_b(t) = k).
\end{align}
Now given $N_b(t) = k$ and the branching time $\bT_k$,  we denote $(e_1, \cdots, e_k)$ to be the jump event of the Dirac branching system at time $T_1, \cdots, T_k$ respectively and denote $N_f(j)$ to be the number of flips of the Dirac branching system over $(T_j, T_{j+1})$ for $j = 0, \cdots, k$ with $T_{k+1}:= t$. By the conditional independence of all $(e_j)_{j=1}^k$ and $(N_f(j))_{j=0}^k$, we have:
\begin{align}
    \E[M_t^{\bbeta} \mid \bT_k , N_b(t) = k] = \prod_{j=1}^k\E[\beta_{e_j}\mid \bT_k , N_b(t) = k] \prod_{j=0}^k\E[\beta_f^{N_f(j)}\mid \bT_k , N_b(t) = k].
\end{align}
Using the fact that $(e_1, \cdots, e_k)\mid \bT_k, N_b(t) = k$ are i.i.d. with each marginal having the outcome of $\mathsf{b}_j$ with probability ${p(\mathsf{b_j})}/({1 - p(\mathsf{f})})$ for $ j = 1, 2$, and that $N_f(j)\mid \bT_k, N_b(t) = k$ is a Poisson random variable with rate $(2j+1)\lambda_f(T_{j+1} - T_j)$, we have 
\begin{align}
      \E[M_t^{\bbeta} \mid \bT_k , N_b(t) = k] =  \bar{\beta_b}^{k}\exp\bigg((\beta_f - 1)\lambda_f\sum_{j=0}^k (2j+1)(T_{j+1} - T_j)\bigg)
\end{align}
It follows that 
\begin{align}
     & \E\Big[\exp \bigg(\eta \int_0^t \abs{\mu_s}ds\bigg)M_t^{\bbeta} B^{\abs{\mu_t}} \mid \bT_k , N_b(t) = k\Big] \\
      &= \exp\bigg(\eta\sum_{j=0}^k (2j+1) (T_{j+1} - T_j)\bigg) \bar{\beta_b}^{k}\exp\bigg((\beta_f - 1)\lambda_f\sum_{j=0}^k (2j+1)(T_{j+1} - T_j)\bigg) B^{2k+1}.
\end{align}
Taking the expectation over $\bT_k$ and using \eqref{eqns: appendix_order_statistics}, we have 
\begin{align}\label{eqns: step_in_V_proof_1}
 B \Big(\frac{2\lambda_b\bar \beta_bB^2}{e^{2\lambda_b t}-1}\Big)^k   k!\int_{0\leq t_1\leq \cdots \leq t_k\leq t} \exp \Big( 2\lambda_b \sum_{j=1}^k t_j + (\eta + \lambda_f(\beta_f -1))\sum_{j=0}^k(2j+1)(t_{j+1} - t_j) \Big)  dt_1\cdots t_k,
\end{align}
in which 
\begin{equation*}
    2\lambda_b \sum_{j=1}^k t_j + (\eta + \lambda_f(\beta_f -1))\sum_{j=0}^k(2j+1)(t_{j+1} - t_j) = [C(\eta, \beta_f) +\lambda_b](2k+1)t -2 C(\eta, \beta_f)\sum_{j=1}^k t_j.
\end{equation*}
Therefore, the integrand in \eqref{eqns: step_in_V_proof_1} is symmetric in $t_1, \cdots, t_k$, hence, we can further simplify \eqref{eqns: step_in_V_proof_1} to be:
\begin{align}
    Be^{C(\eta, \beta_f)t+\lambda_b t} \Big(\frac{2\lambda_b\bar \beta_bB^2}{e^{2\lambda_b t}-1}e^{2tC(\eta, \beta_f)}e^{2\lambda_b t}\int_0^t e^{-2C(\eta, \beta_f)s}ds \Big)^k.
\end{align}
It follows from the initial decomposition at \eqref{eqns: V_proof_initial_decomposition} and the negative binomial distribution \eqref{eqns: negative_binomial} that 
\begin{align}
    V(t, B; \eta, \bbeta) = \sum_{k=0}^\infty \binom{2k}{k} r^k e^{-\lambda_b t}  Be^{C(\eta, \beta_f)t+\lambda_b t},
\end{align}
where 
\begin{align}
    r &= \frac{2\lambda_b\bar \beta_bB^2}{e^{2\lambda_b t}-1}e^{2tC(\eta, \beta_f)}e^{2\lambda_b t}\int_0^t e^{-2C(\eta, \beta_f)s}ds\frac{1-e^{-2\lambda_b t}}{4} \\
    &=\frac{\lambda_b \bar \beta_b B^2}{4C(\eta, \beta_f)}(e^{2tC(\eta, \beta_f)} -1),
\end{align}
which is precisely the definition of $r$ in the definition of \eqref{eqns: definition of r}. The equivalent condition for finiteness follows from the elementary fact that $\sum_{k=0}^\infty \binom{2k}{k}r^k$ converges and is equal to $\frac{1}{\sqrt{1 - 4r}}$ if and only if $r<1/4$.

\subsection{Proof of Differentiability in Proposition \ref{prop: property of semigroup of dbs}}\label{appendix: Sec3_Differentiability}
We want to prove that under the additional differentiability assumption $\fh, \fu, \fm \in C_b^1$, we have $\bP_t \fh \in C_b^1$. By the homogeneity with respect to translations from Proposition \ref{prop: translation_homogeneity}, we have:
\begin{align}
    \bP_t \fh(\bx) = \E[E_t^{x_\bx}M_t^{x_\bx} G_\fh(T_{x_\bx} \mu_t)\mid \mu_0 = \delta_{0, z_\bx, \theta_\bx}],
\end{align}
where:
\begin{align}
    E_t^{x_\bx} &= \exp(\int_0^t \sum_{\by\in \mu_s} \fu(x_\by + x_\bx, z_\by, \theta_{\by})ds).
\end{align}
Differentiating against $x_\bx$ (which we write $\p_x$) gives:
\begin{align}
    \p_x E_t^{x_\bx}= E_t^{x_\bx} \int_0^t \sum_{\by\in \mu_s} \p_x \fu(x_y + x_\bx, z_\by, \theta_\by)ds,
\end{align}
with the upper bound for the absolute value:
\begin{align}\label{eqns: upper_bound_for_Edx}
    \abs{ \p_x E_t^{x_\bx}} \leq   \exp(\sup_{\X}\Re \fu \int_0^t \abs{\mu_s}ds)\norm{\p_x \fu}_\infty \int_0^t \abs{\mu_s} ds.
\end{align}
In addition, if we denote $(\tau_n,\bx_n, e_n)_{n=1}^\infty$ the marked jump sequence of the Dirac branching system, then we can write:
\begin{align}
    M_t^{x_\bx} = \prod_{\tau_n\leq t} \fm (e_n, (x_{\bx_n}+x_\bx, z_{\bx_n}, \theta_{\bx_n})).
\end{align}
Differentiating against $x_\bx$ gives:
\begin{align}
    \p_x M_t^{x_\bx} = \sum_{e\in J} \sum_{\tau_n\leq t,e_n = e} \p_x \fm (e_n, (x_{\bx_n}+x_\bx, z_{\bx_n}, \theta_{\bx_n}))&\times \prod_{\tau_m\leq t, e_m\neq e} \fm (e_m, (x_{\bx_m}+x_\bx, z_{\bx_m}, \theta_{\bx_m}))\\
     &\times  \!\!\!\! \prod_{\tau_m \leq t,e_m = e,m\neq n} \!\!\!\!\!\fm (e_m, (x_{\bx_m}+x_\bx, z_{\bx_m}, \theta_{\bx_m})).
\end{align}
It follows that if we define:
\begin{align}
    N_e(t):= \widehat N([0, t]\times \X\times \{e\}), \quad e\in J,
\end{align}
then we have the absolute value upper bound:
\begin{align}\label{eqns: upper_bound_for_Mdx}
    \abs{\p_x M_t^{x_\bx}} \leq \sum_{e\in J} \norm{\p_x \fm( e, \cdot )}_\infty  N_e(t) m_e^{N_e(t) - 1} \prod_{e'\neq e} m_{e'}^{N_{e'}(t)} 
\end{align}
Lastly:
\begin{align}
    G_\fh (T_{x_\bx}\mu_t) = \prod_{\by\in \mu_t} \fh(x_\by + x_\bx, z_\by, \theta_\by),
\end{align}
with spatial derivative:
\begin{align}
    \p_x  G_\fh (T_{x_\bx}\mu_t) = \sum_{\by\in \mu_t} \p_x \fh(x_\by + x_\bx, z_\by, \theta_\by) \prod_{\by' \in \mu_t - \delta_{\by}}\fh(x_{\by'} + x_\bx, z_{\by'}, \theta_{\by'}).
\end{align}
Taking absolute value gives the upper bound:
\begin{align}\label{eqns: upper_bound_for_Gdx}
    \abs{\p_x  G_\fh (T_{x_\bx}\mu_t)} \leq \norm{\p_x \fh}_\infty \abs{\mu_t} \norm{\fh}_\infty^{\abs{\mu_t}-1}.
\end{align}
Applying the product rule and  using the upper bounds \eqref{eqns: upper_bound_for_Edx}, \eqref{eqns: upper_bound_for_Mdx} and \eqref{eqns: upper_bound_for_Gdx}, we have:
\begin{equation}\label{eqns: final_upper_bound}
    \begin{split}
   & \E[\abs{\p_x (E_t^{x_\bx} M_t^{x_\bx} G_\fh(T_{x_\bx}\mu_t))}]  \\
    \leq & \norm{\p_x \fu}_\infty \p_\eta V(t,\norm{\fh}_\infty;  \sup_{\X} \Re \fu; (m_e)_{e\in J}) + \sum_{e\in J}  \norm{\p_x \fm (e,\cdot)}_\infty\p_{\beta_e}V(t,\norm{\fh}_\infty; \sup_{\X} \Re \fu; (m_e)_{e\in J})   \\
    &\quad +  \norm{\p_x \fh}_\infty \p_B  V(t,\norm{\fh}_\infty; \sup_{\X} \Re \fu; (m_e)_{e\in J}),
\end{split}
\end{equation}
which is finite because $ V(t,\norm{\fh}_\infty;  \sup_{\X}\Re \fu; (m_e)_{e\in J})$ is finite (as we assumed $t<t^*(\norm{\fh}_\infty)$) and is independent of the choice of $x_\bx$. Because of the latter independence, we can apply dominated convergence and interchange expectation with the spatial differential to obtain both differentiability of $\bP_t \fh$ and that $\norm{\p_x\bP_t \fh}_\infty <\infty$ with, for example, an upper bound given by  that in \eqref{eqns: final_upper_bound}.

\subsection{Proof of Proposition \ref{prop: measure_generator}} \label{appendix: Sec3_Generator}
Here we fill in the details of the heuristic proof given in the main text. First for $h>0$ small enough, $\bP_h G_\fh$ is well-defined, Therefore, for such $h$ and by Proposition \ref{prop: semi-martingality_of_dbs} we have:
\begin{equation*}
\begin{split}
    &\bP_h G_{\fh} (\delta_y) - G_{\fh} (\delta_y)  = \E [E_h M_h G_\fh(\mu_h) - E_0M_0G_\fh (\mu_0) \mid \mu_0 = \delta_y]\\
    & = \E\bigg[\int_0^h E_sM_s H(\mu_s)ds \\
    & \qquad {}+ \int_0^hE_{s-}M_{s-}\int_{\X\times J} G_\fh (\mu_{s-} - \delta_\bx)\bigg[\fm(e, \bx) G_{\fh}(\Phi(\delta_\bx; \bx, e)) - \fh(\bx)\bigg]\tilde N(ds, d\bx, de)\bigg], 
\end{split}
\end{equation*}
where
\begin{equation}\label{eqns: definition of H}
\begin{split}
    H(\mu_s)  = \sum_{\bx \in \mu_s} G_\fh (\mu_s - \delta_\bx) &\bigg[-z_\bx \p_x \fh (x_\bx, z_\bx, \theta_\bx) +(\fu(\bx)-\lambda)\fh(\bx)\\
    & \qquad {}+ \lambda\sum_{e\in J} p(e) \fm(e, \bx) G_\fh (\Phi(\delta_\bx; \bx, e))\bigg]
\end{split}\end{equation}
Note that by Corollary \ref{corollary: martingale of product of dbs}, the stochastic integral term is an $L_1$ martingale, if 
\begin{align}\label{eqns: sufficient_martingale_integrability}
    \E\bigg[\int_0^h \sum_{\bx \in \mu_t} \sum_{e\in J}  p(e)\abs{E_t M_t}\abs{G_\fh (\mu_{t-} - \delta_\bx)[\fm(e, \bx) G_\fh (\Phi(\delta_\bx; \bx, e)) - \fh(\bx)]} dt\bigg]<\infty.
\end{align}
Note also that 
\begin{align}
    &\sum_{\bx \in \mu_t} \sum_{e\in J}p(e)\abs{ G_\fh (\mu_{t-} - \delta_\bx)[\fm(e, \bx) G_\fh (\Phi(\delta_\bx; \bx, e)) - \fh(\bx)]} \\
    & \leq  \abs{\mu_t} \norm{\fh}_\infty ^{\abs{\mu_t}-1} (m_f \norm{\fh}_\infty +\sum_{e\in J\backslash\{\mathsf{f}\}} m_e \norm{\fh}_\infty ^3 + \norm{\fh}_\infty ).
\end{align}
Therefore the integral in \eqref{eqns: sufficient_martingale_integrability} is upper bounded by:
\begin{align}
    \int_0^h \p_B V(t, \norm{\fh}_\infty ) dt  (m_f \norm{\fh}_\infty +\sum_{e\in J\backslash\{\mathsf{f}\}} m_e \norm{\fh}_\infty ^3+ \norm{\fh}_\infty ),
\end{align}
which is finite for small enough $h$.  Therefore for such small $h>0$, we have by martingality,
\begin{align}\label{eqns: difference quotient and drift}
    \bP_h G_{\fh} (\delta_y) - G_{\fh} (\delta_y)  =  \E\bigg[\int_0^h E_s M_s H(\mu_s)ds \mid \mu_0 = \delta_y\bigg].
\end{align}
In addition, $\lim_{h\to 0}\frac{1}{h}\int_0^h E_s M_s H(\mu_s)ds  = H(\mu_0) $ almost surely. By standard argument, for instance, by showing that $Y_h:= \frac{1}{h}\int_0^h E_s M_s H(\mu_s)ds$ is uniformly bounded in $L_p$ for some $p>1$ using the explicit expression of $H$ in \eqref{eqns: definition of H} and $V$ in Proposition \ref{prop: property of V} for every sufficiently small $h>0$, then we have:
\begin{align}
    \bL G_\fh (\delta_\bx) &= \lim_{h\to 0}\frac{1}{h}[\bP_h G_{\fh} (\delta_\bx) - G_{\fh} (\delta_\bx)] = H(\delta _\bx)\\
    &=-z_\bx \p_x \fh (\bx) + (\fu(\bx) - \lambda)\fh(\bx) + \lambda\sum_{e\in J}p(e) \fm(e, \bx) G_\fh (\Phi(\delta_\bx; \bx, e)).
\end{align}
Lastly, boundedness of the mapping $\bL \fh(\bx):= \bL G_\fh(\delta_\bx)$ is clear from the explicit formula; for instance we can have:
\begin{align}
    \norm{\bL \fh}_\infty \leq \norm{\p_x \fh}_\infty + \norm{\fu - \lambda}_\infty \norm{\fh}_\infty  + \lambda_f m_f \norm{\fh}_\infty + \lambda_b \bar m_b \norm{\fh}_\infty^3,
\end{align}
where $\bar m_b$ is defined in \eqref{eqns: definition_of_t_star_sim}. 
\subsection{Proof of Lemma \ref{lemma: non_linear_necessary_condition}}\label{appendix: Sec3_Necessary_Condition}
Fix $t\in [0, T)$. We define:
\begin{align}
    \psi^t(s, \bx):= \tilde \psi(t - s, \bx), \quad s\in [0, t],\quad \bx\in \X.
\end{align}
Then $\psi^t$  satisfies:

\begin{equation}\label{eqns: reverse_time_Soler}
\begin{cases}
   -\hfill  \p_s  \psi^t(s, \bx) &= -z_\bx\p_x  \psi^t(s, \bx)\\
   &\quad -im\theta_\bx  G_{  \psi^t(s, \cdot )}(\Phi(\delta_\bx, \mathsf{f}))+i\theta_\bx G_{\psi^t(s, \cdot )}(\Phi(\delta_\bx, \mathsf{b_1}))+i\theta_\bx G_{\psi^t(s, \cdot )}(\Phi(\delta_\bx, \mathsf{b_2}));\\
    \hfill  \psi^t(t, \cdot ) &= \bg(\bx), \quad (s, \bx)\in [0, t]\times \X. 
\end{cases}
\end{equation}
Next, we prove the time-augmented version of the It\^o formulae in Proposition \ref{prop: semi-martingality_of_dbs}. To this end consider the time-augmented local characteristic over $[0, \infty)\times E^\lift$ defined by:
\begin{align}\label{eqns: time_augmented local characteristics}
\tilde \phi(t; (s, \widehat \mu )) &:= (s+t, \phi^\lift(t, \widehat \mu)) \\
\tilde \ell(t; (s, \widehat \mu)) &:= \ell^{\lift} (t; \widehat \mu)\\
     \tilde \Gamma((t, \widehat\mu); d(t', \widehat\mu')) &=  \Gamma^\lift (\widehat\mu; d\widehat\mu')\delta_t(dt'),
\end{align}
for every $\widehat \mu\in E^\lift$. Fix a Dirac branching process $(\mu_t)$. It is straightforward to check that if $(\tau_n, \rho_n, \bx_n, e_n)$ is the embedded Markov chain of $(\mu_t)$, then:
\begin{align}
    (\tau_n, \tau_n, \rho_n, \bx_n, e_n)
\end{align}
is a Markov chain with transition kernel given by the local characteristic defined by \eqref{eqns: time_augmented local characteristics}. Therefore $(t, \mu_t, \xi_t)$ where $(\mu_t, \xi_t)$ is the lifted PDMP associated to $(\mu_t)$, is a PDMP associated to \eqref{eqns: time_augmented local characteristics} with natural filtration given by $\bF^\lift$. 
Next, we extend the domain of $\psi^t$ to $[0, \infty)$ by defining:
\begin{align}
    \Psi^t(s, \bx):= 
\begin{cases}
    \tilde \psi(t - s, \bx), \quad s\in [0, t]\\
    g(\bx), \quad s\geq t,
\end{cases},\quad \bx\in \X.
\end{align}
Then we define the test function $H^t:[0, \infty)\times E\to \C$ by:
\begin{align}
    H^t(s, \mu):= G_{\Psi^t(s,\cdot)}(\mu), \quad s\in [0, \infty).
\end{align}
 Next, we note that for every $s\geq 0, h> 0$, we have:
\begin{align}
  H^t\circ \tilde \phi(h; (s, \mu))=  H^t(h+s, \phi(h;\mu)) = 
\begin{cases}
      \displaystyle\prod_{\bx\in \mu} \psi^t(h+s, x_\bx -z_\bx h, z_\bx, \theta_\bx), \quad s+h \in [0, t]\\
       \displaystyle\prod_{\bx\in \mu} g(x_\bx -z_\bx h, z_\bx, \theta_\bx), \quad s+h\geq t
\end{cases}.
\end{align}
For every $(s, \mu)\in [0, \infty)\times E$, if  $H^t \circ \tilde \phi$ were piece-wise differentiable then we would have:
\begin{equation}\label{eqns: derivative_of_test_function_appendix}
\begin{split}
   & \frac{d}{dh} H^t\circ \tilde \phi(h; (s, \mu))\\
   &= 
\begin{cases}
   \displaystyle \sum_{\bx\in \mu} (\p_s \psi^t - z_\bx \p_x\psi^t)(h+s, x_\bx - z_\bx h, z_\bx, \theta_\bx) \prod_{\by \in \mu - \delta_{\bx}} \psi^t(h+s, x_\by - z_\by h, z_\by, \theta_\by), \quad s+h \in (0, t)\\
    \displaystyle \sum_{\bx\in \mu} - z_\bx \p_x g( x_\bx - z_\bx h, z_\bx, \theta_\bx) \prod_{\by \in \mu - \delta_{\bx}} g(x_\by - z_\by h, z_\by, \theta_\by), \quad s+h>t
\end{cases}.
\end{split}
\end{equation}
In fact, $H^t \circ \tilde \phi$ is  Lipschitz in $h$, so we only have almost-sure differentiability in $h$ and \eqref{eqns: derivative_of_test_function_appendix} should be understood in its integral form.  Since $\psi^t$ satisfies \eqref{eqns: reverse_time_Soler}, we have $\p_s \psi^t - z_\bx \p_x \psi^t$ to be a polynomial in $\psi^t = \psi(t - \cdot)$. By the $L_\infty$ assumption of the original solution $\psi$ at \eqref{eqns: necessary_boundedness_assumption_nonlinear} and the fact that $g\in C_b^1$, the function $H^t\circ \tilde \phi$ is of local finite variation, for every fixed state $(s, \mu)\in [0, \infty)\times E$. Hence, by applying Proposition \ref{prop: semi-martingale_of_PDMP} on the time-augmented PDMP $(t, \mu_t, \xi_t)$ which has the same natural filtration $\bF^\lift$ considered in Proposition \ref{prop: semi-martingality_of_dbs},  we obtain for $s\in [0, t)$ that:
\begin{align}\label{eqns: time_augmented_Ito_appendix}
    \frac{d(E_s M_s H^t(s, \mu_s))}{E_{s-}M_{s-}} &= \sum_{\bx\in \mu_s} H^t(s, \mu_s - \delta_\bx) \bigg[(\p_s \psi^t - z_\bx \p_x \psi^t)(s, \bx) + \lambda \sum_{e\in J} p(e) \fm(e, \bx) H^t(s, \Phi(\delta_\bx; \bx, e))\bigg] ds\\
    &\quad + \int_{\X\times J} H^t (s,\mu_{s-}-\delta_\bx)\bigg[\fm(e, \bx) H^t(s, \Phi(\delta_\bx; \bx, e)) -\psi^t(s, \bx)\bigg] \tilde N(ds, d\bx, de),\quad s>0\\
    E_0 M_0 H^t(0, \mu_0 = \delta_{\bx}) &= \tilde \psi(t, \bx).
\end{align}
By \eqref{eqns: reverse_time_Soler}, the drift term in \eqref{eqns: time_augmented_Ito_appendix} vanishes. Meanwhile, for every $s\in (0, t)$, we have:
\begin{align}\label{eqns: appendix - necessary proof_expectation}
    &\E[\sum_{\bx \in \mu_s} \sum_{e\in J} p(e) \abs{E_s M_s H^t (s,\mu_{s-} - \delta_\bx)[ \fm(e, \bx) H^t (s, \Phi(\delta_\bx; \bx, e)) - \psi^t(s, \bx)]] }]\\
    \leq & C \p_B V(s, B; \norm{\fu}_\infty, (m_e)_{e\in J}), \quad B:= \sup_{t\in [0, T)} \norm{\tilde \psi(t, \cdot)}_\infty ,
\end{align}
where $C$ is some constant depending on $B$ and $(m_e)_{e\in J}$ which is finite, and $V$ is defined in \eqref{eqns: definition_of_V}. By the choice of upper bound in \eqref{eqns: necessary_boundedness_assumption_nonlinear}, the expectation in  \eqref{eqns: appendix - necessary proof_expectation} is finite with:
\begin{align}
    \int_0^t \p_B V(s, B; \norm{\fu}_\infty, (m_e)_{e\in J}) ds<\infty.
\end{align}
Hence, the local martingale term in \eqref{eqns: time_augmented_Ito_appendix} is a true martingale over $[0, t)$ to the filtration $\bF^\lift$. 
It follows that for every $t>r>s\geq 0$, we have:

\begin{align}
    \E[E_r M_r H^t(r, \mu_r) -E_s M_s H^t(s, \mu_s) \mid \F_s] = 0.
\end{align}
Taking $s = 0$ gives:
\begin{align}
    \tilde \psi (t, \bx) =  \E[E_r M_r H^t(r, \mu_r)\mid \mu_0 = \delta_\bx],
\end{align}
for every $0<r<t$. By the \eqref{eqns: necessary_boundedness_assumption_nonlinear} again, we can send $r\to t^-$ and apply dominated convergence to have:
\begin{align}
    \tilde \psi (t, \bx) =  \E[E_t M_{t-} H^t(t, \mu_{t-})\mid \mu_0 = \delta_\bx],
\end{align}
Using the fact that the pure birth process $(N_t)$ defined in \eqref{eqns: pure birth process} satisfies $\P(N_s - N_{s-} = 0) = 1$ for all $s\geq 0$, we have:
\begin{align}\label{eqns: left_limit_sam_as_original_time}
    \tilde \psi(t, \bx) = \E[E_t M_t G_g(\mu_t)\mid \mu_0 = \delta_\bx]
\end{align}
The proof is done since $t<T$ is arbitrary. 
\section{Proof of Lemma \ref{lemma: CliffordPair_Change_of_coordinates}}\label{appendix: Sec4_CliffordPair}
The strategy is to show that every matrix pair $(\gamma_0, \gamma_1)$ satisfying \eqref{eqns: anti_commute_relation} is similar to $(\sigma_3, i\sigma_1)$. By the anti-commutativity condition \eqref{eqns: anti_commute_relation}, we have $\gamma_0^2 -I = 0$. From elementary linear algebra, it follows that $\gamma_0$ is diagonalizable with eigenvalues from $\pm 1$. Since $\gamma_0$ anti-commutes with the invertible matrix $\gamma_1$ (invertible because $-\gamma_1^2 = Id)$,  it must be the case that $\gamma_0$ has distinct eigenvalues and there exists an invertible $S$ such that:
\begin{align} \label{eqns: appendix_CliffordS}
    \gamma_0 = S \sigma_3
S^{-1}, \quad \sigma_3: = \begin{bmatrix}
    1 & 0 \\
    0 & -1
\end{bmatrix}
\end{align}
Meanwhile, since $S^{-1}\gamma_1 S$ anti-commutes with the invertible $\sigma_3$, one can easily check that we must have:
\begin{align}\label{eqns: appendix_Cliffordbc}
    S^{-1}\gamma_1 S = 
    \begin{bmatrix}
        0 & b \\
        c & 0
    \end{bmatrix}, \quad \mbox{for some}\quad b, c\in \C\backslash\{0\}. 
\end{align}
Since $(S^{-1}\gamma_1 S)^2 = -I$ we have $bc = -1$. It follows from explicit calculation that the matrix $D:= 
\begin{bmatrix}
    1 & 0 \\
    0 & i/b
\end{bmatrix}$ is the transformation matrix between $S^{-1}\gamma_1 S$ and $i\sigma_1$, that is:
\begin{align}
    S^{-1}\gamma_1 S = D(i\sigma_1)D^{-1}. 
\end{align}
Since $D$ is diagonal, it commutes with $\sigma_3$, which is also diagonal. It follows that:
\begin{align}
    \gamma_0 = SD\sigma_3 (SD)^{-1}, \quad  \gamma_1 = SD(i\sigma_1) (SD)^{-1}.
\end{align}
In the case where $(\gamma_0, \gamma_1)$ are either Hermitian or anti-Hermitian. By relation \eqref{eqns: anti_commute_relation},  $\gamma_0$ has to be Hermitian and $\gamma_1$ has to be anti-Hermitian by looking at their possible eigenvalues. It follows that the matrix $S$ in \eqref{eqns: appendix_CliffordS} can be chosen to be unitary. It follows that $b, c$ in \eqref{eqns: appendix_Cliffordbc} have to satisfy $\abs{b} = \abs{c} = 1$ and so the resulting $D$ is unitary. This concludes the proof. 

\section{Proof of Lemma \ref{lemma: apriori_L_inf_estimate} - the a priori \texorpdfstring{$L_\infty$}{L infinity} estimate}\label{sec: proof of Lemma_L_inf_estimate}
We follow closely the calculations in Section 2 of \cite{Huh2013GrossNeveu} with slight modifications. 

Define $\bpsi(t, x, z, \theta):= \psi(t, x, z)\I_{\theta = 1}+\overline{\psi(t, x, z)}\I_{\theta = -1}$. Using the product rule on \eqref{eqns: sec4_nonlinear_Dirac}, we can obtain for every $t\geq 0, x\in \R, z, \theta \in \{\pm 1\}$:\ the very useful identity: 
\begin{align}\label{eqns: core_identity}
    \p_t \abs{\psi(t, x, z)}^2 + z\p_x \abs{\psi(t, x, z)}^2 = \sum_{\theta = \pm 1}D(t, x, z, \theta),
\end{align}
where,
\begin{align}\label{eqns: definition_of_D}
    D(t, x, z, \theta) &= i\theta \bpsi(t, x, z, -\theta) \bpsi(t, x, -z, \theta)\\
    &\quad \times [-m + \sum_{z' = \pm 1}\bpsi(t, x, z', 1) \bpsi(t, x, -z', -1)].
\end{align}
The a priori $L_\infty$ bound \eqref{eqns: apriori_L_inf_bound} follows from two manipulations on \eqref{eqns: core_identity}. First, fix $x\in \R$, we define for $z = \pm 1$ the curve $\gamma_z(t) := x+zt$. If $\abs{\psi(t, \gamma_z(t), z)}^2$ were differentiable in $t$, then we have:
\begin{align}\label{eqns: hold_as_1}
    \frac{d}{dt} \abs{\psi(t, \gamma_z(t), z)}^2 = \sum_{\theta = \pm 1} D(t, \gamma_z(t), z, \theta), \quad t\in [0, T).
\end{align}
Define
\begin{align}
    q_z(t) :=  \abs{\psi(t, \gamma_z(t), z)}^2, \quad r_z(t):= \abs{\psi(t, \gamma_z(t), -z)}^2
\end{align} and it follows from the definition of $D$ \eqref{eqns: definition_of_D} that:
\begin{align}\label{eqns: hold_as_2}
    \frac{d}{dt} q_z(t) &\leq 2\sqrt{q_z(t)r_{z}(t)}m +4q_z(t)r_{z}(t)\\
    &\leq mq_z(t) + mr_{z}(t) + 4q_z(t)r_{z}(t)\\
    &=q_z(t) [m+ 4r_{z}(t)] + mr_{z}(t),
\end{align}
in which the second inequality follows from the AM-GM inequality. Therefore we have: 
\begin{align}\label{eqns: hold_as_3}
    \frac{d}{dt} q_z(t) -q_z(t) [m+ 4r_{z}(t)] \leq mr_{z}(t),
\end{align}
which implies, by the method of integrating factor,
\begin{align}\label{eqns: hold_as_4}
    \frac{d}{dt}\bigg[q_z(t) \exp(-\int_0^t m+ 4r_{z}(s)ds)\bigg]&\leq mr_{z}(t)\exp(-\int_0^tm+ 4r_{z}(s)ds)\leq mr_{z}(t), \quad t\in [0, T).
\end{align}
By monotonicity of integration, we have:
\begin{align}\label{eqns: hold_true_1}
    q_z(t) \exp(-\int_0^t m+ 4r_{z}(s)ds) - q_z(0) \leq m\int_0^t r_{z}(s)ds, \quad t\in [0, T)
\end{align}
Indeed $\abs{\psi(t, \gamma_z(t), z)}^2$ is only Lipschitz in $t$, so \eqref{eqns: hold_as_1},  \eqref{eqns: hold_as_2}, \eqref{eqns: hold_as_3}, \eqref{eqns: hold_as_4} should be understood to hold almost surely on $[0, T)$ while we still have \eqref{eqns: hold_true_1} to hold for all $t\in [0, T)$. It follows that:
\begin{align}\label{eqns: first_apriori_estimate}
    q_z(t) \leq [\abs{g(x, z)}^2 +m\int_0^t r_{z}(s)ds]\exp(\int_0^tm+ 4r_{z}(s)ds), \quad t\in [0, T).
\end{align}
The second manipulation follows from the symmetry of $D$ that for every $t>0,x\in \R, z = \pm 1$ and $\theta = \pm 1$
\begin{align}
    D(t, x, -z, -\theta) = -D(t, x, z, \theta).
\end{align}
Hence summing \eqref{eqns: core_identity} over $z= \pm 1$ gives:
\begin{align}\label{eqns: local_conservation_identity}
      \sum_{z = \pm 1}\p_t \abs{\psi(t,x, z)}^2 + z\p_x  \abs{\psi(t,x, z)}^2  = 0. 
\end{align}
We then define  $u:[0, T)\times \R\to \R^2$ by:
\begin{align}
    u(t, x) := (a(t, x), b(t, x))
\end{align}
where:
\begin{align}
    a(t, x):= \sum_{z = \pm 1}\abs{\psi(t, x, z)}^2 , \quad b(t, x):= \sum_{z = \pm 1} z\abs{\psi(t, x, z)}^2, 
\end{align}
so that the left-hand-side of \eqref{eqns: local_conservation_identity} is the divergence of $u$. We remark that $u$ is both differentiable and Lipschitz, thanks to the condition \eqref{eqns: classical_solution_additional_boundedness_condition}. 
Define for every $t>0, x\in \R$ the open domain:
\begin{align}
    \Omega_{t, x}^+:= \{(s, y)\in \R^2: 0<s<t, \quad x+s<y<x+2t-s\}, 
\end{align}
which is the triangle with vertices $(0, x), (t, x+t)$ and $(0, x+2t)$, together with:
\begin{align}
     \Omega_{t, x}^-:= \{(s, y)\in \R^2: 0<s<t, \quad x-2t+s<y<x-s\}, 
\end{align}
which is the triangle with vertices $(0, x), (t, x-t)$ and $(0, x-2t)$. It is straightforward to see that both domains have Lipschitz boundary (see for instance \cite[Definition 9.56]{Leoni2017}) and so by  the Divergence theorem for Lipschitz functions, see for instance \cite[Theorem 9.64]{Leoni2017}, we have:
\begin{align}\label{eqns: divergence theorem_appendix}
    \int_{\p \Omega_{t, x}^+} u\cdot \nu d\mathcal{H}^1=   \int_{\p \Omega_{t, x}^-}u\cdot \nu d\mathcal{H}^1  = 0,
\end{align}
where $\nu$ is the outward unit normal to the boundary of the domain under consideration and $\mathcal{H}^1$ is the one-dimensional Hausdorff measure of $\R^2$. 

In particular, evaluating the integral over $\p\Omega_{t, x}^+$ in \eqref{eqns: divergence theorem_appendix} gives:
\begin{align}
    \int_0^{2t} a(0, x+s)ds &= \int_0^t a(s, x+s) - b(s, x+s)ds + \int_0^t a(t-s, x+t+s)+b(t-s, x+t+s)ds\\
    &=2\int_0^t \abs{\psi(s, x+s, -1)}^2ds +2\int_0^t \abs{\psi(t-s, x+t+s, 1)}^2ds\\
    &\geq 2\int_0^t r_{1}(s) ds.
\end{align}
Hence, we have $\int_0^t r_{1}(s) \leq 2t \norm{g}_\infty^2$ using the fact that $a(0, \cdot) \leq 2\norm{g}_\infty^2$. Similarly by considering the integral over $\p\Omega_{t, x}^-$ in \eqref{eqns: divergence theorem_appendix},  we obtain $\int_0^t r_{-1}(s) \leq 2t \norm{g}_\infty^2$ and so $\int_0^t r_z(s) ds \leq 2t \norm{g}_\infty^2 $ for $z=\pm 1$. Putting this into \eqref{eqns: first_apriori_estimate}, we obtain for all $x\in \R, z= \pm 1$ and $t\in [0,T)$ that:
\begin{align}
    \abs{\psi(t, x+zt, z)}^2=q_z(t) \leq [\norm{g}_\infty^2 + 2mt\norm{g}_\infty^2] \exp(mt + 8t\norm{g}_\infty^2).
\end{align}
Translating $x$ gives the desired inequality \eqref{eqns: apriori_L_inf_bound}. 

\section{Proof details for Sections \ref{subsec: convergence_analysis} - \ref{subsec: global derivative bound and stability}} \label{appendix. Global arguments}
This appendix collects the proof details for  Theorem \ref{thm: convergence_of_global_alg},  and Proposition \ref{prop: global_independent_results}. 

\subsection{Proof of Theorem \ref{thm: convergence_of_global_alg}}\label{appendix: sec5_convergence_Proof}
We pick $\tau>0$ small enough such that:
\begin{align}
    \tau< t^*(L(t, \norm{g}_\infty)),\quad V(\tau, L(t, \norm{g}_\infty)^2, 2\lambda; (m_e^2)_{e\in J})<\infty.
\end{align}
We can also choose $\tau$ such that $K:= [t/\tau]=t/\tau\in \N$ to simplify the presentation of the proof. For any $k\in [K]$, we let $\psi_k:= \psi(k\tau, \cdot)$ for all $k\in [K]$. 

By the flow property, we then have $\psi_k = P_\tau \psi_{k-1}$ for every $k\in [K]$. We also define the constant:
\begin{align}\label{eqns: definition of derivative bound constant}
    C:= \sup_{k\in [K]}\norm{\p_x \psi(\tau k,\cdot)}_\infty,
\end{align}
which is finite by \eqref{eqns: classical_solution_additional_boundedness_condition}. The main observation is that we can formally decompose the error as follows:
\begin{align}\label{eqns: error_decomposition}
    &\hat \psi_k- \psi_k = I_{L,M}[\hat P^{(k)}_{\tau, N, M} \hat \psi_{k-1}] - P_\tau  \psi_{k-1} \\
    =& \bigg(I_{L,M}[\hat P^{(k)}_{\tau, N, M} \hat \psi_{k-1} ]-I_{L,M}[\overrightarrow{P_\tau \hat \psi_{k-1}}]\bigg)+\bigg(I_{L,M}[\overrightarrow{P_\tau \hat \psi_{k-1}}]-I_{L,M}[\overrightarrow{P_\tau \psi_{k-1}}]\bigg)+\bigg(I_{L,M}[\overrightarrow{P_\tau \psi_{k-1}}]- P_\tau  \psi_{k-1}\bigg), 
\end{align}
which gives, from left to right, the Monte Carlo error, the stability error, and the spatial interpolation error. The well-definedness of $P_\tau \hat \psi_{k-1}$ in \eqref{eqns: error_decomposition} will be addressed in the subsequent proof. To analyze the errors, we define for every $k\in [K]$:
\begin{align}
    I_k:= [-L +(\tau+\Delta x) (k-1) , L-(\tau+\Delta x) (k-1)],
\end{align}
so that $I_k^{(\tau+\Delta x)} =I_{k-1}$. We choose  $M$ large enough (so $\Delta x$ small enough) such that $ [-L+t, L-t]\subset I_K $ (one possible choice is $\Delta x \leq \tau/K$); we also have $I_0 = [-L-(\tau + \Delta x), L+(\tau +\Delta x)]$. Now fix $\epsilon>0$ and take $R>0$ such that:
\begin{align}
    L(t-\tau, \norm{g}_\infty) + R<L(t, \norm{g}_\infty).
\end{align}
We then consider the three error terms one by one from right to left. First for every $k\in [K]$, including $ k = 1$,
\begin{align}
    \norm{I_{L,M}[\overrightarrow{P_\tau \psi_{k-1}}]- P_\tau  \psi_{k-1}}_{\infty, I_k} \leq \norm{\p_x P_\tau  \psi_{k-1}}_{\infty, I_{k}^{(\Delta x)}} \frac{2L}{M-1} \leq \frac{2LC}{M-1},
\end{align}
where the mean value theorem is used in the first inequality and $C$ is defined in \eqref{eqns: definition of derivative bound constant}. 
Next, we consider the stability error which vanishes for the case $k = 1$:
\begin{align}
    &\P\Big(\norm{I_{L,M}[\overrightarrow{P_\tau \hat \psi_{k-1}}]-I_{L,M}[\overrightarrow{P_\tau \psi_{k-1}}]}_{\infty, I_k} \geq \epsilon \Big)\\
    \leq& \P\Big(\norm{I_{L,M}[\overrightarrow{P_\tau \hat \psi_{k-1}}]-I_{L,M}[\overrightarrow{P_\tau \psi_{k-1}}]}_{\infty, I_{k}} \geq \epsilon, \norm{\hat \psi_{k-1} - \psi_{k-1}}_{\infty, I_{k-1}}< R\Big) + \P\Big( \norm{\hat \psi_{k-1} - \psi_{k-1}}_{\infty, I_{k-1}} \geq R\Big).
\end{align} 
The choice of $R$ allows us to apply the stability result \eqref{eqns: sec6. stability_local} in Lemma \ref{lemma: a local derivative bound and stability}, giving:
\begin{align}
    &\P\Big(\norm{I_{L,M}[\overrightarrow{P_\tau \hat \psi_{k-1}}]-I_{L,M}[\overrightarrow{P_\tau \psi_{k-1}}]}_{\infty, I_k} \geq \epsilon, \norm{\hat \psi_{k-1} - \psi_{k-1}}_{\infty, I_{k-1}}< R\Big) \\
    \leq &\P\Big(\norm{P_\tau \hat \psi_{k-1}- P_\tau \psi_{k-1}}_{\infty, I_k^{(\Delta x)}} \geq  \epsilon, \norm{\hat \psi_{k-1}-\psi_{k-1}}_{\infty, I_{k-1}}<R \Big)\\
    \leq & \P\Big(\norm{\hat \psi_{k-1} - \psi_{k-1}}_{\infty, I_{k-1}} \geq \frac{\epsilon}{\p_B V(\tau, L(t, \norm{g}_\infty))} \Big).  
\end{align}
Lastly, for the Monte Carlo error, we have the conditional expectation:
\begin{align}
    &\P\Big(\norm{I_{L,M}[\hat P^{(k)}_{\tau, N, M} \hat \psi_{k-1} ]-I_{L,M}[\overrightarrow{P_\tau \hat \psi_{k-1}}]}_{\infty, I_k} \geq \epsilon \mid \hat \psi_{k-1}\Big)\\
    \leq &\P\Big(\norm{I_{L,M}[\hat P^{(k)}_{\tau, N, M} \hat \psi_{k-1} ]-I_{L,M}[\overrightarrow{P_\tau \hat \psi_{k-1}}]}_{\infty, I_k} \geq \epsilon , \norm{\hat \psi_{k-1} - \psi_{k-1}}_{\infty, I_{k-1}}< R \mid \hat \psi_{k-1}\Big) \\
    & \quad + \P\Big(\norm{\hat \psi_{k-1} - \psi_{k-1}}_{\infty, I_{k-1}} \geq R \mid \hat \psi_{k-1}\Big).
\end{align}
By union bound,
\begin{align}
    &\P\Big(\norm{I_{L,M}[\hat P^{(k)}_{\tau, N, M} \hat \psi_{k-1} ]-I_{L,M}[\overrightarrow{P_\tau \hat \psi_{k-1}}]}_{\infty, I_k} \geq \epsilon , \norm{\hat \psi_{k-1} - \psi_{k-1}}_{\infty, I_{k-1}}< R \mid \hat \psi_{k-1}\Big) \\
    \leq & \P\Big(\sup_{x_j\in I_k^{(\Delta x)}, z = \pm 1}\abs{\hat P^{(k)}_{\tau, N, M} \hat \psi_{k-1}(x_j, z) - \overrightarrow{P_\tau \hat \psi_{k-1}}(x_j, z)}\geq \epsilon, \norm{\hat \psi_{k-1}}_{\infty, I_{k-1}} <L(t, \norm{g}_\infty) \mid \hat \psi_{k-1}\Big)\\
    \leq & \sum_{x_j\in I_k^{(\Delta x)}, z = \pm 1} \P\Big(\abs{\hat P^{(k)}_{\tau, N, M} \hat \psi_{k-1}(x_j, z) - \overrightarrow{P_\tau \hat \psi_{k-1}}(x_j, z)} \geq \frac{\epsilon}{2M},  \norm{\hat \psi_{k-1}}_{\infty, I_{k-1}} <L(t, \norm{g}_\infty) \mid \hat \psi_{k-1}\Big)\\
    \leq & \frac{8M^3}{\epsilon^2 N} V(\tau, L(t, \norm{g}_\infty)^2, 2\lambda; (m_e^2)_{e\in J}). 
\end{align}
Hence, taking the expectation over $\hat \psi_{k-1}$ gives:
\begin{align}
  &  \P(\norm{I_{L,M}[\hat P^{(k)}_{\tau, N, M} \hat \psi_{k-1} ]-I_{L,M}[\overrightarrow{P_\tau \hat \psi_{k-1}}]}_{\infty, I_{k}} \geq \epsilon)\\
  \leq& \frac{8M^3}{\epsilon^2 N} V(\tau, L(t, \norm{g}_\infty)^2, 2\lambda; (m_e^2)_{e\in J}) + \P(\norm{\hat \psi_{k-1} - \psi_{k-1}}_{\infty, I_{k-1}} \geq R). 
\end{align}
To sum up, for every $k\in [K]$, we have:
\begin{align}
    \P(\norm{\hat \psi_k - \psi_k}_{\infty, I_k} \geq 3\epsilon ) &\leq \I_{\{\p_B V(\tau, L(t-\tau, \norm{g}_\infty))\frac{2LC}{M-1}\geq \epsilon\}} + 2\P\Big(\norm{\hat \psi_{k-1} - \psi_{k-1}}_{\infty, I_{k-1}}\geq R\Big)\\
    &\quad + \P\Big(\norm{\hat \psi_{k-1} - \psi_{k-1}}_{\infty,I_{k-1} }\geq \frac{\epsilon}{\p_B V(\tau, L(t, \norm{g}_\infty))} \Big)\\
    &\quad + \frac{8M^3}{\epsilon^2 N} V(\tau, L(t, \norm{g}_\infty)^2, 2\lambda; (m_e^2)_{e\in J}).
\end{align}
It follows clearly that the result follows from an inductive argument and by first sending $N\to \infty$ and $M\to \infty$.

\subsection{Proof of Proposition \ref{prop: global_independent_results}}\label{appendix: sec5_global_stability_proof}
   Let $R>0$ and $g\in C_b^1$. Let $h_1,  h_2\in C_b^1$ be such that $\norm{h_1 - g}_\infty , \norm{h_2-g}_\infty\leq R$. We consider a time-marching parameter of:
\begin{align}
    \tau < t^*(L(t, \norm{g}_\infty + R)).
\end{align}
Then by Theorem \ref{thm: global_well_posednesss_L_inf}, we have:
\begin{align}
    \psi_{h_j}(s, \cdot) = P_{s - [s/\tau]\tau}P_\tau^{[s/\tau]}h_j, \quad s\in [0, t], \quad j = 1, 2.
\end{align}
It follows from Lemma \ref{lemma: a local derivative bound and stability} that:
\begin{align}\label{eqns: appendix_stability_inequality}
  &  \norm{\psi_{h_1}(s, \cdot) - \psi_{h_2}(s, \cdot)}_\infty  \\
  \leq & \p_B V(\tau, \max \{\norm{\psi_{h_1}([s/\tau]\tau)}_\infty, \norm{\psi_{h2}([s/\tau]\tau)}_\infty \} ) \norm{\psi_{h_1}([s/\tau]\tau) - \psi_{h_2}([s/\tau]\tau)}_\infty.
\end{align}
Since
\begin{align}
    \max \{\norm{\psi_{h_1}([s/\tau]\tau)}_\infty, \norm{\psi_{h2}([s/\tau]\tau)}_\infty \} \leq\max \{L(t, \norm{h_1}_\infty), L(t, \norm{h_2}_\infty)\} \leq L(t, \norm{g}_\infty+ R),
\end{align}
repeatedly applying \eqref{eqns: appendix_stability_inequality} to itself gives:
\begin{align}
    \norm{\psi_{h_1}(s, \cdot) - \psi_{h_2}(s, \cdot)}_\infty  \leq \p_B V(\tau, L(t, \norm{g}_\infty +R))^{[s/\tau]+1} \norm{h_1 - h_2}_\infty,
\end{align}
for every $s\in [0, t]$. Taking $s = t$ yields the result. 
\bibliographystyle{apalike}
\bibliography{bibtex}

\end{document}